\documentclass[preprint]{elsarticle}

\usepackage[english]{babel}

\usepackage[scale=0.8,a4paper]{geometry}

\usepackage{mathtools}
\usepackage{amsthm}
\usepackage{thmtools}
\usepackage{amssymb}
\usepackage{graphicx}
\usepackage{subcaption}
\usepackage{booktabs}
\usepackage{multirow}
\usepackage{makecell}
\usepackage{proof}
\usepackage{bm}
\usepackage{enumitem}
\usepackage{siunitx}
\usepackage{algorithm}
\usepackage{algpseudocode}
\usepackage{indentfirst}

\usepackage{pgfplots}
\pgfplotsset{compat=1.18}
\usetikzlibrary{arrows.meta}
\usepgfplotslibrary{fillbetween}

\definecolor{graph_1}{RGB}{117,112,179}
\definecolor{graph_2}{RGB}{217,95,2}
\definecolor{graph_3}{RGB}{27,158,119}
\definecolor{graph_4}{RGB}{231,41,138}
\definecolor{review1}{RGB}{117,112,179}
\definecolor{review2}{RGB}{217,95,2}
\definecolor{review3}{RGB}{27,158,119}
\definecolor{allreviews}{RGB}{231,41,138}
\newcommand{\ra}{}
\newcommand{\rabis}{}
\newcommand{\rb}{}
\newcommand{\rc}{}
\newcommand{\rall}{}
\newcommand{\raparagraph}{}
\newcommand{\rallparagraph}{}

\usepackage[colorlinks=true, allcolors=blue]{hyperref}
\usepackage[capitalise,nameinlink,noabbrev]{cleveref}
\Crefformat{appendix}{\S#2#1#3}

\newtheorem{theorem}{Theorem}
\newtheorem{corollary}[theorem]{Corollary}
\newtheorem{lemma}[theorem]{Lemma}
\newtheorem{proposition}[theorem]{Proposition}

\newtheorem{remark}[theorem]{Remark}

\numberwithin{equation}{section}

\DeclareMathOperator*{\argmin}{arg\,min}
\DeclareMathOperator*{\polaratan}{atan2}

\newcommand{\dt}{{\Delta t}}
\newcommand{\dx}{{\Delta x}}
\newcommand*{\boldone}{\text{\usefont{U}{bbold}{m}{n}1}}

\pgfplotsset{
	compare SL and NSL/.style={
            axis lines = left,
            enlarge x limits={abs=10pt},
            enlarge y limits={abs=10pt},
            xlabel style={at={(ticklabel* cs:1.01)},anchor=west},
            ylabel style={at={(ticklabel* cs:1.01)},anchor=west},
            label style={font=\small},
            tick label style={font=\footnotesize},
            xlabel = {\makebox[1pt][c]{$d$}},
            hide obscured x ticks=false,
            hide obscured y ticks=false,
            scaled y ticks = false,
            scaled x ticks = false,
            xtick={1, 2, 3, 4, 5, 6, 7, 8},
            grid=major,
            width=0.33\textwidth,
            height=0.3\textwidth,
		}
}

\usepackage{silence}
\begin{document}

\begin{frontmatter}

    \title{Neural semi-Lagrangian method for high-dimensional advection-diffusion problems}

    \author[Inria]{Emmanuel Franck}
    \ead{emmanuel.franck@inria.fr}

    \author[Inria]{Victor Michel-Dansac\corref{cor1}}
    \ead{victor.michel-dansac@inria.fr}
    \cortext[cor1]{corresponding author}

    \author[IRMA,Inria]{Laurent Navoret}
    \ead{laurent.navoret@math.unistra.fr}

    \author[IRMA,Inria]{Vincent Vigon}
    \ead{vincent.vigon@math.unistra.fr}

    \affiliation[Inria]{organization={Université de Strasbourg, CNRS, Inria, IRMA},
        postcode={F-67000},
        city={Strasbourg},
        country={France}
    }

    \affiliation[IRMA]{organization={IRMA, Université de Strasbourg, CNRS UMR 7501},
        addressline={7 rue René Descartes},
        postcode={67084},
        city={Strasbourg},
        country={France}}

    \begin{abstract}
        This work is devoted to the numerical approximation of high-dimensional advection-diffusion equations.
        It is well-known that classical methods, such as the finite volume method,
        suffer from the curse of dimensionality, and that their time step is constrained by a stability condition.
        The semi-Lagrangian method is known to overcome the stability issue, while recent time-discrete neural network-based approaches overcome the curse of dimensionality.
        In this work, we propose a novel neural semi-Lagrangian method that combines these last two approaches.
        It relies on projecting the initial condition onto a finite-dimensional neural space, and then solving an optimization problem, involving the backwards characteristic equation, at each time step.
        It is particularly well-suited for implementation on GPUs, as it is fully parallelizable and does not require a mesh. We provide rough error estimates, present several high-dimensional numerical experiments to assess the performance of our approach, and compare it to other neural methods.
    \end{abstract}




\end{frontmatter}

\tableofcontents

\section{Introduction}\label{sec:introduction}

In this work, we are interested in solving high-dimensional advection-diffusion equations. The large dimensionality may be inherent, for instance in kinetic equations (kinetic transport, Vlasov, or Fokker-Planck equations), where the unknown typically depends on time, space and velocity. It can also be parametric (examples include equations with parametric source terms, differential operators, or boundary conditions).
Traditional numerical methods, such as finite difference, finite volume, or finite element methods, applied to such high-dimensional problems, require numerous degrees of freedom to achieve a given level of accuracy.
Indeed, the number of degrees of freedom
grows as $N^d$, where $N$ is the number of degrees of freedom in each dimension and $d$ is the number of dimensions.
Thus, $N$ grows exponentially with $d$, leading to a curse of dimensionality.
This exponential growth is a major obstacle to the numerical solution of high-dimensional kinetic equations,
since it leads to a large computational cost and large memory requirements.
\rabis{To address such an issue, several methods have been developed, including
model order reduction~\cite{QuaManNeg2016} and low-rank tensor decomposition~\cite{kolda2009tensor}.
These methods express the solution with a reduced number of degrees of freedom.}

\rabis{In recent years, neural networks have also shown their efficiency in reducing the number of degrees of freedom by enriching classical approximation spaces.
For example, tensor decomposition methods have successfully been combined with neural network approximations, see e.g.~\cite{park2024interpolating,GUO2025118101,guo2025interpolating}.
More generally, it turns out that several classical methods have seen their neural counterparts developed.
The foremost instances of this are the Physics Informed Neural Networks (PINNs) \cite{RAISSI2019686,NEURIPS2021_df438e52,De-Ryck:2022aa} and the Deep Ritz method~\cite{EYu2018}.}
In such methods, the solution itself is expressed as a neural network, whose parameters are adjusted such that the equation is satisfied in strong or weak form.  These methods have produced good results for elliptic Partial Differential Equations (PDEs) in low or high physical or parametric dimensions, see e.g.~\cite{SUKUMAR2022114333,9788008,Hu:2024aa}.
There are also a few applications to high-dimensional transport, although these remain more limited \cite{MISHRA2021107705,Zhang_2023,JinAPNN}. Indeed, as will be detailed later, PINNs is a spacetime method, which treat space and time in the same way and specific techniques have to be considered to incorporate causality in the training to properly capture the solution, see~\cite{WANG2024116813}.
\rc{To improve the accuracy of PINNs for transport problems, some work has been undertaken, either incorporating knowledge of the characteristic curves associated with the problem,
    see~\cite{hu2024physics},
    or using a Lagrangian formulation (which is a similar idea), see~\cite{mojgani2022lagrangian}. In these cases, a space-time approximation is used, contrary to the sequential-in-time methods presented below.}

For time-dependent equations, another strategy is to use neural networks approximations as functions of space only, and make the neural network parameters evolve in time. This led to the development of the discrete PINNs \cite{StiCha2023,biesek2023burgers} or Neural Galerkin \cite{LI2022110958,finzi2023stable,BRUNA2024112588,KAST2024112986} methods.
These methods can be seen as extensions of the classical implicit and explicit Galerkin methods, where the finite-dimensional approximation linear subspace has been replaced with a submanifold of neural networks.

When considering advection equations and even advection-diffusion equations in high dimension, however,
semi-Lagrangian approaches \cite{douglas1982numerical,pironneau1982transport,SemiLagrangianIntegrationSchemesforAtmosphericModelsAReview,SONNENDRUCKER1999201,FILBET2001166,XIU2001658,CROUSEILLES20101927,YANG2021110632}  are known to be more efficient than standard Galerkin method. At each time step, by using the backward integration of the characteristic curves, the solution is transported exactly before being projected onto the spatial approximation space. This projection either relies on interpolation methods, like splines, Lagrange or Hermite interpolations \cite{SONNENDRUCKER1999201,CROUSEILLES20101927,besse2003semi}, or is a Galerkin projection \cite{CroMehVec2011,rossmanith2011positivity,qiu2011positivity}. Dimensional splitting can be deployed to go back to one dimensional projection problems as these projections may be particularly computationally expensive in high dimension. A key feature of semi-Lagrangian schemes is that they have no time step stability constraints, even though they are explicit in time. This reduces the computational cost: the time step is then chosen only to increase the accuracy. Recently, neural network-based methods have been proposed to improve the accuracy of \rb{traditional} grid-based semi-Lagrangian methods \cite{LARIOSCARDENAS2022111623,chen2024conservative,Chen2023ALC}.
However, as explained above, numerical simulations for high-dimensional problems are still hard to perform due to the large number of required degrees of freedom \rb{used by the grid-based methods}.
\rabis{To reduce this number, low-rank tensor methods have recently been developed for kinetic equations, see e.g.~\cite{Kor2015,EinJos2021,ZheHayChrQiu2025} or the review~\cite{EinKorKusMcCQiu2024}.}

In order to tackle transport dynamics in large dimension, we propose in this work to combine the semi-Lagrangian approach with neural approximations. The classical projection onto a linear subspace at each time step of the transported solution is then replaced with an optimization process to fit the parameters. This approach constitutes a semi-Lagrangian variant of the discrete PINN method. Indeed, the first step is to project the initial condition onto the neural approximation space; this is nothing but a nonlinear optimization problem.
Then, at each time step, the approximate solution is transported using the characteristic curves, and the parameters are updated by solving a nonlinear optimization problem.

In the following, we will apply this method to the following parametric linear \textit{advection-diffusion} problem, describing the dynamics of an unknown $u(t,x,\mu) \in \mathbb{R}$, depending on time variable $t \in \mathbb{R}_+$, space variable $x \in \Omega \subset \mathbb{R}^d$ and $\mu \in \mathbb{M} \subset \mathbb{R}^p$ a set of physical parameters. The PDE writes:
\begin{equation}
    \label{eq:advection_diffusion}
    \begin{dcases}
        \partial_t u(t,x,\mu)
        + a(t, x,\mu) \cdot \nabla u(t,x,\mu)
        - \sigma \Delta u(t,x,\mu)
        = 0,
         & \quad x\in \Omega, \, t \in (0, T),           \\
        u(t=0,x,\mu)=u_0(x,\mu),
         & \quad x\in \Omega,                            \\
        u(t,x,\mu)=g(t,x,\mu),
         & \quad x \in \partial \Omega, \, t \in (0, T),
    \end{dcases}
\end{equation}
where
$a(t,x,\mu) \in \mathbb{R}^d$ is the advection field and
$\sigma > 0$ is the constant diffusion coefficient.
The parameters may appear in
the boundary or initial conditions,
or in the advection field.
We will formally write the PDE without initial and boundary conditions as
\begin{equation}
    \label{eq:advection_diffusion_time_operator}
    \partial_t u + \mathcal{L} [u] = 0,
\end{equation}
where the operator $\mathcal{L}$ is linear in $u$ and contains no time derivatives of $u$, or as
\begin{equation}
    \label{eq:advection_diffusion_spacetime_operator}
    \mathcal{T} [u] = 0,
\end{equation}
where the operator $\mathcal{T}$ is linear in $u$ and involves both time and space derivatives.
For the remainder of the paper, we assume that the velocity field is known and is not determined by another equation.

The paper is organized as follows.
In \cref{sec:numerical_methods}, we introduce the classical numerical approaches for this problem and their neural extensions.
This discussion illustrates that our approach naturally follows from previous work.
Then, in \cref{sec:NGSL}, we derive the proposed method for pure advection equations, and extend it to advection-diffusion equations.
Algorithms and \rb{rough} error estimates are provided.
Finally, in \cref{sec:applications}, we provide several validation experiments,
for both pure advection and advection-diffusion equations,
in high-dimensional settings.
\cref{sec:conclusion} concludes the paper.

\section{Classical and neural numerical methods}
\label{sec:numerical_methods}

The goal of this section is to highlight
the key differences between classical and neural methods,
and to cast the latter into a formalism
well-known in the case of classical methods.
To that end, \cref{sec:classical_methods} is
devoted to classical methods,
while \cref{sec:neural_methods} tackles neural methods.
We introduce an
infinite-dimensional Hilbert space~$V$
such that $u \in V$,
where $u$ is the solution
to the partial differential equation
\eqref{eq:advection_diffusion}.
In both \cref{sec:classical_methods}
and \cref{sec:neural_methods},
this function space $V$ will be approximated by a
finite-dimensional subspace,
denoted by $V_N^\text{classical}$ and $\mathcal{V}_N^\text{neural}$
for classical and neural methods, respectively.
In the two above-defined spaces, the subscript $N$
represents the number of degrees of freedom.
The main difference between these two finite-dimensional
subspaces is that,
while~$V_N^\text{classical}$ is a linear subspace
spanned by some given basis functions,
$\mathcal{V}_N^\text{neural}$ is no longer linear.
It may become a submanifold of $V$,
depending on the nonlinearity.
As such, projecting an element of $V$
onto $V_N^\text{classical}$
amounts to solving linear systems
(or quadratic optimization problems),
while projecting it onto $\mathcal{V}_N^\text{neural}$
requires solving nonlinear, usually non-convex,
optimization problems.

First, \cref{sec:classical_methods} introduces
three families of classical methods
for advection-diffusion equations.
These approaches all use Galerkin projections onto a
finite-dimensional subspace to represent the PDE solution.
Namely, we present the spacetime Galerkin method
in \cref{sec:spacetime_galerkin_method},
the Galerkin method in \cref{sec:galerkin_method},
and the semi-Lagrangian Galerkin method in \cref{sec:SL_galerkin_method}.
Second, we discuss neural numerical methods in \cref{sec:neural_methods}.
We briefly recall that PINNs (\cref{sec:PINNs})
and their time-discrete counterparts
(\cref{sec:neural_Galerkin})
can be seen as extensions of classical approaches,
where the projection onto a finite-dimensional subspace
(finite element space, spectral basis, etc.)
is replaced by a projection onto the manifold defined by
neural networks with given architectures.

For simplicity, in this section,
we drop the dependence on the parameter $\mu$ in the solution $u$.
We also do not discuss boundary and initial conditions.
Moreover, throughout this section,
$u_N$ denotes the approximate solution of the PDE,
obtained by the numerical method under consideration
in each specific paragraph, with $N$ dofs.

\subsection{Classical methods}
\label{sec:classical_methods}

As mentioned above,
classical approaches to solve time-dependent PDEs
approximate the solution by projecting the equation
onto a finite-dimensional linear subspace
$V_N^\text{classical}$ of the Hilbert space $V$.
To properly introduce this linear subspace,
let $N \in \mathbb{N}$,
and define $N$ basis functions $\phi = \{\phi_1, \ldots, \phi_N\}$.
Examples of basis functions include
piecewise polynomial functions for finite element methods
or Fourier basis functions for spectral methods.
The approximate solution $u_N$ then
belongs to the $N$-dimensional linear subspace
$V_N^\text{classical}$,
defined using the basis functions by
\begin{equation*}
    V_N^\text{classical} = \left\{
    \sum_{j=1}^N \theta_j \phi_j
    = \bigl \langle \theta, \, \phi \bigr \rangle,
    \;\; \theta \in \Theta \subset \mathbb{R}^N
    \right\} = \operatorname{Span}(\phi_1, \ldots, \phi_N),
\end{equation*}
where $\langle \cdot, \cdot \rangle$
denotes the inner product in $\mathbb{R}^N$
and $\theta = (\theta_1, \ldots, \theta_N)$.
The goal is to find the degrees of freedom (dofs)
$\theta \in \Theta \subset \mathbb{R}^N$ that define $u_N$.
These dofs will depend on the choice of basis functions
and on the equation to solve.
Since we are considering linear equations
and projections onto the
linear subspace $V_N^\text{classical}$,
the dofs~$\theta$ will end up being determined
by solving linear systems.

\subsubsection{Spacetime Galerkin method}
\label{sec:spacetime_galerkin_method}

In the spacetime Galerkin method (see e.g.~\cite{quarteroni2008numerical}),
the basis functions depend on both space and time.
The solution is approximated,
for all $t \in (0, T)$ and $x \in \Omega$, by
\begin{equation}
    \label{eq:spacetime_Galerkin_approximation}
    u_N(t, x)
    = \sum_{j=1}^N \theta_j \phi_j(t, x)
    = \bigl \langle \theta, \, \phi(t, x) \bigr \rangle.
\end{equation}
Now, to derive the equation for the dofs $\theta$,
we plug the spacetime Galerkin approximation
\eqref{eq:spacetime_Galerkin_approximation}
into the PDE \eqref{eq:advection_diffusion_spacetime_operator},
and integrate over the spacetime slab $(0, T) \times \Omega$
against some test function.
We obtain, for all~$i \in \{1, \ldots, N\}$,
\begin{equation*}
    \int_\Omega \int_0^T \mathcal{T}[u_N](t,x) \phi_i(t,x) \, dt \, dx = 0.
\end{equation*}
Using \eqref{eq:spacetime_Galerkin_approximation}, we obtain
\begin{equation*}
    \int_\Omega \int_0^T \sum_{j=1}^N
    \theta_j \mathcal{T}[\phi_j](t,x) \phi_i(t,x)
    \, dt \, dx = 0,
\end{equation*}
and so
\begin{equation*}
    \sum_{j=1}^N \theta_j
    \left( \int_\Omega \int_0^T
    \phi_i(t,x) \mathcal{T}[\phi_j](t,x)
    \, dt \, dx \right) = 0.
\end{equation*}
This leads to the linear system $M \theta = 0$,
where the components of the mass matrix $M$ are defined by
\begin{equation*}
    M_{ij} = \int_\Omega \int_0^T \phi_i(t,x) \mathcal{T}[\phi_j](t,x) \, dt \, dx.
\end{equation*}
The solution of this linear system then gives the dofs $\theta$.

\subsubsection{Galerkin method}
\label{sec:galerkin_method}

Now, assume that the basis functions only depend on space,
and that the dofs $\theta$ depend on time.
The idea is to approximate the solution
at each discrete time using these space-dependent basis functions.
This is equivalent to stating that the linear combination of
the basis functions used to represent the solution depends on time.
Thus, the approximate solution $u_N \in V$
at time $t \in (0, T)$ and position $x \in \Omega$ is given by
\begin{equation}
    \label{eq:Galerkin_approximation}
    u_N(t,x)
    = \sum_{j=1}^N \theta_j(t) \phi_j(x)
    = \bigl \langle \theta(t), \, \phi(x) \bigr \rangle.
\end{equation}
We now plug the Galerkin approximation \eqref{eq:Galerkin_approximation}
into the PDE \eqref{eq:advection_diffusion_time_operator},
and we integrate over $\Omega$ against some test function.
We obtain, for all $i \in \{1, \ldots, N\}$,
\begin{equation*}
    \int_\Omega
    \big( \partial_t u_N(t,x) + \mathcal{L}[u_N](t,x) \big)
    \phi_i(x) \, dx = 0.
\end{equation*}
Using the expression \eqref{eq:Galerkin_approximation} yields
\begin{equation*}
    \int_\Omega \sum_{j=1}^N \left(
    \frac{d \theta_j(t)}{dt} \phi_j(x)
    + \theta_j(t) \mathcal{L}[\phi_j](x)
    \right) \!
    \phi_i(x) \, dx = 0,
\end{equation*}
which we recast as
\begin{equation*}
    \sum_{j=1}^N
    \frac{d \theta_j(t)}{dt}
    \int_\Omega \phi_i(x) \phi_j(x) \, dx
    +
    \sum_{j=1}^N
    \theta_j
    \int_\Omega \phi_i(x) \mathcal{L}[\phi_j](x) \, dx
    = 0.
\end{equation*}
We finally obtain a linear ODE on the dofs, which reads
\begin{equation*}
    M \frac{d \theta(t)}{dt} = L \theta(t),
\end{equation*}
where the components of the matrices $M$ and $L$ are defined by
\begin{equation*}
    M_{ij} = \int_\Omega \phi_i(x) \phi_j(x) \, dx
    \text{\qquad and \qquad}
    L_{ij} = \int_\Omega \phi_i(x) \mathcal{L}[\phi_j](x) \, dx.
\end{equation*}
This linear ODE is then solved using a time-stepping method. We also refer to  \cite{quarteroni2008numerical} for more details.

\subsubsection{Semi-Lagrangian Galerkin Method}
\label{sec:SL_galerkin_method}

The Galerkin semi-Lagrangian method is
a subset of the Galerkin method,
specifically designed for
advection and advection-diffusion equations.
The idea is to avoid the time step restriction by using a Lagrangian approach~\cite{douglas1982numerical,pironneau1982transport}.
We still write the approximate solution
under the form \eqref{eq:Galerkin_approximation}.

We introduce the method
in the specific case of a pure transport equation,
i.e., when $\sigma = 0$ in \eqref{eq:advection_diffusion}.
The exact solution satisfies,
for all $t \in (0, T)$, $x \in \Omega$ and $s \in (0, t)$,
\begin{equation}
    \label{eq:transport_equation_solution}
    u(t,x) = u(s, \mathcal{X}(s; t, x)),
\end{equation}
where $\mathcal{X}$ is the (backwards)
characteristic curve implicitly defined by the ODE
\begin{equation}
    \label{eq:characteristic_curve}
    \begin{dcases}
        \frac{d}{ds} \mathcal{X}(s; t, x) = a(s, \mathcal{X}(s; t, x)) &
        \quad \text{for all } s \in (0, t),                              \\
        \mathcal{X}(t; t, x) = x. \vphantom{\dfrac 1 2}
    \end{dcases}
\end{equation}
Using this, we can formulate the Semi-Lagrangian Galerkin method.
To do this, we approximate $u$ by $u_N$
in \eqref{eq:transport_equation_solution},
and integrate over $\Omega$ against a test function.
We thus impose that our approximate solution must satisfy,
for all $t \in (0, T)$ and $s \in (0, t)$,
\begin{equation*}
    \forall i \in \{1, \ldots, N\},
    \quad
    \int_\Omega u_N(t,x) \phi_j(x) \, dx
    =
    \int_\Omega u_N(s, \mathcal{X}(s; t, x)) \phi_j(x) \, dx.
\end{equation*}
Plugging the Galerkin approximation \eqref{eq:Galerkin_approximation}
into the above equation, we obtain
\begin{equation*}
    \forall i \in \{1, \ldots, N\},
    \quad
    \sum_{j=1}^N \theta_j(t) \int_\Omega \phi_j(x) \phi_i(x) \, dx
    =
    \sum_{j=1}^N \theta_j(s) \int_\Omega \phi_j(\mathcal{X}(s; t, x)) \phi_i(x) \, dx.
\end{equation*}
Thus, from some known degrees of freedom $\theta(s)$ at some time $s < t$,
we can compute the degrees of freedom~$\theta(t)$ at time $t$
by solving the above linear system.
This is achieved without
a stability condition linking $t$ to $s$,
contrary to the previous methods.

In practice, we define a time discretization $(t^n)_{n}$,
with $t^{n+1} = t^n + \dt$,
and we obtain $\theta^0 = \theta(0)$ by projecting the initial condition onto the Galerkin basis.
Then, at each time step, $\theta^{n+1}$ is computed by solving
\begin{equation*}
    \forall j \in \{1, \ldots, N\},
    \quad
    \sum_{j=1}^N \theta_j^{n+1} \int_\Omega \phi_j(x) \phi_i(x) \, dx
    =
    \sum_{j=1}^N \theta_j^n \int_\Omega \phi_j(\mathcal{X}(t^n; t^{n+1}, x)) \phi_i(x) \, dx.
\end{equation*}
The above linear system is rewritten as
\begin{equation*}
    M(t^{n+1}, t^{n+1}) \, \theta^{n+1} = M(t^n, t^{n+1}) \, \theta^{n},
\end{equation*}
where the elements of mass matrix $M(s, t)$ are defined by
\begin{equation*}
    M(s, t)_{ij} = \int_\Omega \phi_j(\mathcal{X}(s; t, x)) \phi_i(x) \, dx.
\end{equation*}
Thanks to this formulation, this method does not have a time step restriction on $\Delta t$.
As an example, if $a$ is a constant vector field,
solving \eqref{eq:characteristic_curve} is immediate,
and we get $\mathcal{X}(s; t, x) = x - (t - s) a$.
In this case, the components of the mass matrix simply read
\begin{equation*}
    M(t^n, t^{n+1})_{ij} = \int_\Omega \phi_i(x) \phi_j(x - a \dt) \, dx.
\end{equation*}
Moreover, if a discontinuous Galerkin basis is used,
the mass matrix $M$ becomes diagonal,
rendering the method particularly efficient \cite{CroMehVec2011,rossmanith2011positivity}.

The method can also be extended to
advection-diffusion equations,
see e.g.~\cite{BokSim2016}.
We do not detail this extension here,
but the same principle applies.
Namely, the feet of the characteristic curves
still have to be computed as in the pure advection case.
Then, in the presence of diffusion,
these feet have to be split
into several diffusion directions.

\subsection{Neural methods}
\label{sec:neural_methods}

We now recall some recent neural numerical methods.
In such methods, as mentioned above,
the approximate function space $\mathcal{V}_N^\text{neural}$
is no longer a linear subspace spanned by basis functions,
but a nonlinear subspace (usually, a submanifold)
defined by parameterized nonlinear functions
(usually, neural networks).
However, it remains finite-dimensional,
and the degrees of freedom are the weights of the nonlinear functions.
To be clear, let
\begin{equation}
    \label{eq:neural_network}
    \begin{aligned}
        \mathcal{N}:
        \mathbb{R}^\mathfrak{d} \times \mathbb{R}^N
         & \to
        \mathbb{R} \\
        (X, \theta)
         & \mapsto
        \mathcal{N}(X, \theta)
    \end{aligned}
\end{equation}
be a parameterized nonlinear function,
typically a neural network.
Its input $X \in \mathbb{R}^\mathfrak{d}$
will either be $(t, x) \in (0, T) \times \Omega$,
to define spacetime neural methods,
or $x \in \Omega$,
to define time-discrete neural methods.
In the first case, $\mathfrak{d} = d+1$;
in the second one, $\mathfrak{d} = d$.
Approximating the PDE solution by such
nonlinear functions leads to
a nonlinear approximation space.
Thus, one now expects to solve
nonlinear optimization problems,
rather than linear systems, to find the dofs $\theta$.
Remark that
nonlinear optimization problems have to be solved
even when approximating linear equations.

Note that the dofs of the nonlinear function
$\mathcal{N}$ remain denoted by
$\theta \in \Theta \subset \mathbb{R}^N$.
For instance, if $\mathcal{N}$ is a neural network,
then the dofs $\theta$ are nothing but the weights and biases of the network.
Similarly to \cref{sec:classical_methods},
the dependence of the dofs in space and time
alters the nature of the neural method.
Namely, dofs that depend on neither space nor time
lead to PINNs, discussed in \cref{sec:PINNs},
while time-dependent dofs lead
to discrete PINNs and the neural Galerkin method,
discussed in \cref{sec:neural_Galerkin}.
\rall{A brief comparison of both space-time and time-discrete neural methods
    is provided in \cref{sec:comparison_neural_methods}.}

\subsubsection{Physics-Informed Neural Networks (PINNs)}
\label{sec:PINNs}

To introduce PINNs \cite{RAISSI2019686,NEURIPS2021_df438e52,De-Ryck:2022aa}, let us define the nonlinear subspace
\begin{equation*}
    \mathcal{V}_N^\text{neural} = \left\{
    (x, t)
    \mapsto \mathcal{N}(x, t, \theta),
    \;\; \theta \in \Theta
    \right\}.
\end{equation*}
The function $\mathcal{N}$ depends on both space and time.
To highlight this dependence,
and to distinguish the parameters from
the space and time variables,
we introduce the notation
\begin{equation*}
    \begin{aligned}
        u_\theta:
        (0, T) \times \Omega
         & \to
        \mathbb{R} \\
        (t, x)
         & \mapsto
        \mathcal{N}(t, x, \theta).
    \end{aligned}
\end{equation*}
Then, a function \smash{$u_N \in \mathcal{V}_N^\text{neural}$} is written as
\begin{equation*}
    u_N(t, x) = u_{\theta}(t, x).
\end{equation*}
Note that $u_N$ depends on both space and time,
similarly to the spacetime Galerkin method.

Finding the dofs $\theta$ then consists in
directly minimizing, over the spacetime slab,
the residual of the PDE~\eqref{eq:advection_diffusion_spacetime_operator}
applied to the approximate solution $u_N$.
This leads to the nonlinear optimization problem
\begin{equation}
    \label{eq:PINN_optimization_problem}
    \theta \in \argmin_{\vartheta \in \Theta}
    \int_{\Omega} \int_0^T \bigl| \mathcal{T}[u_{\vartheta}](t,x) \bigr|^2 \, dt \, dx.
\end{equation}
It is usually solved using stochastic gradient methods,
discretizing the integral with the Monte-Carlo algorithm.

\subsubsection{Discrete PINNs and the Neural Galerkin method}
\label{sec:neural_Galerkin}

Discrete PINNs and the Neural Galerkin method
take up the ideas behind the Galerkin method,
but applied to the context of nonlinear function spaces.
Namely, we replace the approximation space with
\begin{equation*}
    \mathcal{V}_N^\text{neural,$t$} = \left\{
    (x, t)
    \mapsto \mathcal{N}(x, \theta(t)),
    \;\; \theta(t) \in \Theta
    \right\}.
\end{equation*}
The function $\mathcal{N}$ still depends on
space, time, and dofs, but its time dependence
is no longer explicit;
rather, it depends on time because the dofs
themselves are time-dependent.
To emphasize this dependence, we introduce the notation
\begin{equation*}
    \begin{aligned}
        u_{\theta(t)}:
        \Omega
         & \to
        \mathbb{R} \\
        x
         & \mapsto
        \mathcal{N}(x, \theta(t)).
    \end{aligned}
\end{equation*}
The approximate function $u_N$ now depends on space only,
while the dofs $\theta$
(i.e., the weights and biases of the neural network if $u_{\theta}$ is a neural network)
depend on time.
This helps to avoid some issues with
causality and time integration
(outlined in e.g.~\cite{MenLiZhaKar2020}),
and leads to the following approximation:
\begin{equation*}
    u_N(t, x) = u_{\theta(t)}(x).
\end{equation*}
The dofs $\theta(t)$ at time $t$ are then determined such that
the residual of equation \eqref{eq:advection_diffusion_time_operator},
\begin{equation}
    \label{eq:neural_methods_ODE_with_x}
    \int_\Omega \Bigl|\partial_t u_{\theta(t)}(x) + \mathcal{L}[u_{\theta(t)}](x)\Bigr|^2 dx ,
\end{equation}
is minimal in $L^2$ norm.
Of course, this minimizer is only approximated,
since the true solution to the PDE
lies in the infinite-dimensional Hilbert space $V$,
while the approximate solution $u_{\theta(t)}$
is in the finite-dimensional subspace
\smash{$\mathcal{V}_N^\text{neural,$t$}$}.
Then, to evolve the dofs in time,
multiple strategies are available,
two of which are detailed in the remainder of this section.
Before that, we remark that
the initial condition $\theta^0 \coloneqq \theta(0)$ is obtained by fitting
the initial condition of the PDE:
\begin{equation*}
    \theta^0 \in \argmin_{\vartheta \in \Theta}
    \bigl| u_{\vartheta}(x) - u_0(x) \bigr|^2 \, dx.
\end{equation*}

\paragraph{Discrete PINNs \cite{StiCha2023,biesek2023burgers}}
Recall the already-introduced time discretization $(t^n)_n$.
Then, discrete PINNs rely on directly
discretizing the time derivative in
\eqref{eq:neural_methods_ODE_with_x}.
For simplicity, we consider a simple explicit Euler scheme,
but more sophisticated schemes can, and should be, used.
This yields, with $\theta^n \coloneqq \theta(t^n)$,
\begin{equation*}
    \partial_t u_{\theta(t)}
    \simeq
    \frac {u_{\theta^{n+1}} - u_{\theta^n}} \dt
    \text{\qquad and \qquad}
    \mathcal{L}[u_{\theta(t)}]
    \simeq
    \mathcal{L}[u_{\theta^n}].
\end{equation*}
Using the above equation, we define the updated parameters $\theta^{n+1}$ as minimizers of the residual:
\begin{equation}
    \label{eq:dPINN_optimization_problem}
    \theta^{n+1} \in \argmin_{\vartheta \in \Theta}
    \int_{\Omega} \, \Bigl| \,
    u_{\vartheta}(x) - u_{\theta^n}(x) + \dt \, \mathcal{L}[u_{\theta^n}](x)
    \, \Bigr|^2 dx.
\end{equation}

\paragraph{Neural Galerkin method \cite{LI2022110958,finzi2023stable,BRUNA2024112588,KAST2024112986}}
To derive the neural Galerkin method,
we go back to \eqref{eq:neural_methods_ODE_with_x}.
Applying the chain rule to the first term gives
\begin{equation*}
    \partial_t u_{\theta(t)}(x)
    = \left \langle
    \nabla_{\theta} u_{\theta(t)}(x),
    \, \frac{d \theta(t)}{dt}
    \right \rangle,
\end{equation*}
where $\nabla_{\theta} u_{\theta}$
is the gradient of $u_{\theta}$ with respect to $\theta$.
Therefore, we seek the dofs $\theta$
whose time derivative minimize, in the $L^2$ norm,
the residual
\begin{equation*}
    \int_\Omega \left|
    \left \langle
    \nabla_{\theta} u_{\theta(t)}(x),
    \, \frac{d \theta(t)}{dt}
    \right \rangle
    + \mathcal{L}[u_{\theta(t)}](x)
    \right|^2 dx.
\end{equation*}
The time derivative is thus defined
as the minimizer of the above residual:
\begin{equation*}
    \frac{d \theta(t)}{dt} \in \argmin_{\eta \in \mathbb{R}^N}
    \int_{\Omega} \, \Bigl| \,
    \left \langle \nabla_{\theta} u_{\theta(t)}(x), \, \eta \right \rangle
    + \mathcal{L}[u_{\theta(t)}](x)
    \, \Bigr|^2 dx.
\end{equation*}
This is nothing but a quadratic optimization problem,
and its exact solution is
\begin{equation}
    \label{eq:neural_Galerkin_ODE}
    M(\theta) \, \frac{d \theta(t)}{dt} = - L(\theta),
\end{equation}
where the matrix $M$ and the vector $L$ are defined by
\begin{equation}
    \label{eq:neural_Galerkin_M_L}
    M(\theta) = \int_{\Omega}
    \nabla_{\theta} u_{\theta(t)}(x) \otimes
    \nabla_{\theta} u_{\theta(t)}(x) \, dx
    \text{\qquad and \qquad}
    L(\theta) = \int_{\Omega}
    \nabla_{\theta} u_{\theta(t)}(x) \,
    \mathcal{L}[u_{\theta(t)}](x) \, dx,
\end{equation}
where $\otimes$ denotes the outer product of two vectors of $\mathbb{R}^N$.
Therefore, finding the dofs in the neural Galerkin method
amounts to solving the nonlinear ODE \eqref{eq:neural_Galerkin_ODE}.
This is usually done using standard ODE solvers.

\subsubsection{Brief comparison of space-time and time-sequential neural methods}
\label{sec:comparison_neural_methods}

\paragraph{\ra{Space-time vs time-sequential methods}}

\ra{Time-sequential methods are less commonly used in the literature than standard PINNs, likely because they are somewhat more challenging to implement. However, they have consistently produced significantly more accurate results as soon as the PDE becomes even moderately difficult. Comparisons have been carried out in particular in \cite{StiCha2023,biesek2023burgers,BRUNA2024112588}. This is consistent with classical numerical methods, where time-sequential approaches are much more commonly used than space-time methods (though this is less true for certain linear PDEs). In practice, the main limitation to the accuracy of neural methods lies in the optimization process. For moderately difficult PDEs, optimization becomes significantly more challenging for space-time approaches (although this can be improved with progressive time training) compared to sequential methods, where optimization is simpler: for instance, a straightforward projection in the case of an explicit discrete PINN, or a least-squares problem for neural Galerkin. The numerical experiments from \cref{sec:applications} include a comparison of these time-sequential methods with space-time PINNs.}

\paragraph{\rc{Discrete PINNs vs Neural Galerkin}}

\rc{The neural Galerkin and discrete PINNs methods can be interpreted within the same framework of time-sequential methods. Discrete PINNs can be viewed as a ``discretize-then-optimize'' strategy, in which the time variable is first discretized, and the solution at the following time step is then obtained via optimization (see~\cite{chen2024teng}). In contrast, the Neural Galerkin method follows an ``optimize-then-discretize'' paradigm. In some cases, it can be derived as a linearization of the discrete PINNs for small time steps.
    In our experiments, in the explicit setting, we observed that the neural Galerkin method is more accurate than discrete PINNs trained with Adam, but less accurate than discrete PINNs trained using a natural gradient approach (see~\cite{chen2024teng}). This suggests that, in the case of Adam, the optimization error may be larger than the error introduced by the linearization in the neural Galerkin method.
    For discrete PINNs, boundary conditions can be enforced either strongly or weakly, as in standard PINNs. However, in the Neural Galerkin framework, weak enforcement of boundary conditions would introduce a non-convex loss, which would undermine the benefits of the method. Therefore, only strong enforcement is generally applicable, which may, in some cases, reduce efficiency.
    Both methods share the same neural architectures, geometry handling, collocation point sampling, and high-dimensional capabilities. As such, the choice between the two approaches will depend primarily on computational cost and accuracy requirements.
    In the implicit setting, discrete PINNs reduce to solving an elliptic PINN problem at each time step, whereas the neural Galerkin method leads to solving a Galerkin problem in the tangent space. The linearization inherent in the neural Galerkin method may limit its performance in the implicit case. Several test cases comparing the two approaches are available in~\cite{chen2024teng}.}
\section{Neural Semi-Lagrangian method}
\label{sec:NGSL}

In this paper, we propose a new approach to solve advection-diffusion equations.
Similarly to how space-time
and Galerkin methods have their neural equivalents,
discussed in \cref{sec:neural_methods},
we propose to extend
the semi-Lagrangian Galerkin method
from \cref{sec:SL_galerkin_method} to the neural domain.

\rb{This approach is quite different from the ones proposed in~\cite{LARIOSCARDENAS2022111623,chen2024conservative,Chen2023ALC}. In these works, the idea is to maintain a grid-based interpolation while using a neural network to correct the classical interpolation. However, in all cases, the approach still relies on grid-based interpolation, which inherits the challenges of classical methods in high dimensions. Even if such methods make it possible to use coarser grids, in high dimensions the grid size remains large, and both the number of degrees of freedom and the computational cost can still grow explosively. }

Similarly to \cref{sec:neural_Galerkin},
the solution will be approximated by a nonlinear function
whose dofs depend on time:
\begin{equation}
    \label{eq:NGSL_approximation}
    u_N(t, x) = u_{\theta(t)}(x).
\end{equation}
In practice, we will use neural networks to represent $u_{\theta}$.

This extension has several advantages compared
to traditional methods and to other traditional and
neural methods:
\begin{itemize}[nosep]
    \item compared to traditional methods,
          it is able to tackle high-dimensional
          problems with a comparatively low cost,
          since the number of degrees of freedom
          does not increase exponentially;
    \item compared to traditional PINNs,
          it remains sequential in time,
          thus avoiding issues observed with PINNs
          and highlighted in e.g.~\cite{MenLiZhaKar2020};
    \item it eliminates the
          (sometimes very restrictive,
          and for the moment unknown)
          stability condition
          on the time step that discrete PINNs
          and the neural Galerkin method retain
          (if explicit time integration is used);
    \item its optimization problem avoids
          the need to compute an additional PDE residual,
          compared to the discrete PINN and neural Galerkin
          methods with implicit time-stepping.
\end{itemize}
However, it should be noted that the semi-Lagrangian
approach is predominantly suited for
advection-diffusion or Fokker-Planck problems
like \eqref{eq:advection_diffusion}.
This restricts its applicability to a
broader range of phenomena.
Moreover, compared to traditional semi-Lagrangian
methods, the use of a nonlinear representation
makes proving convergence estimates much harder,
and pretty much unattainable at this state.

\subsection{Methodology}
\label{sec:NGSL_methodology}

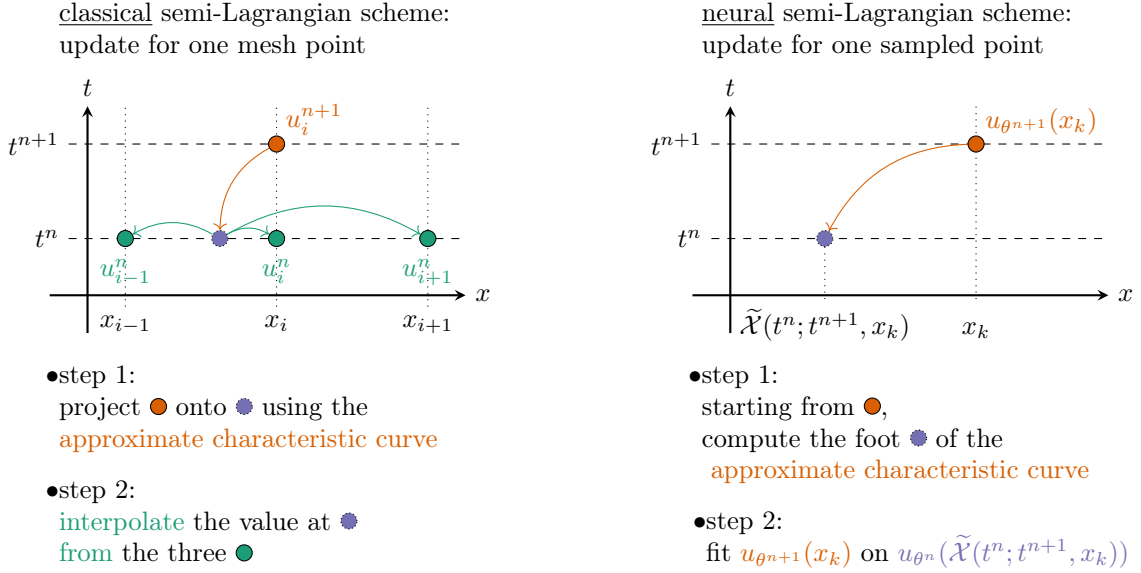
\begin{figure}[!ht]

    \centering

    \begin{tikzpicture}
        \tikzset{
            particle/.style={circle, draw, inner sep=0.75mm},
        }


        \node [anchor=west] at (-0.5,3.5) {\makecell[l]{
                \underline{\smash{classical}} semi-Lagrangian scheme: \\
                update for one mesh point
            }};

        \draw[-stealth, thick] (-0.5,0) -- (5,0) node[right]{$x$};
        \draw[-stealth, thick] (0,-0.5) -- (0,2.5) node[above]{$t$};

        \draw[dashed] (-0.25, 0.75) node[left]{$t^n$} -- (5, 0.75) ;
        \draw[dashed] (-0.25, 2) node[left]{$t^{n+1}$} -- (5, 2) ;

        \draw[dotted] (0.5, -0.15) node[below]{$x_{i-1}$} -- (0.5, 2.5) ;
        \draw[dotted] (2.5, -0.15) node[below]{$x_i$} -- (2.5, 2.5) ;
        \draw[dotted] (4.5, -0.15) node[below]{$x_{i+1}$} -- (4.5, 2.5) ;

        \coordinate (u_i_np1_coord) at (2.5,2);
        \coordinate (u_i_np1_bwd_coord) at (1.75,0.75);
        \coordinate (u_im1_n_coord) at (0.5,0.75);
        \coordinate (u_i_n_coord) at (2.5,0.75);
        \coordinate (u_ip1_n_coord) at (4.5,0.75);

        \node[particle, fill=graph_2] (u_i_np1) at (u_i_np1_coord){};
        \node[above right, graph_2] at (u_i_np1_coord) {$u_{i}^{n+1}$};

        \node[particle, fill=graph_3] (u_im1_n) at (u_im1_n_coord){};
        \node[below, yshift=-4pt, graph_3] at (u_im1_n_coord) {$u_{i-1}^{n}$};
        \node[particle, fill=graph_3] (u_i_n) at (u_i_n_coord){};
        \node[below, yshift=-4pt, graph_3] at (u_i_n_coord) {$u_{i}^{n}$};
        \node[particle, fill=graph_3] (u_ip1_n) at (u_ip1_n_coord){};
        \node[below, yshift=-4pt, graph_3] at (u_ip1_n_coord) {$u_{i+1}^{n}$};

        \node[particle, densely dotted, fill=graph_1] (u_i_np1_bwd) at (u_i_np1_bwd_coord){};

        \draw[->, graph_2] (u_i_np1) to [bend right=30] (u_i_np1_bwd);

        \draw[->, graph_3] (u_i_np1_bwd) to [bend right=30] (u_im1_n);
        \draw[->, graph_3] (u_i_np1_bwd) to [bend left=30] (u_i_n);
        \draw[->, graph_3] (u_i_np1_bwd) to [bend left=30] (u_ip1_n);

        \node [anchor=west] (legend_step_1) at (-0.5,-1.5) {\makecell[l]{
                \makebox[0pt][r]{\textbullet}step 1: \\
                project
                \tikz{\node[particle, fill=graph_2] at (0,0) {}}
                onto
                \tikz{\node[particle, densely dotted, fill=graph_1] at (0,0) {}}
                using the \\
                \textcolor{graph_2}{approximate characteristic curve}
            }};
        \node [below] (legend_step_2) at (legend_step_1.south) {\makecell[l]{
                \makebox[0pt][r]{\textbullet}\sbox0{approximate characteristic curve}\makebox[\wd0][l]{step 2:} \\
                \textcolor{graph_3}{interpolate} the value at
                \tikz{\node[particle, densely dotted, fill=graph_1] at (0,0) {}} \\
                \textcolor{graph_3}{from} the three
                \tikz{\node[particle, fill=graph_3] at (0,0) {}}
            }};


        \pgfmathsetmacro{\dx}{8.5}

        \node [anchor=west] at (-0.5+\dx,3.5) {\makecell[l]{
                \underline{\smash{neural}} semi-Lagrangian scheme: \\
                update for one sampled point
            }};

        \draw[-stealth, thick] (-0.5+\dx,0) -- (5+\dx,0) node[right]{$x$};
        \draw[-stealth, thick] (0+\dx,-0.5) -- (0+\dx,2.5) node[above]{$t$};

        \draw[dashed] (-0.25+\dx, 0.75) node[left]{$t^n$} -- (5+\dx, 0.75) ;
        \draw[dashed] (-0.25+\dx, 2) node[left]{$t^{n+1}$} -- (5+\dx, 2) ;
        \draw[dotted] (3.25+\dx, -0.05) node[below]{$x_k\vphantom{\widetilde{\mathcal{X}}(t^n; t^{n+1}, x_k)}$} -- (3.25+\dx, 2.5) ;

        \coordinate (u_i_np1_coord) at (3.25+\dx,2);
        \coordinate (u_i_np1_bwd_coord) at (1.25+\dx,0.75);

        \node[particle, fill=graph_2] (u_i_np1) at (u_i_np1_coord){};
        \node[above right, graph_2] at (u_i_np1_coord) {$u_{\theta^{n+1}}(x_k)$};

        \node[particle, densely dotted, fill=graph_1] (u_i_np1_bwd) at (u_i_np1_bwd_coord){};

        \draw[->, graph_2] (u_i_np1) to [bend right=30] (u_i_np1_bwd);
        \draw[dotted] (1.25+\dx, -0.05) node[below]{$\widetilde{\mathcal{X}}(t^n; t^{n+1}, x_k)$} -- (u_i_np1_bwd) ;

        \node [anchor=west] (legend_step_1) at (-0.5+\dx,-1.705) {\makecell[l]{
                \makebox[0pt][r]{\textbullet}step 1: \\
                starting from
                \tikz{\node[particle, fill=graph_2] at (0,0) {}}, \\
                compute the foot
                \tikz{\node[particle, densely dotted, fill=graph_1] at (0,0) {}}
                of the \\
                \sbox0{fit $u_{\theta^{n+1}}(x_k)$ on \smash{$u_{\theta^n}(\widetilde{\mathcal{X}}(t^n; t^{n+1}, x_k))$}}
                \makebox[\wd0][l]{\textcolor{graph_2}{approximate characteristic curve}}
            }};
        \node [below] (legend_step_2) at (legend_step_1.south) {\makecell[l]{
        \makebox[0pt][r]{\textbullet}\sbox0{approximate characteristic curve}\makebox[\wd0][l]{step 2:} \\
        fit $\color{graph_2}u_{\theta^{n+1}}(x_k)$
        on \smash{$\color{graph_1}u_{\theta^n}(\widetilde{\mathcal{X}}(t^n; t^{n+1}, x_k))$}
        }};

    \end{tikzpicture}

    \caption{\rc{Visual comparison of the classical (left)
            and neural (right) semi-Lagrangian schemes
            in one space dimension.
            In the mesh-based classical scheme, the update is performed on a mesh,
            and thus an interpolation onto the mesh is needed after
            computing the foot of the characteristic curve.
            Conversely, in the meshless neural scheme,
            the approximation at the previous time step
            is directly evaluated at the foot of the characteristic curve.}}
    \label{fig:NSL_scheme_drawing}

\end{figure}

Recall that the solution to the advection equation
satisfies \eqref{eq:transport_equation_solution}.
Hence, we seek an approximate solution~$u_N$ that also satisfies it:
\begin{equation*}
    \forall t \in (0, T), \,
    \forall s \in (0, t), \,
    \forall x \in \Omega, \qquad
    u_N(t, x) \simeq u_N(s, \mathcal{X}(s; t, x)).
\end{equation*}
Plugging \eqref{eq:NGSL_approximation}, we obtain
\begin{equation}
    \label{eq:NGSL_approximate_solution}
    \forall t \in (0, T), \,
    \forall s \in (0, t), \,
    \forall x \in \Omega, \qquad
    u_{\theta(t)}(x) \simeq u_{\theta(s)}(\mathcal{X}(s; t, x)).
\end{equation}
The discrete version of
\eqref{eq:NGSL_approximate_solution} reads,
using the time discretization~$(t^n)_n$,
\begin{equation*}
    \forall n \geqslant 0, \,
    \forall x \in \Omega, \qquad
    u_{\theta(t^{n+1})}(x)
    \simeq
    u_{\theta(t^{n})}(\mathcal{X}(t^n; t^{n+1}, x)).
\end{equation*}
We wish to approximately satisfy
the above equation in the $L^2$ norm.
Therefore, we define the dofs
$\theta^{n+1} = \theta(t^{n+1})$
at time $t^{n+1}$
as the minimizer of the following optimization problem:
\begin{equation}
    \label{eq:NGSL_nonlinear_optimization}
    \theta^{n+1} \in \argmin_{\vartheta \in \Theta}
    \int_{\Omega} \bigl|
    u_{\vartheta}(x) - u_{\theta^n}(\mathcal{X}(t^n; t^{n+1}, x))
    \bigr|^2 \, dx.
\end{equation}
Just like the semi-Lagrangian Galerkin method,
this method does not have a stability condition on the time step.
\rc{The newly introduced neural (meshless) semi-Lagrangian scheme
    is visually compared to a classical (mesh-based) one
    in \cref{fig:NSL_scheme_drawing}.}

In the linear advection case,
we recall that the characteristic curve is given by
$\mathcal{X}(s; t, x) = x - (t - s) a$.
Therefore, the optimization problem reads
\begin{equation*}
    \theta^{n+1} \in \argmin_{\vartheta \in \Theta}
    \int_{\Omega} \bigl|
    u_{\vartheta}(x) - u_{\theta^n}(x - a \Delta t)
    \bigr|^2 \, dx.
\end{equation*}
In other cases, the ODE \eqref{eq:characteristic_curve}
defining the characteristic curves has to be solved,
either exactly or approximately, to compute $\mathcal{X}$.
This is detailed in the following section.

\subsection{Using the method in practice}
\label{sec:NGSL_algorithms}

In practice, we do not exactly solve
the optimization problem \eqref{eq:NGSL_nonlinear_optimization}.
Indeed, we have to make three adjustments
to make this problem computationally feasible.
\begin{enumerate}
    \item As usual, the integrals are approximated using the Monte-Carlo method.
          Given a number $N_c$ of samples, called ``collocation points''
          and denoted by $(x_k)_{k \in \{1, \dots, N_c\}}$,
          we write
          \begin{equation*}
              \frac{1}{N_c} \sum_{k=1}^{N_c} \bigl|
              u_{\vartheta}(x_k) - u_{\theta^n}(\mathcal{X}(t^n; t^{n+1}, x_k))
              \bigr|^2
              \simeq
              \int_{\Omega} \bigl|
              u_{\vartheta}(x) - u_{\theta^n}(\mathcal{X}(t^n; t^{n+1}, x))
              \bigr|^2 \, dx.
          \end{equation*}
          This approximation comes with an error in $\mathcal{O}(1/\sqrt{N_c})$.
    \item A time-stepping algorithm (e.g.\ a Runge-Kutta method)
          is employed to solve the
          ODE~\eqref{eq:characteristic_curve}
          describing the characteristic curves.
          This leads to an approximate characteristic curve,
          which is denoted by
          \begin{equation*}
              \widetilde{\mathcal{X}} \approx \mathcal{X}.
          \end{equation*}
          Note that \smash{$\widetilde{\mathcal{X}} = \mathcal{X}$}
          as soon as the ODE has a closed-form solution.
          \ra{Because of the Monte-Carlo integration,
              this characteristic solver is applied to integration points that are randomly sampled.
              When the velocity field has an analytic expression,
              this part is fully vectorizable and easily parallelizable
              by assigning subsets of points to different processes as needed.
              If it does not have an analytic expression
              (e.g.\ given by another differential equation, see \cref{sec:vlasov_poisson}),
              it has to be computable at any point in the domain
              to benefit from the vectorization and parallelization.
              For instance, it could be given by a spline or a neural network.
              In practice, solving the ODE represents a
              negligible part of the overall computation cost.}
    \item The optimization problem is solved using
          a stochastic gradient descent method.
          This is done by computing the gradient
          of the objective function with respect to $\theta$.
          This gradient is computed using the chain rule,
          and it requires computing the gradient of $u_{\theta}$ with respect to $\theta$,
          which is done using automatic differentiation.
\end{enumerate}

We now summarize the Neural Semi-Lagrangian (NSL)
method in the following algorithm,
in the case of a parametric advection equation.
It is detailed in \cref{alg:NSL_known_characteristic},
while the computation of the approximate characteristic curve is more classical,
and recalled in \cref{alg:characteristic_curve} for completeness.

\begin{algorithm}[!ht]
    \caption{Neural Semi-Lagrangian (NSL) method for advection equations}
    \label{alg:NSL_known_characteristic}
    \begin{algorithmic}[1]
        \State \textbf{Input:} Initial condition $u_0(x, \mu)$, time discretization $(t^n)_n$, neural network architecture, final time $T$, number of collocation points $N_c$
        \State \textbf{Output:} Approximate solution $u_N(t^n, x, \mu)$ for all $t^n$, all $x \in \Omega$, and all $\mu \in \mathbb{M}$
        \medskip
        \State \textbf{Initialization:}
        Compute the initial dofs $\theta^0$ by randomly sampling collocation points
        $(x_k)_k \in \Omega^{N_c}$
        and collocation parameters
        $(\mu_k)_k \in \mathbb{M}^{N_c}$ and solving
        \[
            \theta^0 \in \argmin_{\vartheta \in \Theta}
            \sum_{k=1}^{N_c} \,
            \bigl| u_{\vartheta}(x_k, \mu_k) - u_0(x_k, \mu_k) \bigr|^2
        \]
        \medskip
        \While{$t < T$}
        \State Solve the nonlinear optimization problem
        \[
            \theta^{n+1} \in \argmin_{\vartheta \in \Theta}
            \sum_{k=1}^{N_c} \,
            \bigl|
            u_{\vartheta}(x_k, \mu_k) -
            u_{\theta^n}(\widetilde{\mathcal{X}}(t^n; t^{n+1}, x_k, \mu_k))
            \bigr|^2
        \]
        \hspace{\algorithmicindent}
        \For{each iteration of the nonlinear
            optimization algorithm}
        \State
        Randomly sample collocation points
        $(x_k)_k \in \Omega^{N_c}$
        and collocation parameters
        $(\mu_k)_k \in \mathbb{M}^{N_c}$
        \State
        Compute the approximate foot
        \smash{$\widetilde{\mathcal{X}}(t^n; t^{n+1}, x_k, \mu_k)$}
        of the characteristic curve
        using \cref{alg:characteristic_curve}
        \EndFor
        \State Update the time step: $t \leftarrow t + \Delta t$
        \EndWhile
        \medskip
        \State \textbf{Return:} $u_N(t, x, \mu) = u_{\theta(t)}(x)$
    \end{algorithmic}
\end{algorithm}

\begin{remark}
    \cref{alg:NSL_known_characteristic} corresponds
    to the backwards semi-Lagrangian method,
    where the characteristic curves
    are integrated backwards in time.
    While the backwards method is well-known
    within the realm of traditional numerical methods,
    the forward version of the semi-Lagrangian method
    also exists, see~\cite{CroResSon2009}.
    It consists in integrating the characteristic curves
    forwards in time rather than backwards.
    It would also be possible to develop a forwards version
    of the present method, simply by
    replacing~\eqref{eq:NGSL_nonlinear_optimization} with
    \begin{equation*}
        \label{eq:backwards_NGSL_nonlinear_optimization}
        \theta^{n+1} \in \argmin_{\vartheta \in \Theta}
        \int_{\Omega} \bigl|
        u_{\vartheta}({\mathcal{X}^\text{fwd}}(t^{n+1}; t^{n}, x)) - u_{\theta^n}(x)
        \bigr|^2 \, dx,
    \end{equation*}
    where \smash{${\mathcal{X}^\text{fwd}}$}
    is the (forwards) characteristic curve
    originating in point $x$ at time $t^n$, and
    \smash{${\mathcal{X}^\text{fwd}}(t^{n+1}; t^{n}, x)$}
    is the position of the head at time $t^{n+1}$.
    Both forwards and backwards versions of
    the present method should perform the same,
    and we focus on the backwards version
    to simplify the presentation.
\end{remark}

\subsection{\rb{Rough error estimates}}
\label{sec:NGSL_error}

\rb{In this section, we present a rough estimate of the error
    made by the NSL method
    in order to decompose it into several contributions:
    integration, optimization, approximation and characteristic errors.} \rc{Indeed, at each iteration, a projection problem onto the neural network function space is solved. The error thus depends on the approximation capacity of neural networks (approximation error) as well as the error at the previous step. Other sources of error come from the fact that the projection problem is actually modified: the integral is replaced by a Monte Carlo estimate (integration error) and the characteristic curves are solved numerically if they are not explicitly known (characteristic error). Finally, the optimization problem is only solved approximately (optimization error). Errors then accumulate over the course of iterations.}
\rb{We do not expect the resulting estimate to give practical control over the error;
    rather, we aim at a better understanding of the error sources.}

We do not consider the dependence in the parameters $\mu$, for simplicity.
Moreover, we consider the pure transport case,
i.e., we take $\sigma = 0$ in \eqref{eq:advection_diffusion}.
We seek to estimate
\begin{equation*}
    \mathcal{E}^{n+1} = \|
    u_{\theta^{n+1}} - u(t^{n+1}, \cdot)
    \|_{L^2(\Omega)},
\end{equation*}
where $u$ is the exact PDE solution,
and where $u_{\theta^{n+1}}$ is the approximate solution
at time~$t^{n+1}$.
It is obtained by solving, in an approximate manner,
the nonlinear optimization problem
from \cref{alg:NSL_known_characteristic}, which reads
\begin{equation}
    \label{eq:NGSL_nonlinear_optimization_in_error_estimate}
    \argmin_{{u \in \mathcal{V}_N^\text{neural,\smash{$t^{n+1}$}}}}
    \frac 1 {N_c} \sum_{k=1}^{N_c} \bigl|
    u(x_k) - u_{\theta^n}(\widetilde{\mathcal{X}}(t^n; t^{n+1}, x_k))
    \bigr|^2,
\end{equation}
where \smash{$\mathcal{V}_N^\text{neural,$t^{n+1}$} \subset L^2(\Omega)$}
is the nonlinear approximation space, whose solution is denoted $u_{\theta^{n+1}_{\ast}}$.
Let us emphasize that, up to an optimization error to be discussed later,
we have
\begin{equation*}
    u_{\theta^{n+1}} \approx u_{\theta^{n+1}_{\ast}},
\end{equation*}
since $u_{\theta^{n+1}}$ is an approximate minimizer of
\eqref{eq:NGSL_nonlinear_optimization_in_error_estimate},
while $u_{\theta^{n+1}_{\ast}}$ is the true minimizer.
Further, note that
\begin{equation*}
    u^{0} \in \argmin_{{u \in \mathcal{V}_N^\text{neural,\smash{$0$}}}}
    \frac 1 {N_c} \sum_{k=1}^{N_c} \bigl| u(x_k) - u_0(x_k) \bigr|^2.
\end{equation*}

To bound the error $\mathcal{E}^{n+1}$,
we define three quantities,
respectively corresponding to the
integration, optimization and approximation errors.
In these definitions,
$(x_k)_{k \in \{1, \dots, N_c\}} \in \Omega^{N_c}$
form a set of $N_c$ collocation points. For any $f$ in $V$, the \emph{integration error} is defined by
\begin{equation}
    \varepsilon_\text{int} \bigl(f; (x_k)_k \bigr)
    =
    \left|
    \bigl\|
    f
    \bigr\|_{L^2(\Omega)}
    -
    \frac 1 {N_c} \sum_{k=1}^{N_c} \bigl|
    f(x_k)
    \bigr|^2
    \right|.\label{def:integration_error}
\end{equation}
This is nothing
but the error made when approximating the integral of $f$
using the Monte-Carlo method. We next define the (discrete)
\emph{optimization error} between the exact minimizer $f_N^\star$ and the  approximate minimizer~$f_N$
of the same optimization problem:
\begin{equation}
    \varepsilon_\text{opt}
    \bigl(f_N, f_N^\star; (x_k)_k \bigr)
    =
    \frac 1 {N_c} \sum_{k=1}^{N_c} \bigl|
    f_N(x_k) - f_N^\star(x_k)
    \bigr|^2.\label{def:optimization_error}
\end{equation}
This measures
how close the approximate minimizer
is to the exact minimizer of the optimization problem,
using a Monte-Carlo approximation.
The \emph{approximation error} is defined by
\begin{equation}
    \varepsilon_\text{approx} \bigl(f_N, f \bigr)
    =
    \bigl\|
    f_N - f
    \bigr\|_{L^2(\Omega)}.\label{def:approximation_error}
\end{equation}
This measures how well a function in
the nonlinear subspace $\mathcal{V}_N^\text{neural, t}$
approximates a given function in the infinite-dimensional space $V$.
It is expected to decrease as the number $N$ of dofs increases.

Equipped with these definitions,
we are almost ready to state the main result of this section.
Before that, we state and prove the following lemma.

\begin{lemma}
    \label{lem:characteristic_is_diffeomorphism}
    Let $n \geqslant 0$, and define the function
    $\mathcal{X}^{n+1}: x \mapsto \mathcal{X}(t^n; t^{n+1}, x)$.
    Then, for all $f \in V$,
    \begin{equation*}
        \bigl\|
        f \circ \mathcal{X}^{n+1}
        \bigr\|_{L^2(\Omega)}
        \leqslant
        \frac 1 {d_\mathcal{X}^{n+1}}
        \bigl\|
        f
        \bigr\|_{L^2(\Omega)}
        \text{, \quad with }
        d_\mathcal{X}^{n+1} = \sqrt{
            \inf_{y \in \Omega}
            \bigl| \det
            \left(
            D \mathcal{X}^{n+1}
            \right) (y) \bigr|
        },
    \end{equation*}
\end{lemma}

\begin{proof}
    The function $\mathcal{X}^{n+1}$ is a diffeomorphism
    from $\Omega$ to itself, since it is the flow
    of the ODE \eqref{eq:characteristic_curve},
    see e.g.~\cite{HirSmaDev2013}.
    Therefore, we have
    \begin{equation*}
        \int_\Omega f(\mathcal{X}^{n+1}(x))^2 \, dx
        =
        \int_{\mathcal{X}^{n+1}(\Omega)} f(y)^2
        | \det (D \mathcal{X}^{n+1})^{-1}(y) | \, dy
        =
        \int_{\Omega} f(y)^2
        \frac {dy} {| \det (D \mathcal{X}^{n+1})(y) |}.
    \end{equation*}
    Introducing $d_\mathcal{X}^{n+1}$ and taking the square root
    of the above equation gives the desired result.
\end{proof}

\begin{proposition}
    \label{thm:NGSL_error}
    Let $u$ be the exact solution of
    the transport equation
    (i.e., equation \eqref{eq:advection_diffusion} with $\sigma = 0$)
    on a bounded domain~$\Omega$,
    and let $u_{\theta^n}$ be the approximate solution
    at time $t^{n}$.
    Let $(x_k)_{k \in \{1, \dots, N_c\}} \in \Omega^{N_c}$
    be a set of $N_c$ collocation points. Let $p$ and $n_\tau$ be, respectively, the order and the number of sub-time steps of the characteristic curve solver.
    We have, for all~$n$ such that $t^n < T$,
    \begin{equation}
        \label{eq:NGSL_error_estimate}
        \|
        u_{\theta^n} - u(t^{n}, \cdot)
        \|_{L^2(\Omega)}
        \leqslant
        \sum_{\mathfrak{n}=0}^{n}
        \left(
        \varepsilon_\text{int}^{\mathfrak{n}}
        + \varepsilon_\text{opt}^{\mathfrak{n}}
        + \varepsilon_\text{approx}^{\mathfrak{n}}
        \right)
        \mathcal{D}^{\mathfrak{n}}
        +
        \biggl(
        \sum_{\mathfrak{n}=1}^{n}
        \mathcal{C}^{\mathfrak{n}}
        \mathcal{D}^{\mathfrak{n}}
        \biggr)
        \left( \frac {\Delta t}{n_\tau} \right)^{p+1} + \mathcal{O}\bigg(
        \left( \frac {\Delta t}{n_\tau} \right)^{2p+2}\bigg),
    \end{equation}
    which involves the error terms at each iteration
    \begin{equation*}
        \begin{aligned}
            \varepsilon_\text{int}^{\mathfrak{n}}
             & =
            \varepsilon_\text{int}
            \bigl(u_{\theta^{\mathfrak{n}}} - u_{\theta_{\ast}^{\mathfrak{n}}}; (x_k)_k \bigr), \\
            \varepsilon_\text{opt}^{\mathfrak{n}}
             & =
            \varepsilon_\text{opt}
            \bigl(u_{\theta^{\mathfrak{n}}}, u_{\theta_{\ast}^{\mathfrak{n}}}; (x_k)_k \bigr),
        \end{aligned}\qquad\qquad
        \varepsilon_\text{approx}^{\mathfrak{n}}
        =
        \begin{dcases}
            \varepsilon_\text{approx}
            \bigl(
            u_{\theta_{\ast}^{\mathfrak{n}}},
            u_{\theta^{\mathfrak{n}-1}} \circ \widetilde{\mathcal{X}}(t^{\mathfrak{n}-1}; t^{\mathfrak{n}}, \cdot)
            \bigr)
             & \text{ if } \mathfrak{n} > 0, \\
            \varepsilon_\text{approx}
            \bigl(
            u_{\theta_{\ast}^{\mathfrak{n}}},
            u_0
            \bigr)
             & \text{ otherwise,}            \\
        \end{dcases}
    \end{equation*}
    where the errors are defined by \eqref{def:integration_error}-\eqref{def:optimization_error}-\eqref{def:approximation_error} and the terms
    \begin{equation*}
        \mathcal{C}^{\mathfrak{n}}
        =
        \frac {M_p} {d_\mathcal{X}^{\mathfrak{n}}}
        \bigl\| (\mathcal{X}^{\mathfrak{n}})^{(p+1)} \bigr\|_{L^\infty(\Omega)},
        \bigl\| \nabla u_{\theta^{\mathfrak{n}-1}} \bigr\|_{L^2(\Omega)},\quad
        \quad\text{ and }\quad
        \mathcal{D}^{\mathfrak{n}} = \prod_{m=\mathfrak{n}+1}^{n} \frac {1}{d_\mathcal{X}^{m}},
    \end{equation*}
    where $M_p$ is a constant associated with the characteristic curve solver and $d_\mathcal{X}^{\mathfrak{n}}$ is defined in \cref{lem:characteristic_is_diffeomorphism}.
\end{proposition}

\begin{proof}
    We wish to bound the error $\mathcal{E}^{n+1} = \|
        u_{\theta^{n+1}} - u(t^{n+1}, \cdot)
        \|_{L^2(\Omega)}$.
    We first note that,
    since $u$ is the exact solution of the advection equation,
    with characteristic curve $\mathcal{X}$,
    \begin{equation*}
        \mathcal{E}^{n+1}
        = \|
        u_{\theta^{n+1}} - u(t^{n+1}, \cdot)
        \|_{L^2(\Omega)}
        = \|
        u_{\theta^{n+1}} - u(t^{n}, \cdot) \circ \mathcal{X}^{n+1}
        \|_{L^2(\Omega)},
    \end{equation*}
    where we have defined $\mathcal{X}^{n+1}: x \mapsto \mathcal{X}(t^n; t^{n+1}, x)$.
    Now, we split the error into a telescopic sum and obtain:
    \begin{equation*}
        \begin{aligned}
            \mathcal{E}^{n+1}
            \leqslant
             & \bigl\| u_{\theta^{n+1}} - u_{\theta^{n+1}_{\ast}}  \bigr\|_{L^2(\Omega)}                                             \\
             & \phantom{\bigl\| u_{\theta^{n+1}}}
            +\bigl\| u_{\theta^{n+1}_{\ast}} - u_{\theta^n} \circ \widetilde{\mathcal{X}}^{n+1}        \bigr\|_{L^2(\Omega)}         \\
             & \phantom{\bigl\| u_{\theta^{n+1}} + u_{\theta^{n+1}_{\ast}}}
            + \bigl\| u_{\theta^n} \circ \widetilde{\mathcal{X}}^{n+1} - u_{\theta^n} \circ \mathcal{X}^{n+1}  \bigr\|_{L^2(\Omega)} \\
             & \phantom{\bigl\|  u_{\theta^{n+1}} + u_{\theta^{n+1}_{\ast}} + u_{\theta^n} \circ \widetilde{\mathcal{X}}^{n+1}}
            + \bigl\| u_{\theta^n} \circ \mathcal{X}^{n+1} - u(t^{n}, \cdot) \circ \mathcal{X}^{n+1}
            \bigr\|_{L^2(\Omega)} =
            \alpha + \beta + \gamma + \delta,
        \end{aligned}
    \end{equation*}
    %
    with $\alpha, \beta, \gamma, \delta$ respectively referring to the four terms of the upper bound. We now bound each of them separately.

    First, with a triangle inequality, $\alpha$ satisfies
    \begin{equation*}
        \alpha \leqslant
        \left|
        \bigl\|
        u_{\theta^{n+1}} - u_{\theta^{n+1}_{\ast}}
        \bigr\|_{L^2(\Omega)}
        -
        \frac 1 {N_c} \sum_{k=1}^{N_c} \bigl|
        u_{\theta^{n+1}}(x_k) - {u}^{n+1}(x_k)
        \bigr|^2
        \right|
        +
        \frac 1 {N_c} \sum_{k=1}^{N_c} \bigl|
        u_{\theta^{n+1}}(x_k) - {u}^{n+1}(x_k)
        \bigr|^2.
    \end{equation*}
    Thus, $\alpha$ is directly bounded by
    the integration and optimization errors,
    and we obtain
    \begin{equation*}
        \alpha \leqslant
        \varepsilon_\text{int}^{n+1} + \varepsilon_\text{opt}^{n+1}.
    \end{equation*}
    Furthermore, we directly note that
    $\beta = \varepsilon_\text{approx}^{n+1}$.

    Then, the term $\gamma$ corresponds to the error
    made when approximating the characteristic curve.
    The local truncation error of the numerical solver writes:
    \begin{equation*}
        \widetilde{\mathcal{X}}(t^n;t^{n+1},x) =
        \mathcal{X}(t^n;t^{n+1},x) + M_p\, \mathcal{X}^{(p+1)}(t^n;t^{n+1},x) \left(\frac{\Delta t}{n_\tau}\right)^{p+1} + \mathcal{O}\bigg(
        \left( \frac {\Delta t}{n_\tau} \right)^{2p+2}\bigg)
    \end{equation*}
    where $M_p$ is a positive real constant and $(\mathcal{X}^{n+1})^{(p+1)}$ is the $(p\!+\!1)$-th derivative of
    the function $t \mapsto \mathcal{X}(t; t^{n+1}, x)$. Thus, performing a Taylor expansion, we have
    \begin{align*}
        \gamma & = \Bigl\|
        u_{\theta^n} \circ \widetilde{\mathcal{X}}^{n+1} - u_{\theta^n} \circ \mathcal{X}^{n+1}
        \Bigr\|_{L^2(\Omega)} \vphantom{\dfrac 1 2}          \\
               & = \Bigl\|
        u_{\theta^n} \bigg(\mathcal{X}^{n+1} + M_p\, \mathcal{X}^{(p+1)}(t^n;t^{n+1},\cdot) \left(\frac{\Delta t}{n_\tau}\right)^{p+1} + \mathcal{O}\bigg(
            \left( \frac {\Delta t}{n_\tau} \right)^{2p+2} \bigg) \bigg) - u_{\theta^n} \circ \mathcal{X}^{n+1}
        \Bigr\|_{L^2(\Omega)}                                \\
               & \leqslant \Bigl\|
        \Big(\nabla u_{\theta^n} \circ \mathcal{X}^{n+1}\Big)  \, M_p\, \mathcal{X}^{(p+1)}(t^n;t^{n+1},\cdot)  \left(\frac{\Delta t}{n_\tau}\right)^{p+1}
        \Bigr\|_{L^2(\Omega)} + \mathcal{O}\bigg(
        \left( \frac {\Delta t}{n_\tau} \right)^{2p+2}\bigg) \\
               & \leqslant
        \Bigl\| \nabla u_{\theta^n} \circ \mathcal{X}^{n+1} \Bigr\|_{L^2(\Omega)} \, M_p\, \Bigl\| \mathcal{X}^{(p+1)}(t^n;t^{n+1},\cdot)\Bigr\|_{L^\infty(\Omega)}    \left(\frac{\Delta t}{n_\tau}\right)^{p+1} + \mathcal{O}\bigg(
        \left( \frac {\Delta t}{n_\tau} \right)^{2p+2}\bigg),
    \end{align*}
    where the higher order terms are lumped into the last term. Arguing \cref{lem:characteristic_is_diffeomorphism},
    we further bound the gradient norm:
    \begin{equation*}
        \bigl\| \nabla u_{\theta^n} \circ \mathcal{X}^{n+1} \bigr\|_{L^2(\Omega)}
        \leqslant
        \frac 1 {d_\mathcal{X}^{n+1}}
        \bigl\| \nabla u_{\theta^n} \bigr\|_{L^2(\Omega)}.
    \end{equation*}
    All in all, we obtain the bound
    \begin{equation*}
        \gamma
        \leqslant
        \frac {M_p} {d_\mathcal{X}^{n+1}}
        \bigl\| (\mathcal{X}^{n+1})^{(p+1)} \bigr\|_{L^\infty(\Omega)}
        \bigl\| \nabla u_{\theta^n} \bigr\|_{L^2(\Omega)}
        \left( \frac {\Delta t}{n_\tau} \right)^{p+1} + \mathcal{O}\bigg(
        \left( \frac {\Delta t}{n_\tau} \right)^{2p+2}\bigg).
    \end{equation*}

    Lastly, using \cref{lem:characteristic_is_diffeomorphism} on $\delta$,
    we obtain
    \begin{equation*}
        \delta =
        \bigl\|u_{\theta^n} \circ \mathcal{X}^{n+1} - u(t^{n}, \cdot) \circ \mathcal{X}^{n+1}
        \bigr\|_{L^2(\Omega)}
        \leqslant
        \frac 1 {d_\mathcal{X}^{n+1}}
        \bigl\|u_{\theta^n} - u(t^{n}, \cdot) \bigr\|_{L^2(\Omega)}
        =
        \frac 1 {d_\mathcal{X}^{n+1}}
        \mathcal{E}^{n}.
    \end{equation*}

    Putting everything together, we obtain
    \begin{equation*}
        \mathcal{E}^{n+1}
        \leqslant
        \varepsilon_\text{int}^{n+1}
        + \varepsilon_\text{opt}^{n+1}
        + \varepsilon_\text{approx}^{n+1}
        + \mathcal{C}^{n+1} \left( \frac {\Delta t}{n_\tau} \right)^{p+1}
        + \mathcal{O}\bigg(
        \left( \frac {\Delta t}{n_\tau} \right)^{2p+2}\bigg)
        + \frac {\mathcal{E}^{n}}{d_\mathcal{X}^{n+1}},
    \end{equation*}
    where we have defined
    \begin{equation*}
        \mathcal{C}^{n+1}
        =
        \frac {M_p} {d_\mathcal{X}^{n+1}}
        \bigl\| (\mathcal{X}^{n+1})^{(p+1)} \bigr\|_{L^\infty(\Omega)}
        \bigl\| \nabla u_{\theta^n} \bigr\|_{L^2(\Omega)}.
    \end{equation*}
    Using the discrete version of Grönwall's lemma, we obtain
    \begin{equation*}
        \mathcal{E}^{n}
        \leqslant
        \sum_{\mathfrak{n}=1}^{n}
        \left(
        \varepsilon_\text{int}^{\mathfrak{n}}
        + \varepsilon_\text{opt}^{\mathfrak{n}}
        + \varepsilon_\text{approx}^{\mathfrak{n}}
        + \mathcal{C}^{\mathfrak{n}} \left( \frac {\Delta t}{n_\tau} \right)^{p+1}
        + \mathcal{O}\bigg(
        \left( \frac {\Delta t}{n_\tau} \right)^{2p+2}\bigg)
        \right)
        \mathcal{D}^{\mathfrak{n}}
        + \mathcal{D}^{1} \mathcal{E}^{0},
    \end{equation*}
    where we have set
    \begin{equation*}
        \mathcal{D}^{\mathfrak{n}} = \prod_{m=\mathfrak{n}+1}^{n} \frac {1}{d_\mathcal{X}^{m}}.
    \end{equation*}
    Noting that $
        \mathcal{E}^{0} \leqslant
        \varepsilon_\text{int}^{0}
        + \varepsilon_\text{opt}^{0}
        + \varepsilon_\text{approx}^{0}
    $,
    we finally obtain
    \begin{equation*}
        \mathcal{E}^{n}
        \leqslant
        \sum_{\mathfrak{n}=0}^{n}
        \left(
        \varepsilon_\text{int}^{\mathfrak{n}}
        + \varepsilon_\text{opt}^{\mathfrak{n}}
        + \varepsilon_\text{approx}^{\mathfrak{n}}
        \right)
        \mathcal{D}^{\mathfrak{n}}
        +
        \biggl(
        \sum_{\mathfrak{n}=1}^{n}
        \mathcal{C}^{\mathfrak{n}}
        \mathcal{D}^{\mathfrak{n}}
        \biggr)
        \left( \frac {\Delta t}{n_\tau} \right)^{p+1} + \mathcal{O}\bigg(
        \left( \frac {\Delta t}{n_\tau} \right)^{2p+2}\bigg),
    \end{equation*}
    which is the desired result.
\end{proof}

\begin{corollary}
    \label{cor:NGSL_error}
    Under the hypotheses of \cref{thm:NGSL_error},
    assume that the velocity field $a$ is divergence-free.
    Let
    $C
        =
        \max_{\mathfrak{n} \in \{ 1, \dots, n \}}
        \mathcal{C}^{\mathfrak{n}}
    $.
    Then
    \begin{equation*}
        \mathcal{E}^{n+1}
        \leqslant
        \sum_{\mathfrak{n}=0}^{n} \varepsilon_\text{int}^{\mathfrak{n}}
        + \sum_{\mathfrak{n}=0}^{n} \varepsilon_\text{opt}^{\mathfrak{n}}
        + \sum_{\mathfrak{n}=0}^{n} \varepsilon_\text{approx}^{\mathfrak{n}}
        + \frac{C T}{n_\tau} \left( \frac {\Delta t}{n_\tau} \right)^{p}.
    \end{equation*}
\end{corollary}

\begin{proof}
    If the velocity field $a$ is divergence-free,
    then invoking Liouville's theorem,
    the determinant of the Jacobian matrix of the characteristic curve
    is constant and equal to one.
    Therefore, for all $n$, $\mathcal{D}^{n} = 1$.
    Plugging this into \cref{thm:NGSL_error} gives the desired result.
\end{proof}

Compared to classical error estimates,
the interpolation error in the semi-Lagrangian schemes (see e.g.~\cite{falcone1998convergence,ChaDesMeh2013,besse2008convergence,ferretti2020stability}) or the projection error in the semi-Lagrangian discontinuous Galerkin schemes (see e.g.~\cite{restelli2006semi,rossmanith2011positivity,qiu2011positivity}) is replaced with the approximation error,
and we have additional optimization and integration errors
(which are already present in semi-Lagrangian Galerkin methods).
\ra{The characteristic error,
    depending on the accuracy of the characteristic curve solver,
    is also present in classical schemes.
    Let us discuss these four errors in more detail.}
\begin{itemize}[nosep]
    \item The integration error
          depends on the chosen quadrature method.
          In our case, since we choose the Monte-Carlo method,
          it is bounded by $1 / \sqrt{N_c}$,
          where $N_c$ is the number of collocation points,
          up to a constant independent of $N_c$.
          Thus, this error can be made as small as desired.
          In addition, we can further simplify~\eqref{eq:NGSL_error_estimate}
          by remarking that
          \begin{equation*}
              \sum_{\mathfrak{n}=1}^{n} \varepsilon_\text{int}^{\mathfrak{n}}
              \leqslant
              \frac{C_\text{int} T} {\dt \sqrt{N_c}},
          \end{equation*}
          where $C_\text{int}$ is a constant independent of $\dt$ and $N_c$.
    \item The optimization error
          is the threshold at which the optimization algorithm is stopped.
          Depending on the algorithm used,
          this threshold may or may not be computationally reachable.
          Moreover, it may strongly depend on time
          since the optimization problem may be harder to solve
          at some time steps than at others.
    \item The approximation error
          depends only on the expressiveness
          of the approximation space
          \smash{$u_{\theta^{n+1}_{\ast}} \in \mathcal{V}_N^\text{neural,\smash{$t$}}$}.
          To decrease this error, we must increase the number of dofs
          (i.e., the number of neurons in the neural network),
          but no generic quantitative estimates
          are available for the moment. Note also that, if the approximation error is supposed to be bounded by $\varepsilon_N$, then the corresponding error satisfies
          \begin{equation*}
              \sum_{\mathfrak{n}=1}^{n} \varepsilon_\text{approx}^{\mathfrak{n}}
              \leqslant
              \frac{T \varepsilon_N} {\dt}.
          \end{equation*}
          In classical semi-Lagrangian schemes, $\varepsilon_N$ is generally replaced with an approximation error of the form $\mathcal{O}(\Delta x^{q+1})$, where $q$ is the spatial order of approximation. Here, on the contrary, we might expect that $\varepsilon_N$ is very small compared with $\Delta t$.
    \item \ra{The characteristic error directly depends on
              $\Delta t$ and on the characteristic curve solver.
              As such, it is the easiest to control,
              either by improving the order of the solver,
              or by adding more sub-time steps.
              Let us not that, as soon as an exact expression of the
              characteristic curve is known,
              $M_p$ is equal to zero,
              and so $C_n$ vanishes for all $n \geqslant 0$,
              leading to the whole characteristic error vanishing,
              as expected.}
\end{itemize}
All in all, supposing that these three errors can be controlled,
they must be at least in $O(\Delta t^2)$ to ensure the convergence of the method.
In pratice, we do not have fine control over these three errors.
If these errors become large,
larger times steps are recommended
to limit the error accumulating at each time step.
Note, however, that it would also lead to more difficult optimization problems.
This is expanded upon in \cref{rem:NGSL_time_step}.

\begin{remark}
    \label{rem:NGSL_time_step}
    \cref{thm:NGSL_error} might suggest that, in the case of pure advection, more time steps would lead to errors accumulating, and therefore one should take a time step as large as possible. But things are not that simple: the most computationally expensive step of our algorithm is to find a $\vartheta$ that allows $u_{\vartheta}$ (the target) to best fit $u_{\theta^n}(\mathcal{X}(t^n;t^{n+1}, \cdot))$ (the source). To achieve this, we perform a gradient descent starting from $\vartheta=\theta^n$. When the time step $\Delta t = t^{n+1} - t^n$ is very small, the characteristic curve $\mathcal{X}(t^n;t^{n+1},\cdot)$ is very short, so the source and the target are very close: this implies that the gradient descent is fast and leads to a very good fit. Conversely, if we take an extremely large $\Delta t$, the source and target are very far apart, and the minimization problem then becomes difficult. For problems where the complexity of the target solution increases with time (such as the Vlasov equation, which leads to filamentation), the fitting between a ``simple'' source and a ``very complex'' target could be technically very challenging. The choice of the time step is therefore necessarily a compromise. In the case of advection-diffusion, as explained in the next section, the error on the diffusion scales as $\sqrt{\Delta t}$. Therefore, it outweighs the fitting errors of the neural network, and the time step should be taken as small as possible.
\end{remark}

\subsection{Extension to advection-diffusion equations}
\label{sec:advection_diffusion}

Equipped with the NSL method applied to advection equations
with a space-, time- and parameter-dependent
advection field,
we extend it to the full advection-diffusion equation
\eqref{eq:advection_diffusion}.
In this context, recall that $\sigma$
represents a constant diffusion coefficient.

To perform this extension,
it turns out that it is sufficient to modify
the computation of the foot of the characteristic curves.
Indeed, after e.g.~\cite{Fer2010},
the advection-diffusion equation can be reformulated
as a stochastic differential equation (SDE)
with a Gaussian noise,
whose variance depends on the dimension $d$,
the time step~$\dt$,
and the diffusion coefficient~$\sigma$.
\rc{To discretize this SDE, an appropriate discretization
    of the Gaussian noise must be provided,
    see e.g.~\cite{HamMehSelSon2016,BonFerRoc2018}.
    For instance, on a 2D Cartesian mesh,
    it is shown in~\cite{BonCalCarFer2021} that
    4 directions are required for a scheme of order $1$
    and 9 directions for a scheme of order $2$.
    The first-order scheme is easily extended
    to higher dimensions,
    and thus requires $2d$ directions in dimension $d$.}

\rc{In this work, we follow the same approach,
    and define the $2d$ diffusion directions by
    \begin{equation*}
        \forall i \in \{1, \dots, d\}, \qquad
        v_i = e_i
        \text{\quad and \quad}
        v_{d+i} = -e_i,
    \end{equation*}
    where for all $i \in \{1, \dots, d\}$,
    $e_i$ is the \smash{$i^\text{th}$} vector
    of the canonical basis of $\mathbb{R}^d$.}

\rc{Equipped with the $2d$ directions,
    we can now compute the foot of the characteristic curves
    in the presence of a constant diffusion,
    with a scheme of time order $1$.
    In this case, the nonlinear optimization problem
    in step~5 of \cref{alg:NSL_known_characteristic}
    is replaced with
    \begin{equation*}
        \label{eq:NGSL_nonlinear_optimization_diffusion}
        \theta^{n+1} \in \argmin_{\vartheta \in \Theta}
        \sum_{k=1}^{N_c} \,
        \left|
        u_{\vartheta}(x_k, \mu_k) -
        \sum_{i=1}^{2d}
        u_{\theta^n}(\widetilde{\mathcal{X}}_i(t^n; t^{n+1}, x_k, \mu_k))
        \right|^2.
    \end{equation*}
    Compared to the pure advection case,
    this problem contains $2d$ characteristic feet
    \smash{$(\widetilde{\mathcal{X}}_i)_{i \in \{1, \dots, 2d\}}$}
    to be computed instead of one.
    They are calculated, using the $2d$ directions, as follows:
    \begin{equation*}
        \forall i \in \{1, \dots, 2d\}, \quad
        \widetilde{\mathcal{X}}_i(t^n; t^{n+1}, x, \mu)
        = \widetilde{\mathcal{X}}(t^n; t^{n+1}, x, \mu)
        + \sqrt{2 d \, \sigma \dt} \, v_i,
    \end{equation*}
    where \smash{$\widetilde{\mathcal{X}}(t^n; t^{n+1}, x)$}
    is the result of \cref{alg:characteristic_curve},
    i.e., the approximate foot of the characteristic
    curve in the presence of advection only.}

\rc{We remark that our approach is compatible with any already-developed strategy to improve the approximation of diffusion in stochastic differential equations or semi-Lagrangian schemes. For instance, using higher-order spherical quadrature rules could provide finer and more uniform angular coverage than the canonical basis of $\mathbb{R}^d$, see e.g.~\cite{SarSol2019}. Similarly, adaptive strategies in which directions are dynamically aligned with the local eigenspectrum of the diffusion tensor could better capture anisotropic effects, see e.g.~\cite{BokSim2016}.
    These variants could be integrated into our method without altering its general structure.}

\subsection{Further improvements}
\label{sec:improving_accuracy}

In this section,
motivated by the error estimate from \cref{thm:NGSL_error},
we present several techniques
to further improve the results of the NSL method.
Namely, we first focus on improving the resolution
of the optimization problem
\eqref{eq:NGSL_nonlinear_optimization}
in \cref{sec:natural_gradient_preconditioning}
via natural gradient preconditioning,
to decrease the optimization error
$\varepsilon_\text{opt}$.
Then, we introduce an adaptive sampling strategy
in \cref{sec:adaptive_sampling} that,
while not directly related to the NSL method,
allows us to improve the quality of the solution
by improving the Monte-Carlo approximation
of the integral in~\eqref{eq:NGSL_nonlinear_optimization},
thus decreasing the integration error
$\varepsilon_\text{int}$.
\ra{Finally, we briefly mention in
    \cref{sec:boundary_conditions}
    how to tackle boundary conditions.}

\subsubsection{Natural gradient preconditioning}
\label{sec:natural_gradient_preconditioning}

To properly introduce natural gradient descent,
we go back to the optimization problem
\eqref{eq:NGSL_nonlinear_optimization}
and rewrite it,
dropping the time indices for simplicity, as
\begin{equation}
    \label{eq:NGSL_nonlinear_optimization_rewritten}
    \theta \in \argmin_{\vartheta \in \Theta}
    \int_{\Omega} \bigl|
    \mathcal{N}(x, \vartheta) - u^\star(x)
    \bigr|^2 \, dx.
\end{equation}
In \eqref{eq:NGSL_nonlinear_optimization_rewritten},
we have reused the nonlinear parametric
function $\mathcal{N}$ defined
in \eqref{eq:neural_network},
typically a neural network,
and we have introduced the function
$u^\star$ to concisely represent
the approximate foot of the characteristic curve.
It should be noted that
\eqref{eq:NGSL_nonlinear_optimization_rewritten}
is nothing but a nonlinear least-squares problem
with objective function $u^\star$.

Natural gradient descent
was introduced in~\cite{AmaDou1998}.
Now, we present the method,
alongside its main advantages
over standard gradient descent,
and why it is particularly well-suited
to solving problems such as
\eqref{eq:NGSL_nonlinear_optimization_rewritten}.

To solve \eqref{eq:NGSL_nonlinear_optimization_rewritten},
one usually constructs a sequence of
iterates $(\vartheta^n)_n$,
which converges to the solution $\theta$.
These iterates are updated using the
gradient descent method\footnote{Usually, a stochastic version of the gradient method is used; the deterministic version is presented here for simplicity.}.
This method amounts to finding a so-called
descent direction $\eta$,
and writing $(\vartheta^n)_n$ as the time discretization
of the ordinary differential equation
$\dot{\vartheta} = \eta$.
To find the descent direction,
it is convenient to rewrite the optimization problem
\eqref{eq:NGSL_nonlinear_optimization_rewritten} as
\begin{equation*}
    \theta \in \argmin_{\vartheta \in \Theta}
    \mathfrak{l}(\mathfrak{n}(\vartheta)),
\end{equation*}
with $\mathfrak{l}$ the loss function
depending on the neural network $\mathfrak{n}$,
which itself depends on the dofs $\vartheta$.
More precisely, we have
\begin{equation}
    \label{eq:definitions_l_and_n}
    \begin{aligned}
        \mathfrak{l}: \;
        \mathcal{H} & \to \mathbb{R}               \\
        f           & \mapsto \int_{\Omega} \bigl|
        f(x) - u^\star(x)
        \bigr|^2 \, dx
    \end{aligned}
    \text{\qquad and \qquad}
    \begin{aligned}
        \mathfrak{n}: \;
        \mathbb{R}^N & \to \mathcal{H}                                                         \\
        \vartheta    & \mapsto (x \mapsto \mathcal{N}(x, \vartheta)). \vphantom{\int_{\Omega}}
    \end{aligned}
\end{equation}
In \eqref{eq:definitions_l_and_n},
similarly to \eqref{eq:NGSL_nonlinear_optimization_in_error_estimate},
we have defined the finite-dimensional
nonlinear function space
\begin{equation*}
    \mathcal{H} = \{
    x \in \Omega
    \mapsto
    \mathcal{N}(x, \vartheta) \in \mathbb{R},
    \;\; \vartheta \in \Theta
    \}
    \subset L^2(\Omega).
\end{equation*}
Equipped with this notation,
the gradient descent (GD) methods
corresponds to the steepest descent
in terms of the vector~$\vartheta$,
i.e., the descent direction is given by
\begin{equation*}
    \eta^\text{GD}
    =
    - \nabla_\vartheta \, \mathfrak{l}(\mathfrak{n}(\vartheta)),
\end{equation*}
where $\nabla_\vartheta \, \mathfrak{l}(\mathfrak{n}(\vartheta))$
is the gradient of $\vartheta \mapsto \mathfrak{l}(\mathfrak{n}(\vartheta))$.
This disregards the fact that $\mathcal{H}$
is a nonlinear space.

Conversely, natural gradient descent (NGD)
corresponds to an approximation of the steepest
descent direction in terms of
the nonlinear function $\mathfrak{n}$.
In the general case, one needs to consider the
tangent space of the manifold~$\mathcal{H}$
at the point $\mathfrak{n}(\vartheta)$.
However, in our specific case,
the system under consideration
is the nonlinear least-squares problem
\eqref{eq:NGSL_nonlinear_optimization_rewritten}.
Thus, NGD amounts to
modifying the gradient descent direction to
\begin{equation*}
    \eta^\text{NGD}
    =
    G(\vartheta)^\dagger \, \eta^\text{GD},
\end{equation*}
where $G(\vartheta)$ is the so-called
\emph{Fisher information matrix},
with $G(\vartheta)^\dagger$ its pseudo-inverse.
Note that this matrix is identical
to the Neural Galerkin mass matrix,
defined in \eqref{eq:neural_Galerkin_M_L}.
The $(i,j)$ coefficient of this matrix is defined by
\begin{equation*}
    G_{i,j}(\vartheta)
    =
    \left\langle
    \partial_{\vartheta_i} \mathfrak{n}(\vartheta),\,
    \partial_{\vartheta_j} \mathfrak{n}(\vartheta)
    \right\rangle_{L^2(\Omega)}
    =
    \int_\Omega
    \partial_{\vartheta_i} \mathcal{N}(x, \vartheta) \,\,
    \partial_{\vartheta_j} \mathcal{N}(x, \vartheta) \,\,
    dx.
\end{equation*}
In the end, NGD amounts to preconditioning
the gradient descent direction
using the gradient of the neural network $\mathcal{N}$
with respect to its dofs.

More sophisticated methods have been proposed
to further improve NGD, most notably
in the context of physics-informed learning,
see for instance~\cite{NurLeiYan2023,MueZei2023,SchFur2025}.
In all of these cases, information from the loss function
is used to improve the convergence of the optimization
problem.
In physics-informed learning,
the loss function involves the residual of the PDE,
like \eqref{eq:PINN_optimization_problem} for PINNs.
As a result, the preconditioner
becomes much more complicated to compute,
since derivatives with respect to $\vartheta$ of the
space and time derivatives of the residual become involved.
One major advantage of the proposed method
is that the optimization
problem~\eqref{eq:NGSL_nonlinear_optimization}
is nothing but a nonlinear least-squares problem,
for which classical natural gradient preconditioning
can directly be applied with minimal
implementation effort and computational cost.

\subsubsection{Adaptive sampling}
\label{sec:adaptive_sampling}

We now introduce a technique to
decrease the integration error $\varepsilon_\text{int}$.
Namely, rather than uniformly sampling
the collocation points in step 7 of
\cref{alg:NSL_known_characteristic},
we propose an adaptive sampling strategy,
in the framework of rejection sampling.
For simplicity, we assume that we wish to
sample additional points in the vicinity
of the zeros of some function $f: \Omega \to \mathbb{R}$.
This is motivated by e.g.
the transport of level-set functions,
or the approximation of localized functions
(where $f$ could be
the inverse of the norm of the gradient of the function
to be approximated).

\cref{alg:adaptive_sampling}
details the adaptive sampling strategy.
It relies on three hyperparameters
$\sigma_1$, $\sigma_2$ and $\sigma_3$,
whose specific choice turns out not to be crucial.
The algorithm performs three iterations,
each corresponding to a different value of the parameter
$\sigma \in \{\sigma_1, \sigma_2, \sigma_3\}$.
In each iteration, it uniformly samples $10 N_c$ candidate points from the domain
$\Omega \times \mathbb{M}$,
evaluates the function $f$ at those points,
and computes weights $\omega_k$
that reflect proximity to the zero set of $f$.
From these, it selects up to $N_c/4$
points where $\omega_k > 0.75$,
which are close to the zeros of $f$,
and adds them to the growing set of collocation points.
After the three passes,
it supplements the set with up to
$N_c$ additional uniformly sampled points
to ensure some coverage of the entire domain.
The final output is a collection of
$N_c$ points in $\Omega \times \mathbb{M}$,
strategically biased towards regions where
$f$ is close to $0$.

\begin{algorithm}[!ht]
    \caption{Adaptively sampling $N_c$ points concentrated around the zeros of $f$}
    \label{alg:adaptive_sampling}
    \begin{algorithmic}[1]
        \State \textbf{Input:} Number $N_c$ of collocation points, function $f$, parameters $\sigma_1$, $\sigma_2$, $\sigma_3$
        \State \textbf{Output:} $N_c$ points $(x_k, \mu_k)_{k \in \{1, \dots, N_c\}} \in \Omega^{N_c} \times \mathbb{M}^{N_c}$

        \medskip

        \For{$\sigma \in \{\sigma_1, \sigma_2, \sigma_3\}$}

        \State Uniformly sample $10 N_c$ points in $\Omega \times \mathbb{M}$, denoted by $(\bar x_k, \bar \mu_k)_{k \in \{1, \dots, N_c\}}$

        \State Compute $\omega_k = \exp(-f(\bar x_k, \bar \mu_k)^2 / (2 \sigma^2))$ for all $k \in \{1, \dots, 10 N_c\}$

        \State Add at most $N_c / 4$ random points $(\bar x_k, \bar \mu_k)$ such that $\omega_k > 0.75$ to the set of collocation points $(x_k, \mu_k)$

        \EndFor

        \medskip

        \State Add at most $N_c$ uniformly sampled points to the set of collocation points $(x_k, \mu_k)$

        \medskip

        \State \textbf{Return:} $(x_k, \mu_k)_{k \in \{1, \dots, N_c\}}$
    \end{algorithmic}
\end{algorithm}

Let us emphasize that more sophisticated adaptive sampling
could be used to further decrease $\varepsilon_\text{int}$.
The proposed strategy was chosen for its simplicity.
One idea could be to advect the collocation points
using the characteristic curves.
However, since the optimization problem
\eqref{eq:NGSL_nonlinear_optimization}
is a simple nonlinear least-squares problem
and does not involve the residual of the PDE
(compared to e.g.\ discrete PINNs),
we would not need more complex physics-informed sampling
relying on the PDE residual,
see e.g.\ the review paper~\cite{WuZhuTanKarLu2023}.
In our case,
another possibility would be to use optimal least-squares
sampling~\cite{CohMig2017}.

\subsubsection{Boundary conditions}
\label{sec:boundary_conditions}

{\raparagraph In neural methods such as PINNs, we distinguish two main approaches to handling boundary conditions.

\paragraph{Weak boundary conditions}
In the so-called weak approach,
initially proposed in~\cite{RAISSI2019686},
a loss term enforcing the boundary conditions
is added to the loss function.
The same strategy is applicable to our neural semi-Lagrangian framework. In \eqref{eq:advection_diffusion}, the boundary condition is given by $u=g$ on $\partial \Omega$. Then, at each step, we solve:
\begin{equation*}
    \theta^{n+1} \in \argmin_{\vartheta \in \Theta}
    \int_{\Omega} \bigl|
    u_{\vartheta}(x) - u_{\theta^n}(\mathcal{X}(t^n; t^{n+1}, x))
    \bigr|^2 \, dx
    +
    \omega^\text{BC} \int_{\partial \Omega} \bigl|
    B(u_{\vartheta}(x)) - g(x)
    \bigr|^2 \, dx
\end{equation*}
instead of \eqref{eq:NGSL_nonlinear_optimization},
where $\omega^\text{BC}$ is a hyperparameter that
controls the weight of the boundary condition loss term.
Neumann and periodic boundary conditions are handled
in a similar way.

\paragraph{Strong boundary conditions}
The strong approach, generally preferable, consists in enforcing the boundary conditions directly within the model. It was originally proposed in \cite{LagLikFot1998}, and it is the one we used in the paper. It can be done by using a representation of the domain with a level-set function, see~\cite{SUKUMAR2022114333}, or with eigenmodes of the Laplace operator associated with the domain, see~\cite{KAST2024112986}. For a Dirichlet boundary condition ($u=g$ on $\partial \Omega$), given a level-set function~$\Phi$ for the domain, the boundary condition is imposed by replacing $u_{\theta}$ with
$
    \Phi u_{\theta} + g.
$
This formulation clearly equals $g$ on the boundary.
Neumann boundary conditions are treated similarly;
the reader is referred to~\cite{SUKUMAR2022114333}
for a detailed breakdown.
To impose periodic boundary conditions
for rectangular domains in the NSL method,
we simply place the foot of the characteristic curve
back into the domain $\Omega$ when it leaves through
a periodic boundary.
This amounts to adding a step between steps 14 and 15 of
\cref{alg:characteristic_curve},
replacing $\mathcal{X}$
with~$x_{\text{min}} + (\mathcal{X} - x_{\text{min}}) \% X$,
where $x_{\text{min}}$
is the lower bound of the domain $\Omega$,
$X$ is its extent, and $\%$ denotes the modulo operator.
For example,
in the case of a domain $\Omega = [0, 1] \times [-1, 1]$
with periodic boundary conditions,
$x_\text{min} = (0, -1)$ and $X = (1, 2)$.
}

\section{Validation}
\label{sec:applications}

This section is dedicated to the presentation of some numerical experiments,
to validate the proposed method.

We first test the method on an advection equation
with constant advection coefficient
in one space dimension in \cref{sec:1D_transport_constant}.
Then, we move on to nonconstant advection coefficients
in a multidimensional setting in \cref{sec:nD_transport}.
In these cases, we compare the proposed Neural Semi-Lagrangian (NSL) method
from \cref{sec:NGSL}
to other neural methods,
namely PINNs from \cref{sec:PINNs},
discrete PINNs (dPINNs) from \cref{sec:neural_Galerkin},
and the neural Galerkin~(NG) method,
also from \cref{sec:neural_Galerkin}.
\rall{In all these cases, we deliberately avoid showing
    the variability of the results with respect
    to the random initialization and optimization of the neural networks.
    Indeed, the method has a very low variability,
    and showing it would only clutter the figures.}

Afterwards, in \cref{sec:3D_cylinder},
we test the several avenues of improvement
proposed in \cref{sec:improving_accuracy}.
The improved NSL scheme is then used
on challenging level-set transport problems
in \cref{sec:level_sets}.
Lastly, the scheme is tested on high-dimensional
advection-diffusion problems in
\cref{sec:advection_diffusion_experiments}.

Unless otherwise mentioned,
all ODEs (in e.g.\ dPINN and NG,
or in the NSL method when
the foot of the characteristic curve is unknown)
are solved using the Runge-Kutta 4 (RK4) method.
Moreover, all nonlinear optimization problems
in \cref{sec:1D_transport_constant,sec:nD_transport}
are solved by first applying the Adam optimizer,
and then switching to LBFGS for the last $10\%$ of the epochs.
Afterwards, natural gradient preconditioning is applied.
The problems are solved using the \texttt{PyTorch} library~\cite{PyTorch2}
and the \texttt{ScimBa}\footnote{\url{https://gitlab.inria.fr/scimba/scimba}}
scientific machine learning library;
to report computation times,
we use a single AMD Instinct MI210 GPU.

\subsection{Constant advection in 1D}
\label{sec:1D_transport_constant}

The first numerical experiment we run consists in solving
\begin{equation}
    \label{eq:advection_1D}
    \begin{dcases}
        \partial_t u + a \partial_x u = 0,
         & \qquad x \in (0, 2), \; t \in (0, 1), \, \mu \in \mathbb{M}, \\
        u(0, x, \mu) = u_0(x, \mu),
         & \qquad x \in (0, 2), \, \mu \in \mathbb{M},                  \\
        u(t, 0, \mu) = u(t, 2, \mu),
         & \qquad t \in (0, 1), \, \mu \in \mathbb{M},
    \end{dcases}
\end{equation}
where we have introduced two parameters: the variance $\nu$ of the Gaussian pulse,
and the advection coefficient~$a$.
The two parameters are grouped into a single vector $\mu = (\nu, a) \in \mathbb{M}$,
and we set the parameter space to be
$\mathbb{M} = [0.05, 0.15] \times [0.5, 1] \subset \mathbb{R}^2$.
The parametric initial condition is defined by
\begin{equation*}
    u_0(x, \mu) = \exp\left(-\frac{(x - 0.5)^2}{2 \nu^2}\right),
\end{equation*}
and the exact solution is given by
$u_\text{ex}(t, x, \mu) = u_0(x - a t, \mu)$.
Since $a$ is constant in space,
the first branch of \cref{alg:characteristic_curve}
is applied.

\subsubsection{Non-parametric case}
\label{sec:1D_transport_constant_non_parametric}

We first focus on the non-parametric case,
where the advection coefficient is fixed to $a = 1$
in \eqref{eq:advection_1D}.
The variance $v$ is also set to a constant value.
We compare the PINN, dPINN, NG and NSL methods.
All hyperparameters are given in
\cref{tab:hyperparameters_1D_transport_constant_non_parametric}.
Note that the PINN solves a 2D problem
(1D in space and 1D in time),
while the neural network inherent in dPINN, NG and NLS
solves a 1D problem in space.
For comparison purposes,
we denote the pointwise error by $e_x$
and the $L^2$ error by $e_t$.
They are respectively given by
\begin{equation*}
    e_x(x) = u(1, x) - u_N(1, x)
    \text{\qquad and \qquad}
    e_t(t) = \int_0^2 \bigl| u(t, x) - u_N(t, x) \bigr|^2 dx.
\end{equation*}
\ra{They are computed using analytic solutions
    and a Monte-Carlo estimate of the $L^2$
    error using many more points than during
    the training process.}

\begin{figure}[!ht]
    \centering
    \begin{subfigure}{\textwidth}
        \centering
        \includegraphics{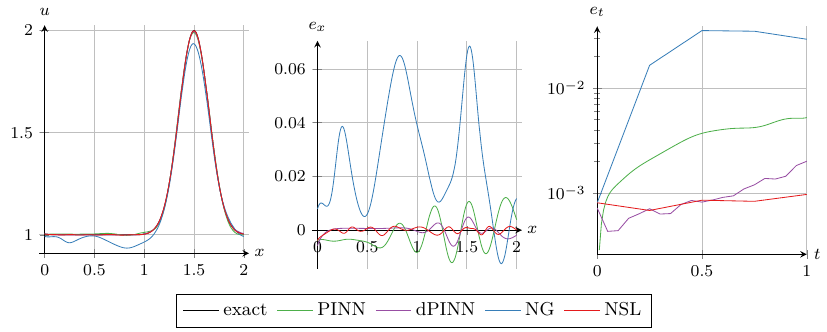}
        \caption{1D constant advection: $\nu = 0.15$.}
        \label{fig:1D_constant_advection_0.15}
    \end{subfigure}
    \begin{subfigure}{\textwidth}
        \centering
        \includegraphics{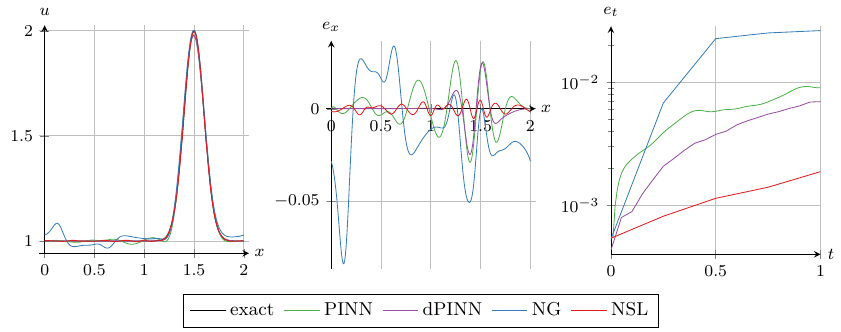}
        \caption{1D constant advection: $\nu = 0.1$.}
        \label{fig:1D_constant_advection_0.1}
    \end{subfigure}
    \begin{subfigure}{\textwidth}
        \centering
        \includegraphics{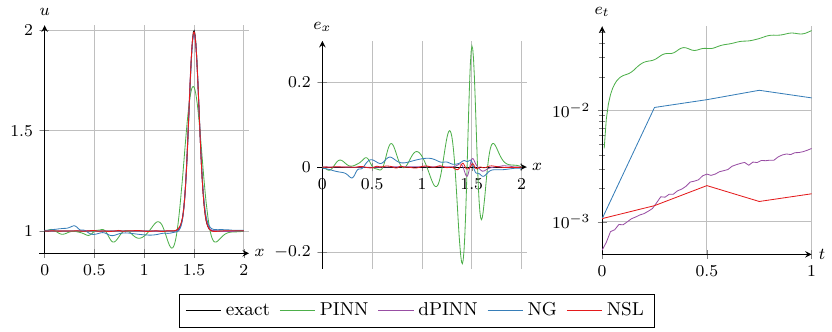}
        \caption{1D constant advection: $\nu = 0.05$.}
        \label{fig:1D_constant_advection_0.05}
    \end{subfigure}
    \caption{1D constant advection from \cref{sec:1D_transport_constant_non_parametric}: comparison of the different methods for several values of the variance $\nu$, and for a fixed time step ($\dt = 0.05$, i.e., $20$ time steps in the dPINN method; $\dt = 0.25$, i.e., $4$ time steps in the NG and SL methods). From left to right: prediction of the solution at $t=1$, pointwise errors between the predicted solutions and the exact solution at $t=1$, and time evolution of the $L^2$ error.}
    \label{fig:1D_constant_advection}
\end{figure}

\cref{fig:1D_constant_advection} displays the results of the four approaches,
for four distinct values of the variance $v$.
Note that each value of the variance corresponds
to a new problem to solve,
with a new neural network to train.
The truly parametric case is treated in the following section.
For dPINN, we take $\dt = 0.05$
(corresponding to $20$ time steps)
for $\nu \geqslant 0.1$, and $\dt = 0.02$
(corresponding to $50$ time steps) otherwise.
For NG and NSL, we take $\dt = 0.25$
(corresponding to $4$ time steps).
We have chosen $150$ epochs for the inner optimization
problems \eqref{eq:NGSL_nonlinear_optimization}
of the NSL method, since the problem converged quickly
and more epochs did not bring noticeable improvements.

We observe that the NSL method outperforms the other methods,
by almost an order of magnitude in some cases.
Indeed, the NSL method is able to capture
both the peak of the Gaussian pulse
and the constant state far away from the peak,
with good accuracy.

On the contrary, the other methods tend to either diffuse the peak
or create oscillations in the constant state, or both.
Notably, the dPINN method required quite a restrictive
time step to perform well;
larger time steps led to numerical instabilities
and oscillations destroying the simulation quality.
Moreover, it also required better solving the optimization problems,
which led to an increased number of epochs
during the initialization phase
and the inner epochs, compared to the NG and NSL methods.
Therefore, even if the results of the dPINN method
are quite good, they come at a large computational cost,
at least $10$ times as large as the NG and NSL methods.

\begin{figure}[!ht]
    \centering
    \includegraphics{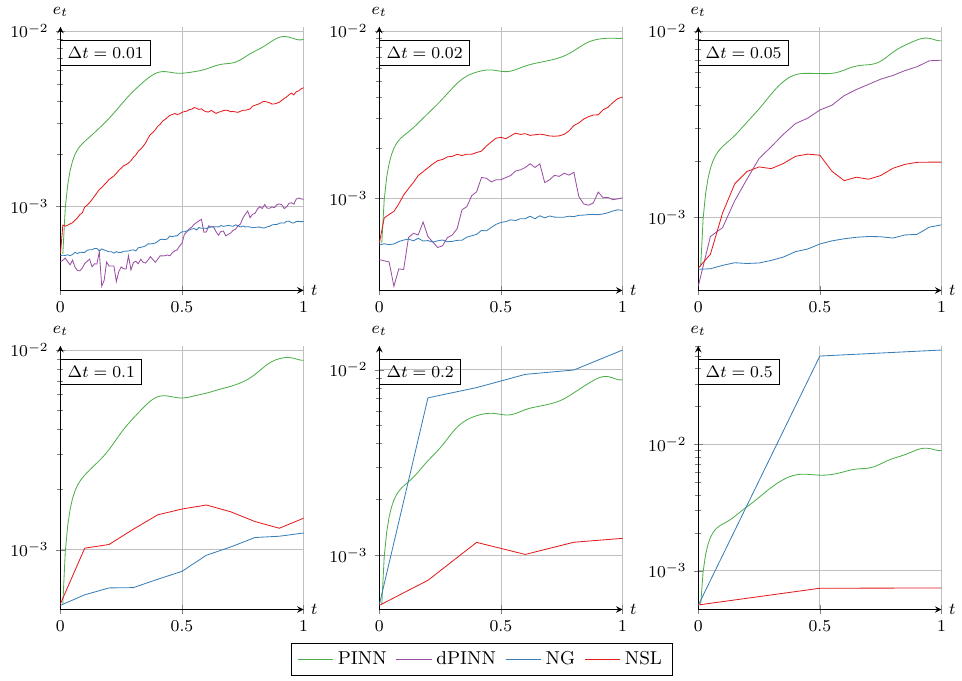}
    \caption{1D constant advection from \cref{sec:1D_transport_constant_non_parametric}: comparison of the different methods for several values of the time step $\dt$ in the NG and SL methods and a fixed variance, equal to $0.1$. From top to bottom and left to right, we take $\dt \in \{0.01, 0.02, 0.05, 0.1, 0.2, 0.5\}$, and we display the time evolution of the $L^2$ error $e_t$.}
    \label{fig:1D_constant_advection_wrt_dt}
\end{figure}

To better understand the role of $\dt$,
we display in \cref{fig:1D_constant_advection_wrt_dt}
the results of the methods for $\nu=0.1$ and for several values of $\dt$.
Namely, we show the $L^2$ error $e_t$ as a function of time,
for $\dt \in \{0.01, 0.02, 0.05, 0.1, 0.2, 0.5\}$,
which corresponds to $100$, $50$, $20$, $10$, $5$, and $2$ time steps, respectively.
For time steps $\dt \leqslant 0.05$,
we perform $300$ epochs for the inner optimization problems
instead of $150$, to help convergence.
The results of the dPINN method are only shown
for $\dt \leqslant 0.05$, since otherwise it does not converge.

On the one hand, for both the NG and dPINN methods,
we observe that the accuracy is improved
when using smaller time steps,
overtaking the PINN for $\dt \leqslant 0.1$.
On the other hand,
interestingly, we note that smaller time steps do not necessarily lead to better results for the NSL method,
at least in this simpler case.
\ra{Indeed, the exact characteristic curves are known for constant advection,
    and so the error due to the characteristic solver vanishes.
    We are thus left with only the integration, optimization, and approximation errors.
    In this simple 1D example, the integration and approximation errors remain very low,
    and so the whole error is driven by the optimization one.
    These optimization errors will accumulate during the time iterations,
    because each optimization problem is not solved exactly.
    This naturally means that more iterations will lead to a larger total error,
    unless they are also accompanied by a reduction of the optimization error at each iteration.
    In short, this expected behavior is observed since the optimization errors
    are much larger than the time errors;
    this was already the subject of \cref{rem:NGSL_time_step}.
    We expect this observation to change when
    considering time-dependent advection coefficients,
    or advection-diffusion equations,
    or even pure advection equations whose solutions exhibit strong space variations.}
Further, note that, in this case,
we have considered a fixed number $N_e$ of epochs for each time step.
Since, for smaller time steps, the optimization problems are easier to solve,
we could have reduced $N_e$ and obtained almost the same results
for a smaller computational cost.
However, this reduction is not enough to justify taking small time steps,
because it remains cheaper to solve a few hard optimization problems.

\subsubsection{Parametric case}
\label{sec:1D_transport_constant_parametric}

We now turn to the parametric version of
\eqref{eq:advection_1D},
where this time the variance $\nu$ and the advection coefficient~$a$
are both parameters.
The hyperparameters are defined in
\cref{tab:hyperparameters_1D_transport_constant_parametric}.
Compared to the previous case, the PINN is now solving a 4D problem
(1D in space, 1D in time, and 2D in the parameter space),
and the dPINN, NG and NSL neural networks are solving 3D problems
(1D in space and 2D in the parameter space).
This makes the problem quite a bit more complex,
which explains the increased number of epochs,
collocation points and neurons in the hidden layers.

An added difficulty for the dPINN, NG and NSL methods
is that the initial condition does not depend on~$a$,
so the network has to learn to ignore the second parameter
for learning the initial condition,
but to use it for the time evolution.
To overcome that difficulty, on the one hand,
we set a rather small time step for NG and dPINN
($\dt = 0.025$).
On the other hand,
the NSL time step is still large ($\dt = 0.25$),
but the number of epochs in the inner optimization problems
has been increased to $\num{2500}$.
With this setup, on the one hand, for NG and NSL,
it takes $15$ seconds to train
the network approximating the initial condition,
while the whole time stepping takes about $1$ minute
(meaning that NSL iterations are about $10$ times slower than NG ones).
\ra{This can be compared to the non-parametric case,
    where the NG and NSL methods took about $10$ seconds to train the initial condition
    and~$5$ seconds for the time-stepping.
    This is due to the larger network and
    higher number of epochs and collocation points.}
On the other hand, it takes about $5$ minutes to train the PINN.
The dPINN method is much slower than the other ones:
indeed, training and solving the optimization problems
totals about $20$ minutes.
This is due to the larger dimension of the problem,
requiring additional collocation points.

\begin{figure}[!ht]
    \centering
    \includegraphics{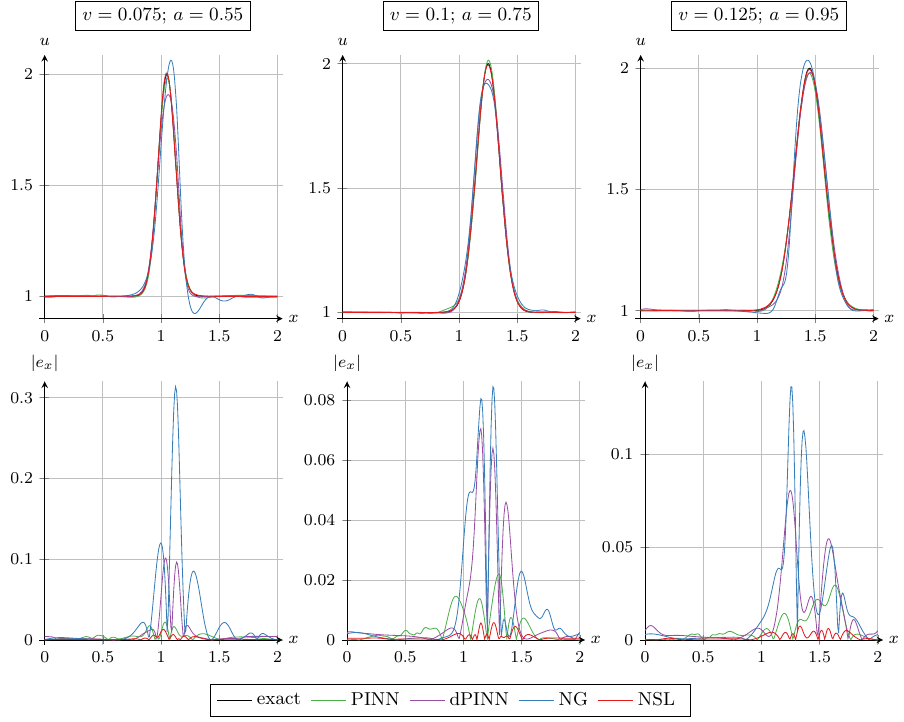}
    \caption{1D parametric advection from \cref{sec:1D_transport_constant_parametric}: comparison of the different methods for several values of the parameters and with a fixed time step. Top panels: approximate solutions; bottom panels: absolute values of the errors. From left to right: $\mu = (0.075, 0.55)$, $\mu = (0.1, 0.75)$, and $\mu = (0.125, 0.95)$.}
    \label{fig:1D_param_advection}
\end{figure}

As a first test, we take three values of the parameters
$\mu \in \mathbb{M}$,
and we draw the solutions (top panels)
and associated errors (bottom panels)
in \cref{fig:1D_param_advection}.
We observe that NSL consistently
yields a better approximation
than the other methods;
namely, it is less oscillatory.

To get a more quantitative estimation
of the error with respect to the parameter set $\mathbb{M}$,
we compute the~$L^2$ norm of $e_x$
(computed at the final time) for $\num{2000}$
randomly sampled parameters in $\mathbb{M}$.
The results are presented in \cref{tab:1D_param_advection},
where we report the minimum, average and maximum
of the~$L^2$ norm of $e_x$ over $\mathbb{M}$,
as well as its standard deviation.
We observe that, over the whole parameter set $\mathbb{M}$,
the NSL method consistently gives very good results.
On average, the NSL results
are quite close to the PINN.
However, this comes at a much lower computational cost,
as the NSL method is about
$5$ times faster than the PINN method.

\begin{table}[!ht]
    \centering
    \begin{tabular}{lcccc}
        \toprule
            & PINN            & dPINN           & NG              & NSL             \\
        \cmidrule(lr){1-5}
        min & $\num{2.92e-3}$ & $\num{4.47e-3}$ & $\num{1.63e-2}$ & $\num{5.03e-3}$ \\
        avg & $\num{9.00e-3}$ & $\num{2.69e-2}$ & $\num{5.31e-2}$ & $\num{1.08e-2}$ \\
        max & $\num{4.29e-2}$ & $\num{9.76e-2}$ & $\num{5.00e-1}$ & $\num{9.09e-2}$ \\
        std & $\num{5.78e-3}$ & $\num{1.36e-2}$ & $\num{1.14e-1}$ & $\num{5.10e-3}$ \\
        \bottomrule
    \end{tabular}
    \caption{1D parametric advection from \cref{sec:1D_transport_constant_parametric}: statistics on the errors $e_x$ for $\num{2000}$ randomly sampled parameters in $\mathbb{M}$.}
    \label{tab:1D_param_advection}
\end{table}

\subsection{Nonconstant advection in 2D}
\label{sec:nD_transport}

We now provide some numerical experiments
on advection equations with nonconstant advection coefficients,
in two and three space dimensions,
and with several parameters.
The first test case consists in
a 2D advection equation with a rotating transport,
presented in \cref{sec:2D_transport_rotating_parametric}.
The second one,
presented in \cref{sec:2D_fake_Vlasov}
is a 2D advection equation,
without parameters, mimicking the Vlasov equation
representing the movement of charged particles in an electric field.

In this section,
we do not consider the dPINN method any longer,
since it is not competitive with NSL and NG for the problems at hand,
especially in terms of stability condition,
and to simplify the presentation.

\subsubsection{2D parametric rotating transport}
\label{sec:2D_transport_rotating_parametric}

We first consider a 2D advection equation
with a nonconstant advection coefficient,
leading to a rotating advection equation
within the unit disk $\Omega = \mathbb{D}^1$.
The advection equation reads
\begin{equation*}
    \begin{dcases}
        \partial_t u + a(x) \cdot \nabla u = 0,
         & \quad x \in \Omega, \;
        t \in (0, 1), \;
        \mu \in \mathbb{M},                \\
        u(0, x, \mu) = u_0(x, \mu),
         & \quad x \in \Omega, \;
        \mu \in \mathbb{M},                \\
        u(t, x, \mu) = 0,
         & \quad x \in \partial \Omega, \;
        t \in (0, 1), \;
        \mu \in \mathbb{M},
    \end{dcases}
\end{equation*}
with the divergence-free advection vector field given by
\begin{equation}
    \label{eq:2D_rotating_advection_vector_field}
    a(x) = \begin{pmatrix}
        2 \pi x_2 \\ - 2 \pi x_1
    \end{pmatrix}.
\end{equation}
The parametric exact solution is given by
\begin{equation*}
    u_\text{ex}(t, x, \mu)
    =
    1 + \exp \left(
    -\frac {1} {2 v^2} \Big(
    \big(x_1 - c \, x_1^0(t)\big)^2 +
    \big(x_2 - c \, x_2^0(t)\big)^2
    \Big)
    \right),
\end{equation*}
where $x_1^0(t) = \cos(2 \pi t)$ and $x_2^0(t) = \sin(2 \pi t)$,
and with the parameters $\mu = (v, c)$,
where $\mu$ is in the parameter set
$\mathbb{M} = [0.05, 0.1] \times [0.2, 0.4]$.
Equipped with the exact solution,
the initial condition is then simply given by
$u_0(x, \mu) = u_\text{ex}(0, x, \mu)$.

Note that the advection field
\eqref{eq:2D_rotating_advection_vector_field}
is no longer constant in space,
although it still leads to an exact solution
to the ODE \eqref{eq:characteristic_curve}
governing the characteristic curves.
This exact solution is given by
\begin{equation*}
    \mathcal{X}_\text{ex}(s; t, x) =
    \begin{pmatrix}
        k_2(s; x) \cos(2 \pi t) +
        k_1(s; x) \sin(2 \pi t) \\
        k_1(s; x) \cos(2 \pi t) -
        k_2(s; x) \sin(2 \pi t)
    \end{pmatrix},
\end{equation*}
where we have set
\begin{equation*}
    k_1(s; x) = x_1 \sin(2 \pi s) + x_2 \cos(2 \pi s)
    \text{\qquad and \qquad}
    k_2(s; x) = x_1 \cos(2 \pi s) - x_2 \sin(2 \pi s).
\end{equation*}
This allows us to apply the second branch of
\cref{alg:characteristic_curve}
when computing the characteristic curves
in the NSL method.

The hyperparameters are given in
\cref{tab:hyperparameters_2D_transport_rotating_parametric}.
We note that, for this experiment,
more collocation points are used to approximate
the integrals.
Indeed, the solution is a Gaussian function,
highly localized in both parameter space and physical space,
which requires a better Monte-Carlo approximation.

We run the PINN, NG and NSL methods on this problem.
On the one hand, we set
a fixed time step $\dt = 0.02$ for the NG method,
leading to $50$ time steps
and $\dt = 0.5$ for the NSL method,
leading to $2$ time steps.
\ra{Since the characteristic curves are known exactly, like in
    \cref{sec:1D_transport_constant_non_parametric},
    we have chosen a small number of time steps
    to avoid the accumulation of optimization errors.}
With this setup, it takes about $10$ seconds to
train the network approximating the initial condition,
and about $20$ seconds to run
both the NG and the NSL methods.
On the other hand, since the PINN is now solving a 5D problem,
a larger network has to be used,
leading to a training time of about $5$ minutes.

\begin{figure}[!ht]
    \centering
    \includegraphics[width=\textwidth]{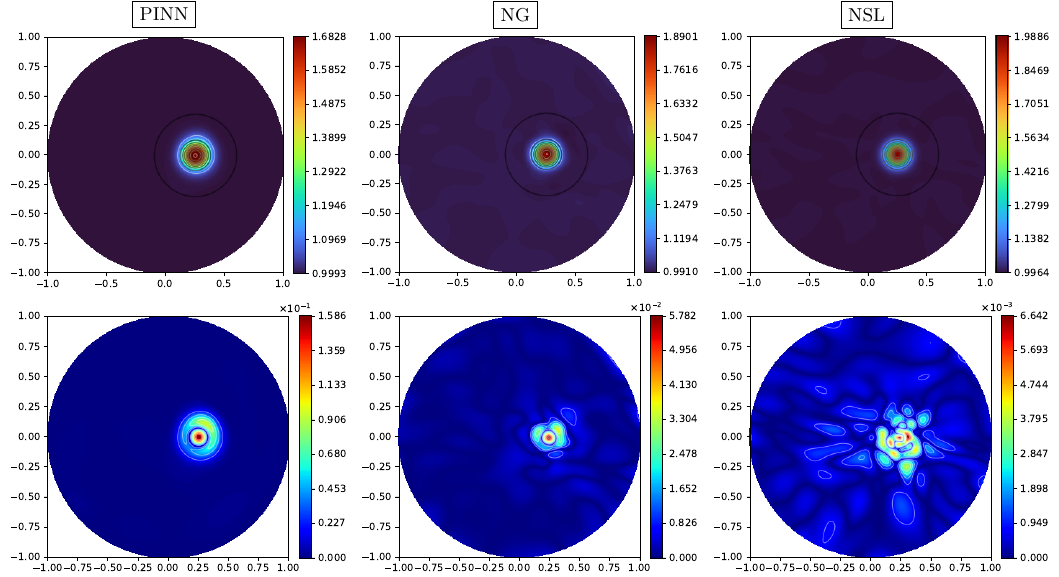}
    \caption{2D parametric rotating advection from \cref{sec:2D_transport_rotating_parametric}: from left to right, results of the PINN, NG and NSL schemes for $\mu = (0.06, 0.25)$, with a fixed time step for NG ($\dt = 0.02$) and NSL ($\dt = 0.5$). Top panels: approximate solutions; bottom panels: absolute value of the pointwise error.}
    \label{fig:2D_transport_rotating_parametric}
\end{figure}

The results are displayed on \cref{fig:2D_transport_rotating_parametric},
where we show the approximate solutions (top panels)
and the associated errors (bottom panels).
We observe a very good agreement of the NSL solution,
especially close to the peak of the Gaussian bump,
which is much better approximated with the NSL method
than with the other methods.

For a refined comparison of the methods,
we run the three methods on $\num{2000}$ randomly sampled parameters in $\mathbb{M}$,
and report the statistics on the errors
in \cref{tab:2D_transport_rotating_parametric}.
In both $L^2$ and $L^\infty$ norms,
we observe a significant gain in accuracy
for the NG and NSL methods over the PINN.
Furthermore, the NSL method outperforms the NG method,
by a factor of $3$ in $L^2$ norm,
and up to $5$ in $L^\infty$ norm.
\ra{For all cases, we note that the maximum errors
    are concentrated close to the boundaries of the parameter space.}

\begin{table}[!ht]
    \centering
    \begin{tabular}{lcccccc}
        \toprule
            & \multicolumn{2}{c}{NSL} & \multicolumn{2}{c}{NG} & \multicolumn{2}{c}{PINN}                                                        \\
        \cmidrule(lr){2-3} \cmidrule(lr){4-5} \cmidrule(lr){6-7}
            & $L^2$ error             & $L^\infty$ error       & $L^2$ error              & $L^\infty$ error & $L^2$ error    & $L^\infty$ error \\
        \cmidrule(lr){1-7}
        min & \num{3.79e-04}          & \num{5.73e-03}         & \num{1.05e-03}           & \num{1.89e-02}   & \num{1.20e-03} & \num{5.25e-02}   \\
        avg & \num{4.99e-04}          & \num{1.08e-02}         & \num{1.57e-03}           & \num{5.78e-02}   & \num{2.73e-03} & \num{1.36e-01}   \\
        max & \num{9.47e-04}          & \num{4.28e-02}         & \num{2.80e-03}           & \num{1.52e-01}   & \num{5.26e-03} & \num{2.70e-01}   \\
        std & \num{9.76e-05}          & \num{6.24e-03}         & \num{3.72e-04}           & \num{2.89e-02}   & \num{7.97e-04} & \num{5.20e-02}   \\
        \bottomrule
    \end{tabular}
    \caption{2D parametric rotating advection from \cref{sec:2D_transport_rotating_parametric}: statistics on the errors $e_x$ for $\num{2000}$ randomly sampled parameters in $\mathbb{M}$.}
    \label{tab:2D_transport_rotating_parametric}
\end{table}

\subsubsection{1D1V Vlasov equation}
\label{sec:2D_fake_Vlasov}

The next test case in this series is
a Vlasov equation
in one space dimension and one velocity dimension.
It is nothing but a 2D advection equation without parameters,
whose setup is described in e.g.~\cite{BerPeh2023}.
It is given by the following set of equations,
with periodic boundary conditions in $x$ and $v$:
\begin{equation}
    \label{eq:fake_Vlasov}
    \begin{dcases}
        \partial_t u + v \, \partial_x u + \sin(x) \, \partial_v u = 0,
         & \quad x \in (0, 2 \pi), \;
        v \in (-6, 6), \;
        t \in (0, 4.5),               \\
        u(0, x, v) = u_0(x, v),
         & \quad x \in \Omega, \;
        v \in (-6, 6),                \\
        u(t, 0, v) = u(t, 2 \pi, v),
         & \quad v \in (-6, 6), \;
        t \in (0, 4.5),               \\
        u(t, x, -6) = u(t, x, 6),
         & \quad x \in (0, 2 \pi), \;
        t \in (0, 4.5),               \\
    \end{dcases}
\end{equation}
with the initial condition
\begin{equation*}
    u_0(x, v) =
    \frac{1}{\sqrt{2 \pi}}
    \exp\left(-\frac{v^2}{2}\right).
\end{equation*}
This time, neither the advection equation
\eqref{eq:fake_Vlasov}
nor the characteristic ODE~\eqref{eq:characteristic_curve}
have a closed-form solution.
To obtain a reference solution,
we numerically solve \eqref{eq:fake_Vlasov} using
a classical semi-Lagrangian scheme\footnote{The semi-Lagrangian code is inspired from the one developed by Pierre Navaro, available at \url{https://pnavaro.github.io/python-notebooks/19-LandauDamping.html}.}
on a Cartesian grid,
with $512$ in space and $8192$ in velocity,
and with $200$ time steps per second
(leading to e.g.\ $900$ time steps for $t = 4.5$).
Moreover, to solve the characteristic ODE,
we use the third branch of \cref{alg:characteristic_curve},
with $n_\tau = 5$ sub-time steps.

\begin{figure}[!ht]
    \centering
    \includegraphics[width=\textwidth]{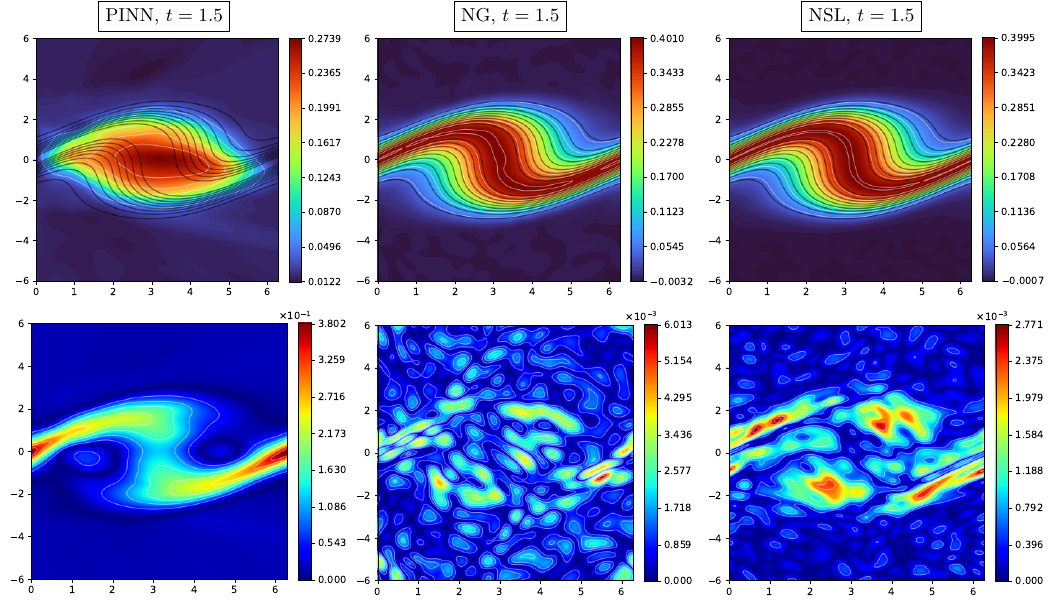}
    \caption{2D Vlasov equation from \cref{sec:2D_fake_Vlasov}: from left to right, results of the PINN, NG and NSL schemes at time $t=1.5$. Top panels: approximate solutions and contour lines of the exact solution (in black); bottom panels: absolute value of the pointwise error.}
    \label{fig:2D_fake_Vlasov_1.5}
\end{figure}

\begin{figure}[!ht]
    \centering
    \includegraphics[width=\textwidth]{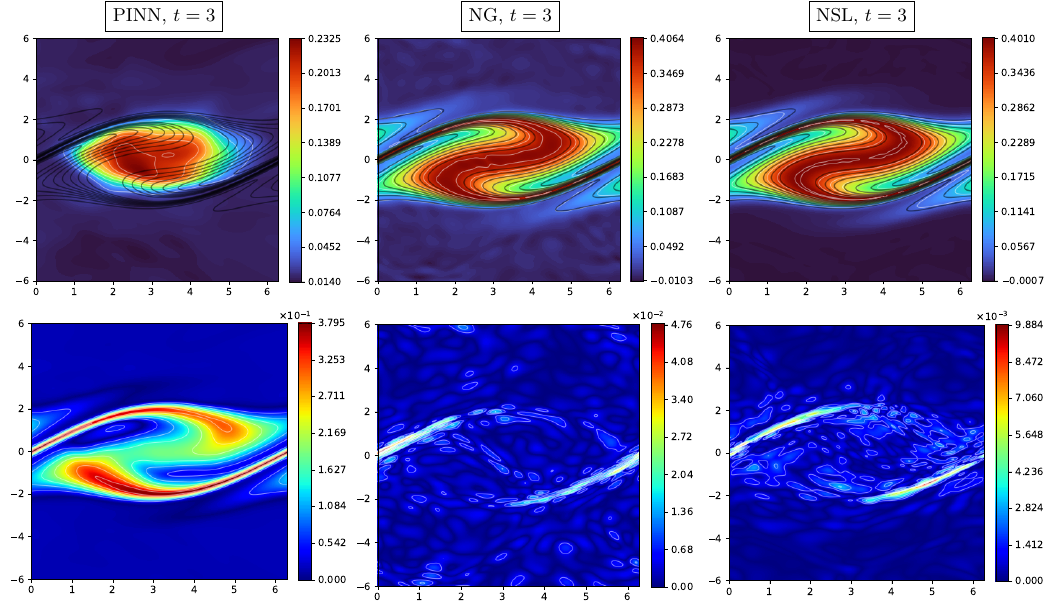}
    \caption{2D Vlasov equation from \cref{sec:2D_fake_Vlasov}: from left to right, results of the PINN, NG and NSL schemes at time $t=3$. Top panels: approximate solutions and contour lines of the exact solution (in black); bottom panels: absolute value of the pointwise error.}
    \label{fig:2D_fake_Vlasov_3.0}
\end{figure}

\begin{figure}[!ht]
    \centering
    \includegraphics[width=\textwidth]{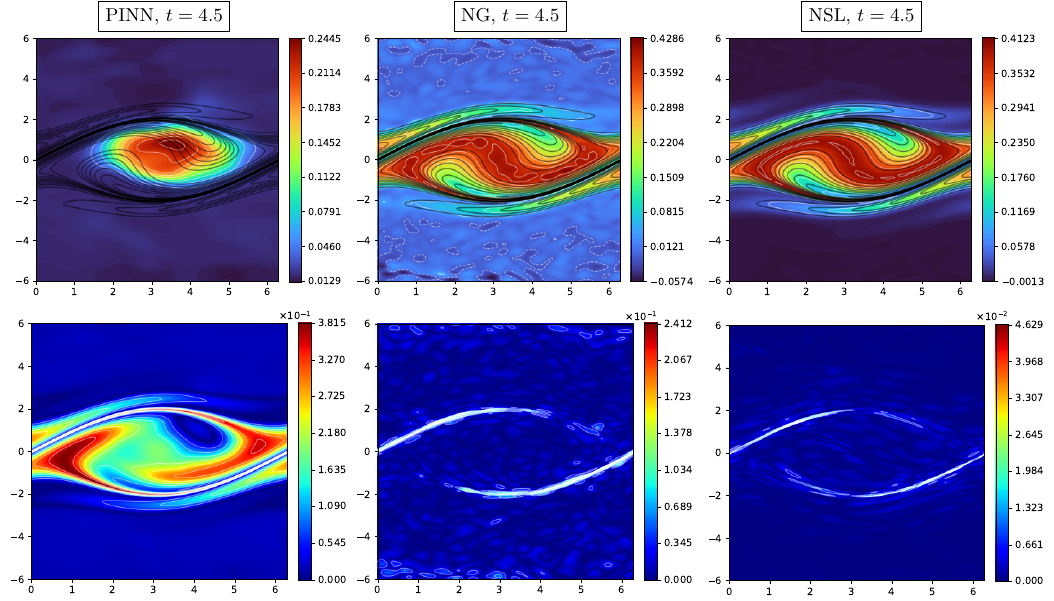}
    \caption{2D Vlasov equation from \cref{sec:2D_fake_Vlasov}: from left to right, results of the PINN, NG and NSL schemes at time $t=4.5$. Top panels: approximate solutions and contour lines of the exact solution (in black); bottom panels: absolute value of the pointwise error.}
    \label{fig:2D_fake_Vlasov_4.5}
\end{figure}

Since this problem is more complex than the previous ones,
in that it leads to a sheared solution with filaments,
we ran the simulation in double precision.
Moreover, the hyperparameters are adjusted accordingly.
Namely, more collocation points are used to approximate
the integrals.
All hyperparameters are given in
\cref{tab:hyperparameters_2D_fake_Vlasov}.
We elect to use a fixed time step
$\dt = 10^{-2}$ for the NG method,
corresponding to $450$ time iterations,
and $\dt = 1.5$ for the NSL method,
corresponding to $3$ time iterations.
This leads to a computation time of about
$5$ minutes for the PINN,
$1.5$ minute for the NSL method and
$4$ minutes for the NG method
(plus $30$ seconds for the initial condition).
The results are displayed in
\cref{fig:2D_fake_Vlasov_1.5,fig:2D_fake_Vlasov_3.0,fig:2D_fake_Vlasov_4.5}.

This experiment shows the limits of the PINN
when approximating a solution with fine structures.
Indeed, the results of the PINN
(displayed in the left panels of the figures)
are quite poor, and the solution is diffused away.
Since the PINN learns all times at once,
there is no reason why the early times
should be better approximated than the later times,
and this is confirmed by observing the results.

The NG scheme, on the other hand,
provides a better approximation
of the solution than the PINN.
However, even though the solution is less diffusive,
it remains more oscillatory,
especially close to the filaments.
In particular, we see that the solution becomes
negative in some regions, which is not physical.

Finally, the NSL method provides the best approximation
among the three tested methods.
Namely, the fine filaments present in
the solution are well-captured,
although the approximation quality decreases with time
since the filamentation process increases.
This decrease in quality
is not improved when decreasing the time step,
even locally in time
(e.g.\ taking two time steps between $t=3$ and $t=4.5$).
This motivates the next section,
where we will introduce preconditioning
and other techniques to further improve
the accuracy and efficiency of our method.

\begin{table}[!ht]
    \centering
    \begin{tabular}{ccccccc}
        \toprule
        \multirow{2}{*}{time $t$} & \multicolumn{2}{c}{PINN error} & \multicolumn{2}{c}{NG error} & \multicolumn{2}{c}{NSL error}                                                          \\
        \cmidrule(lr){2-3}  \cmidrule(lr){4-5}  \cmidrule(lr){6-7}
                                  & $L^2$ error                    & $L^\infty$ error             & $L^2$ error                   & $L^\infty$ error & $L^2$ error      & $L^\infty$ error \\
        \cmidrule(lr){1-7}
        $1.5$                     & $\num{7.98e-02}$               & $\num{3.80e-01}$             & $\num{1.19e-03}$              & $\num{6.01e-03}$ & $\num{6.76e-04}$ & $\num{2.77e-03}$ \\
        $3$                       & $\num{1.06e-01}$               & $\num{3.80e-01}$             & $\num{3.39e-03}$              & $\num{4.76e-02}$ & $\num{8.22e-04}$ & $\num{9.88e-03}$ \\
        $4.5$                     & $\num{1.17e-01}$               & $\num{3.82e-01}$             & $\num{1.48e-02}$              & $\num{2.41e-01}$ & $\num{1.50e-03}$ & $\num{4.63e-02}$ \\
        \bottomrule
    \end{tabular}
    \caption{Comparison of $L^2$ and $L^\infty$ errors for the PINN and the NG and NSL methods at different times.}
    \label{tab:errors_wrt_time}
\end{table}

Lastly, we report in \cref{tab:errors_wrt_time}
the $L^2$ and $L^\infty$ errors
of each method,
bearing in mind that the $L^\infty$ error
may not be very informative in this case,
since we are comparing solutions with sharp gradients.
These errors confirm our observations from
\cref{fig:2D_fake_Vlasov_1.5,fig:2D_fake_Vlasov_3.0,fig:2D_fake_Vlasov_4.5},
namely that PINN is the least accurate method,
while NSL is the most accurate one,
improving the $L^2$ error by a factor of about $10$
and the $L^\infty$ error by a factor of about~$5$
compared to the NG method,
all for a third of the computational cost.

\subsection{Transport equation in a cylinder}
\label{sec:3D_cylinder}

From this section onwards,
we will only consider the NSL method,
and check it on several benchmark problems.
The present section is devoted to testing
the improvements proposed in \cref{sec:improving_accuracy}.
To that end, we consider an advection equation
in three space dimensions and two parameter dimensions.
The space domain is a cylinder
$\Omega = \mathbb{D}^1 \times [0, 2]$.
The solution $u$ is governed by the advection equation
\begin{equation*}
    \partial_t u + a(x) \cdot \nabla u = 0,
    \quad x = (x_1, x_2, x_3)^\intercal \in \Omega, \;
    t \in (0, 2), \;
    \mu \in \mathbb{M},
\end{equation*}
supplemented with periodic boundary conditions
in the third variable,
and a suitable initial condition.
Namely, we set the advection field to
\begin{equation}
    \label{eq:3D_cylinder_advection_vector_field}
    a(x) = \begin{pmatrix}
        - 2 \pi x_2 \\ 2 \pi x_1 \\ x_3
    \end{pmatrix},
\end{equation}
corresponding to a rotation in the first two variables
and a constant advection in the third one.
Moreover, the (parameter-dependent) exact solution
is a bump function, which reads
\begin{equation*}
    u_\text{ex}(t, x, \mu)
    =
    1 + \exp \left(
    -\frac {1} {2 v^2} \Big(
    \big(x_1 - c \, x_1^0(t)\big)^2 +
    \big(x_2 - c \, x_2^0(t)\big)^2 +
    x_3^0(t)^2
    \Big)
    \right),
\end{equation*}
where we have set
\begin{equation*}
    x_1^0(t) = \cos(2 \pi t),
    \qquad
    x_2^0(t) = \sin(2 \pi t),
    \qquad
    x_3^0(t) = (x_3 - t) \% 2 - 1,
\end{equation*}
where $a \% b$ denotes the modulo operation.
The two parameters
$c \in (0.3, 0.5)$ and $v \in (0.05, 0.15)$
respectively represent
the position of the center of the Gaussian bump,
and its variance.
The exact solution corresponds to a rotation of the bump
around the $x_3$-axis,
with a constant advection in the $x_3$ direction.
We take a time step $\dt = 0.5$,
corresponding to $4$ time steps.

We compare the results of the NSL method,
as described before,
to the improved version equipped with
natural gradient preconditioning (\cref{sec:natural_gradient_preconditioning}) and
adaptive sampling (\cref{sec:adaptive_sampling}).
The adaptive sampling is based on the zeroes of the
gradient of the approximate solution, i.e.,
\begin{equation*}
    f = \frac 1 {10^{-2} + \| \nabla u_\theta \|^2}
\end{equation*}
in \cref{alg:adaptive_sampling}, where $10^{-2}$
is a safety factor to avoid division by zero.
Furthermore, we take $\sigma_1 = 0.5$,
$\sigma_2 = 10$ and $\sigma_3 = 500$
in \cref{alg:adaptive_sampling}.
We test three configurations for the improved NSL method:
first with a small network (configuration (a)),
then with a larger network
but few epochs (configuration (b)),
and lastly with a larger network,
as well as additional epochs
and collocation points (configuration (c)).
The hyperparameters are summarized in
\cref{tab:hyperparameters_3D_rotation},
for each configuration.

\begin{table}[!ht]
    \centering
    \begin{tabular}{ccccc}
        \toprule
        \multirow{2}{*}{Hyperparameter}
                      & \multirow{2}{*}{NSL}
                      & \multicolumn{3}{c}
        {improved NSL}
        \\
        \cmidrule(lr){3-5}
                      &
                      & configuration (a)
                      & configuration (b)
                      & configuration (c)
        \\
        \cmidrule(lr){1-5}
        $N_e$ (init.) & $\num{1000}$
                      & $\num{150}$          & $\num{150}$    & $\num{500}$    \\
        $N_e$ (iter.) & $\num{1000}$
                      & $\num{50}$           & $\num{50}$     & $\num{200}$    \\
        $N_c$         & $\num{50000}$
                      & $\num{50000}$        & $\num{50000}$  & $\num{150000}$ \\
        $\ell$        & $[60, 60, 60]$
                      & $[20, 40, 20]$       & $[60, 60, 60]$ & $[60, 60, 60]$ \\
        \bottomrule
    \end{tabular}
    \caption{Hyperparameters (number of epochs $N_e$ and layers $\ell$) used for initialization (init.) and NSL iterations (iter.) in \cref{sec:3D_cylinder}. In every case, we use a $\tanh$ activation function.}
    \label{tab:hyperparameters_3D_rotation}
\end{table}

To report the results,
we perform a parametric study
over $\num{2000}$ randomly sampled parameters
in $\mathbb{M}$,
and report the statistics in \cref{tab:3D_cylinder},
along with the computation time.
The first configuration of the improved NSL method
decreases the computation time by a factor of roughly $10$,
and decreases the average error
by a factor of $10$ in $L^2$ norm,
$20$ in $L^\infty$ norm.
These substantial increases in accuracy and efficiency
make the improved NSL method even more competitive
with other neural methods.
We remark that the network is really quite small,
and natural gradient preconditioning shows that
solving the optimization problem is an important bottleneck
in the training of the neural network.
Moving on to larger networks and additional epochs,
we observe a further decrease of the error,
alongside a (substantial) increase in the computation time.
We observe that, in configuration (b),
the increase in epochs in not sufficient
to improve the accuracy and solve the optimization problems
more efficiently than in configuration (a).
However, in configuration (c),
the further increase in both epochs and collocation points
allows us to further decrease the error
by a factor of around $4$.
The downside is that the computation time is increased
by a factor of about $45$.
All in all, we see that the improved NSL method
is able to provide a very accurate solution
for a small computational cost,
and that this accuracy can be improved further
should one be willing to pay the computation time.

\begin{table}[!ht]
    \centering
    \begin{tabular}{llcccc}
        \toprule
         &
         & min                                   & avg            & max            & std            \\
        \cmidrule(lr){1-6}
        \multirow{4}{*}{non-improved}
         & $L^2$ error
         & \num{2.72e-03}                        & \num{5.22e-03} & \num{1.08e-02} & \num{1.69e-03} \\
         & $L^\infty$ error
         & \num{5.55e-02}                        & \num{1.91e-01} & \num{3.75e-01} & \num{9.97e-02} \\
         & computation time (init.)
         & \multicolumn{4}{c}{\qty{121.78}{\s}}                                                     \\
         & computation time (iter.)
         & \multicolumn{4}{c}{\qty{296.58}{\s}}                                                     \\
        \cmidrule(lr){1-6}
        \multirow{4}{*}{configuration (a)}
         & $L^2$ error
         & \num{3.44e-04}                        & \num{4.36e-04} & \num{9.39e-04} & \num{1.23e-04} \\
         & $L^\infty$ error
         & \num{3.01e-03}                        & \num{1.03e-02} & \num{7.06e-02} & \num{1.23e-02} \\
         & computation time (init.)
         & \multicolumn{4}{c}{\qty{19.06}{\s}}                                                      \\
         & computation time (iter.)
         & \multicolumn{4}{c}{\qty{26.35}{\s}}                                                      \\
        \cmidrule(lr){1-6}
        \multirow{4}{*}{configuration (b)}
         & $L^2$ error
         & \num{5.85e-04}                        & \num{7.21e-04} & \num{1.26e-03} & \num{1.30e-04} \\
         & $L^\infty$ error
         & \num{3.16e-03}                        & \num{8.49e-03} & \num{6.03e-02} & \num{9.58e-03} \\
         & computation time (init.)
         & \multicolumn{4}{c}{\qty{110.20}{\s}}                                                     \\
         & computation time (iter.)
         & \multicolumn{4}{c}{\qty{149.87}{\s}}                                                     \\
        \cmidrule(lr){1-6}
        \multirow{4}{*}{configuration (c)}
         & $L^2$ error
         & \num{1.04e-04}                        & \num{1.36e-04} & \num{3.10e-04} & \num{3.33e-05} \\
         & $L^\infty$ error
         & \num{9.84e-04}                        & \num{2.74e-03} & \num{3.00e-02} & \num{4.08e-03} \\
         & computation time (init.)
         & \multicolumn{4}{c}{\qty{774.73}{\s}}                                                     \\
         & computation time (iter.)
         & \multicolumn{4}{c}{\qty{1267.12}{\s}}                                                    \\
        \bottomrule
    \end{tabular}
    \caption{3D parametric rotating advection from \cref{sec:3D_cylinder}: statistics on the errors $e_x$ for $\num{2000}$ randomly sampled parameters in $\mathbb{M}$.}
    \label{tab:3D_cylinder}
\end{table}

\begin{figure}[!ht]
    \centering
    \includegraphics{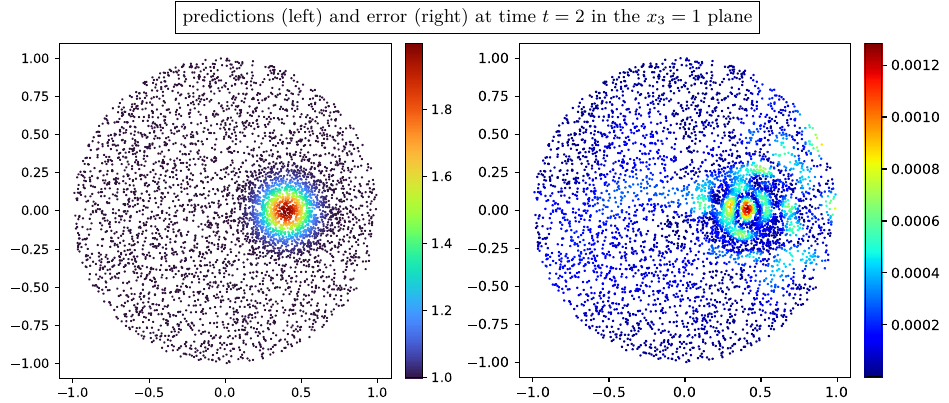}
    \caption{3D parametric advection from \cref{sec:3D_cylinder}, in configuration (c): approximate solution (left) and error (right) in the $x_3=1$ plane at the final time $t=2$, for the physical parameters $c = 0.4$ and $v = 0.1$.}
    \label{fig:3D_cylinder_z}
\end{figure}

\begin{figure}[!ht]
    \centering
    \includegraphics{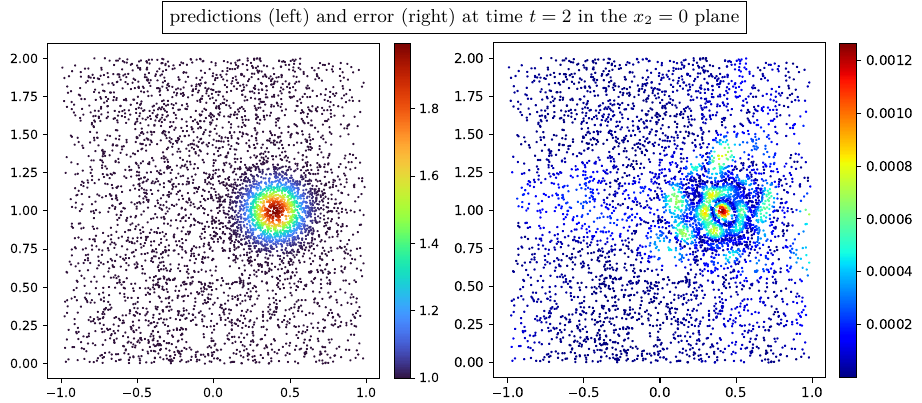}
    \caption{3D parametric advection from \cref{sec:3D_cylinder}, in configuration (c): approximate solution (left) and error (right) in the $x_2=0$ plane at the final time $t=2$, for the physical parameters $c = 0.4$ and $v = 0.1$.}
    \label{fig:3D_cylinder_y}
\end{figure}

To give an idea of the adaptively sampled points,
we represent in \cref{fig:3D_cylinder_z,fig:3D_cylinder_y}
the predictions (left panels)
and associated errors (right panels)
of the improved method, in configuration (c),
at the final time $t=2$,
for the physical parameters $c = 0.4$ and $v = 0.1$,
and in two planes where the Gaussian bump is localized
($x_3=1$ for \cref{fig:3D_cylinder_z}
and $x_2=0$ for \cref{fig:3D_cylinder_y}).
We have displayed $\num{5000}$ adaptively sampled points
out of the $\num{150000}$ collocation points
actually used in this test case.
In both cases,
in addition to noticing that the errors produced
by the improved NSL method are quite low,
we observe that the adaptively sampled points
are concentrated around the peak of the Gaussian bump,
where the solution is the most challenging to approximate,
and where the error is the largest,
even though the adaptive sampling procedure
is not aware of the exact solution.

\begin{remark}
    \label{rmk:higher_cpu_time}
    For the remainder of the numerical experiments,
    we elect to use hyperparameters
    similar to configuration (c).
    This leads to higher computation times,
    but better accuracy for the optimization problems.
    Using hyperparameters similar to configuration (a)
    would lead to significantly
    lower computation times,
    at the cost of a decrease in accuracy.
\end{remark}

\subsection{Deformation of level-set functions}
\label{sec:level_sets}

Equipped with the improved version of NSL,
we now test it on two challenging test cases,
namely the deformation of level-set functions,
\rall{in 2D (\cref{sec:2D_level_set}),
    in 3D (\cref{sec:3D_level_set}),
    and in 3D with 2 parameters (\cref{sec:3D_level_set_param})}.
These test cases,
described in e.g.~\cite{LeV1996, BuiDapFre2011},
consist in the transport of a level-set function
by a time-dependent advection field.
They are particularly interesting since,
while the transient solution is not known,
the final solution is the initial condition itself,
which allows for
a better assessment of the accuracy of the method.
Since we are approximating a level-set function,
we use $f = u_\theta$ in \cref{alg:adaptive_sampling};
we take $\sigma_1 = \num{2e-3}$, $\sigma_2 = \num{1e-2}$ and $\sigma_3 = \num{5e-2}$.
In this section, we take $n_\tau = 10$ sub-time steps
in the characteristic curve solver.

\subsubsection{Deformation of a 2D level-set function}
\label{sec:2D_level_set}

The initial condition of this two-dimensional level-set
advection test case is,
accordingly, a level-set function.
It is defined in the space domain
$\Omega = [0, 1]^2$, and it represents the disk
of radius $0.15$ centered at $(0.5, 0.75)$, i.e.,
\begin{equation*}
    u_0(x) =
    \left(x_1 - \frac 1 2 \right)^2
    + \left(x_2 - \frac 3 4 \right)^2
    - 0.15^2.
\end{equation*}
The shape then undergoes a deformation,
according to the following time-dependent
advection field:
\begin{equation*}
    a(t, x) = \begin{pmatrix}
        -\sin^2(\pi x_1) \sin(2 \pi x_2) \cos(\frac{\pi t}{T}) \\
        \sin^2(\pi x_2) \sin(2 \pi x_1) \cos(\frac{\pi t}{T})
    \end{pmatrix},
\end{equation*}
where $T$ is the final time.
The shape is most deformed at time $t = T/2$,
and goes back to the initial shape at time $t = T$.
The hyperparameters for this experiment are given in
\cref{tab:hyperparameters_level_set}.
We take a time step $\dt = 0.2$,
corresponding to $40$ time steps.
Note that additional time steps are performed
compared to previous experiments;
this is to obtain a good discretization of
the time-dependent advection field.
With this configuration,
the full computation takes about one hour.

\begin{figure}[!ht]
    \centering
    \includegraphics{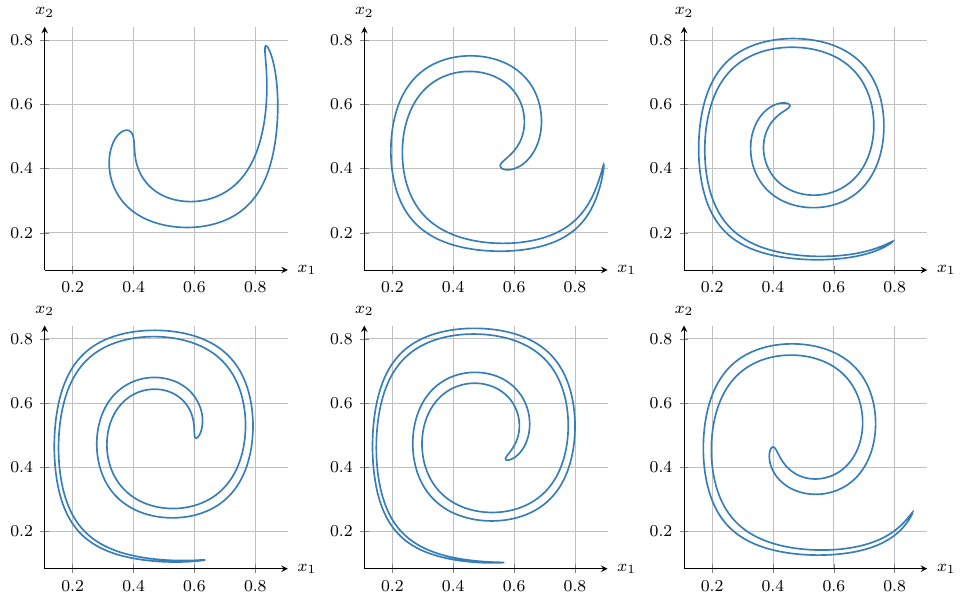}
    \caption{2D level-set deformation from \cref{sec:2D_level_set}: approximate solution at several times (from left to right and top to bottom, $t=0.8$, $t=1.6$, $t=2.4$, $t=3.2$, $t=4$ and $t=6$). Only the zero contour of the level-set function is displayed.}
    \label{fig:2D_level_set}
\end{figure}

\cref{fig:2D_level_set} depicts the zero contour
of the approximate level-set function
at several times.
We observe a good agreement with the reference solution
from \cite{LeV1996, BuiDapFre2011}.
Namely, the deformation of the shape is well-captured,
even at its most deformed, for $t = T/2 = 4$.

\begin{figure}[!ht]
    \centering
    \includegraphics{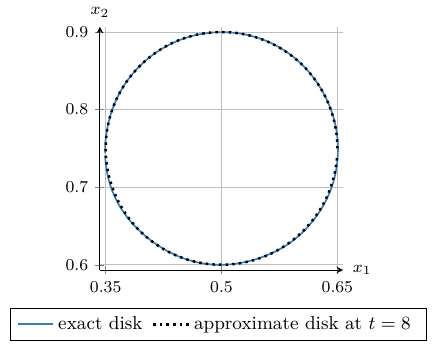}
    \caption{2D level-set deformation from \cref{sec:2D_level_set}: comparison of the exact solution (dotted black line) and the NSL solution (solid blue line) at the final time $t = T$. Only the zero contour of the level-set function is displayed.}
    \label{fig:2D_level_set_comparison}
\end{figure}

To get a better view of the accuracy of the NSL method,
we compare the NSL solution at $t=8$ to the exact solution,
which is nothing but the initial condition.
Namely, the zero contour of the exact solution
is the disk of radius $0.15$ centered at $(0.5, 0.75)$.
This comparison is carried out in
\cref{fig:2D_level_set_comparison},
where we observe a very good agreement
between the two solutions,
further validating the accuracy of the NSL method.
\ra{For a more quantitative assessment,
    we compute the error between the true volume
    ($0.15^2 \pi \approx \num{0.070685}$)
    and the approximate one,
    computed by the Monte-Carlo method.
    We obtain
    \begin{equation*}
        \int_{\Omega} \boldone_{\{u_\theta < 0\}} \, dx
        \approx \num{0.070924},
    \end{equation*}
    which leads to a relative error of about \qty{0.338}{\percent}.}

\subsubsection{Deformation of a 3D level-set function}
\label{sec:3D_level_set}

After the 2D level-set deformation,
we move on to a more complex three-dimensional test case.
This time, the initial condition is a level-set function,
defined in $\Omega = [0, 1]^3$,
of the sphere of radius $0.15$
centered at $(0.35, 0.35, 0.35)$.
This initial condition is given by
\begin{equation*}
    u_0(x) =
    \left(x_1 - 0.35 \right)^2
    + \left(x_2 - 0.35 \right)^2
    + \left(x_3 - 0.35 \right)^2
    - 0.15^2.
\end{equation*}
It is then advected
with the time-dependent advection field
\begin{equation}
    \label{eq:3D_level_set_advection_vector_field}
    a(t, x) =
    \begin{pmatrix}
        \sin^2(\pi x_1) \sin(2 \pi x_2) \sin(2 \pi x_3) \cos(\frac{\pi t}{3}) \\
        \sin^2(\pi x_2) \sin(2 \pi x_1) \sin(2 \pi x_3) \cos(\frac{\pi t}{3}) \\
        \sin^2(\pi x_3) \sin(2 \pi x_1) \sin(2 \pi x_2) \cos(\frac{\pi t}{3})
    \end{pmatrix}.
\end{equation}
Similarly to \cref{sec:2D_level_set},
we recover the initial condition for $t = 3$.
We take a time step $\dt = 0.3$,
which corresponds to $10$ time steps.
This leads,
together with the hyperparameters
reported in \cref{tab:hyperparameters_level_set},
lead to a computation time of about $90$ minutes.
Same as the 2D level-set deformation,
additional points are sampled around
the zeroes of the approximate solution.

\begin{figure}[!ht]
    \centering
    \includegraphics[width=0.75\textwidth]{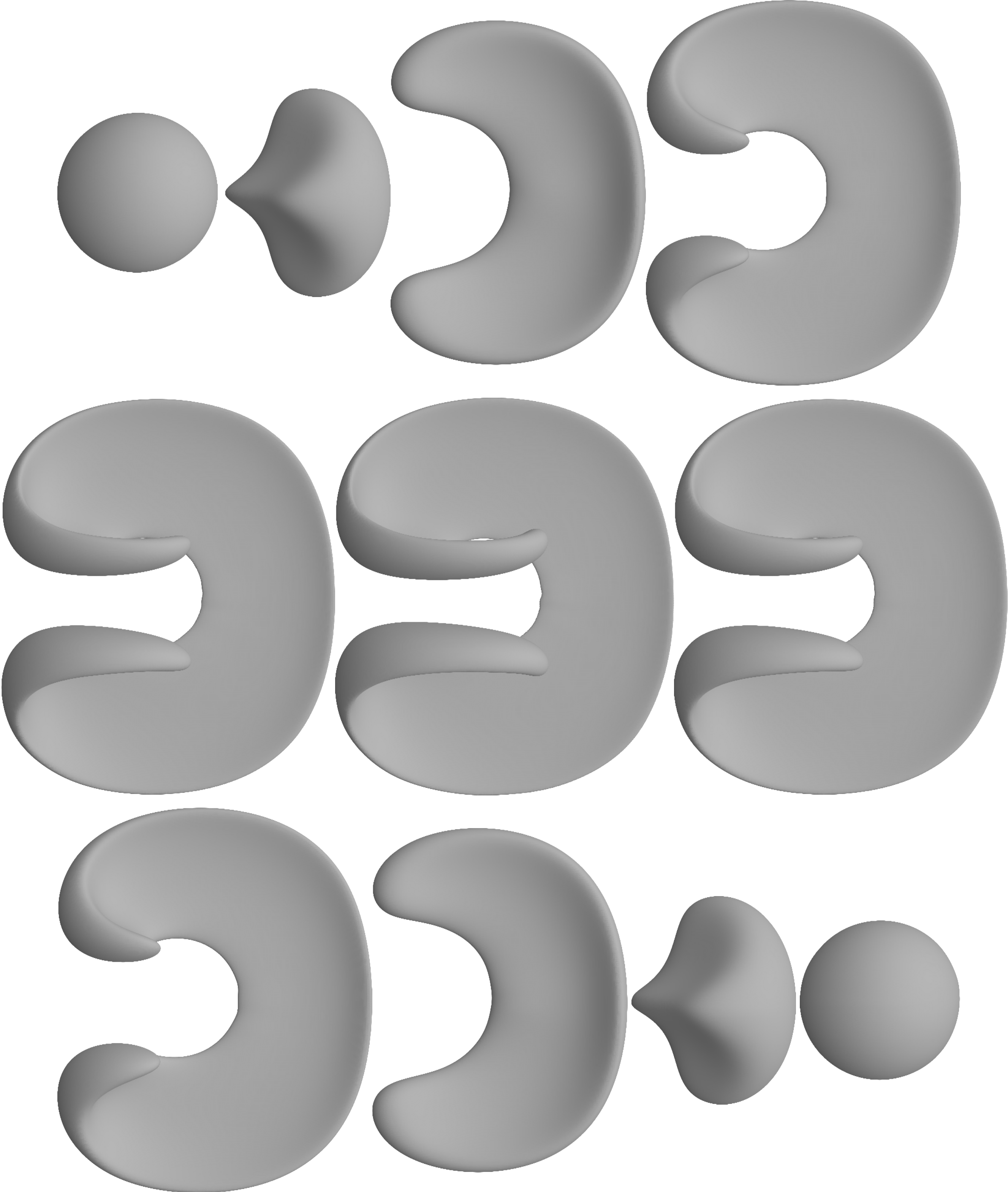}
    \caption{3D level-set deformation from \cref{sec:3D_level_set}: zero contour of the approximate solution at all $10$ computed times (from left to right and top to bottom, starting from $t=0$ on the top left and finishing with $t=3$ on the bottom right).}
    \label{fig:3D_level_set}
\end{figure}

First, in \cref{fig:3D_level_set}, we display
the zero contour of the approximate level-set function
at all $10$ computed times between $0$ and $3$,
with a step of $\dt$.
The approximate zero contour is well-captured
by the NSL scheme,
as can be seen by comparing it to the reference
solution from \cite{LeV1996,BuiDapFre2011}.

\begin{figure}[!ht]
    \centering
    \includegraphics{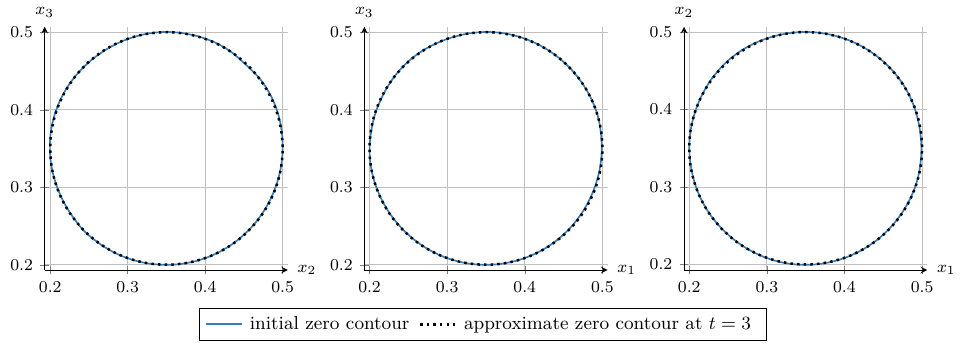}
    \caption{3D level-set deformation from \cref{sec:3D_level_set}: comparison of the exact solution (dotted black line) and the NSL solution (solid blue line) at the final time $t = 3$. Only the zero contour of the level-set function is displayed. From left to right, we display the sphere sliced by planes $x_1 = 0.35$, $x_2 = 0.35$ and $x_3=0.35$.}
    \label{fig:3D_level_set_comparison}
\end{figure}

Second,
\cref{fig:3D_level_set_comparison} shows the zero contour
of the approximate level-set function,
comparing it to the exact solution $u_0$.
To do so, we display the 3D sphere
sliced by the planes $x = 0.35$, $y = 0.35$ and $z = 0.35$,
leading to three 2D graphs.
Once again, we observe a very good agreement
between the NSL solution and the exact solution,
even for this more complex three-dimensional problem.
\ra{Quantitatively, the exact volume of the sphere is
$\frac 4 3 \pi (0.15)^3 \approx \num{0.014137}$,
while the approximate volume computed by the NSL method
is roughly equal to \num{0.014108}.
The relative error is thus of about \qty{0.209}{\percent}.}

{\rallparagraph
\subsubsection{Deformation of a 3D level-set function with two parameters}
\label{sec:3D_level_set_param}

Lastly, we add two parameters to the previous 3D test case.
The initial condition becomes a level-set function
of a perturbed ellipsoid centered at $(0.35, 0.35, 0.35)$.
The space domain remains $\Omega = [0, 1]^3$,
and the parameters $\mu = (\alpha_1, \epsilon)$
live in the parameter domain
$\mathbb{M} = [0.8, 1.2] \times [0, 0.25]$.
The first parameter, $\alpha_1$,
is a measure of the $x$-major axis of the ellipsoid,
while the second parameter, $\epsilon$,
dictates how much the ellipsoid is perturbed.
To give the initial condition, let us first define
\begin{equation*}
    x_{1,c} = x_1 - 0.35,
    \quad
    x_{2,c} = x_2 - 0.35,
    \quad
    x_{3,c} = x_3 - 0.35.
\end{equation*}
The radius and polar angle are defined by
\begin{equation*}
    r^2 =
    \bigg(\frac{x_{1,c}}{\alpha_1}\bigg)^2
    + \bigg(\frac{x_{2,c}}{\alpha_1^{-1}}\bigg)^2
    + \bigg(\frac{x_{3,c}}{1.2}\bigg)^2
    \text{\quad and \quad}
    \varphi = \polaratan\left(
    \sqrt{x_{1,c}^2 + x_{2,c}^2}, \, x_{3,c}
    \right),
\end{equation*}
and the initial condition is then
\begin{equation*}
    u_0(x, \mu) = 2 \left(
    r^2 \left(
        1 + \epsilon \sin(9 \varphi) e^{-1.5 \left( \varphi - \frac{\pi}{2} \right)^2}
        \right)
    - 0.15^2
    \right).
\end{equation*}
This initial condition is advected
with the same velocity field
\eqref{eq:3D_level_set_advection_vector_field}
as in \cref{sec:3D_level_set},
and the initial condition is recovered at $t = 3$.
We still take a time step $\dt = 0.3$ and $10$ time steps.
The hyperparameters are given in
\cref{tab:hyperparameters_level_set}.
The total computation time is around $3$ hours,
twice as long as the non-parametric case.
Same as the other two level-set deformations,
the adaptive sampling focuses on
the zeroes of the approximate solution.

\begin{figure}[!ht]
    \centering
    \includegraphics{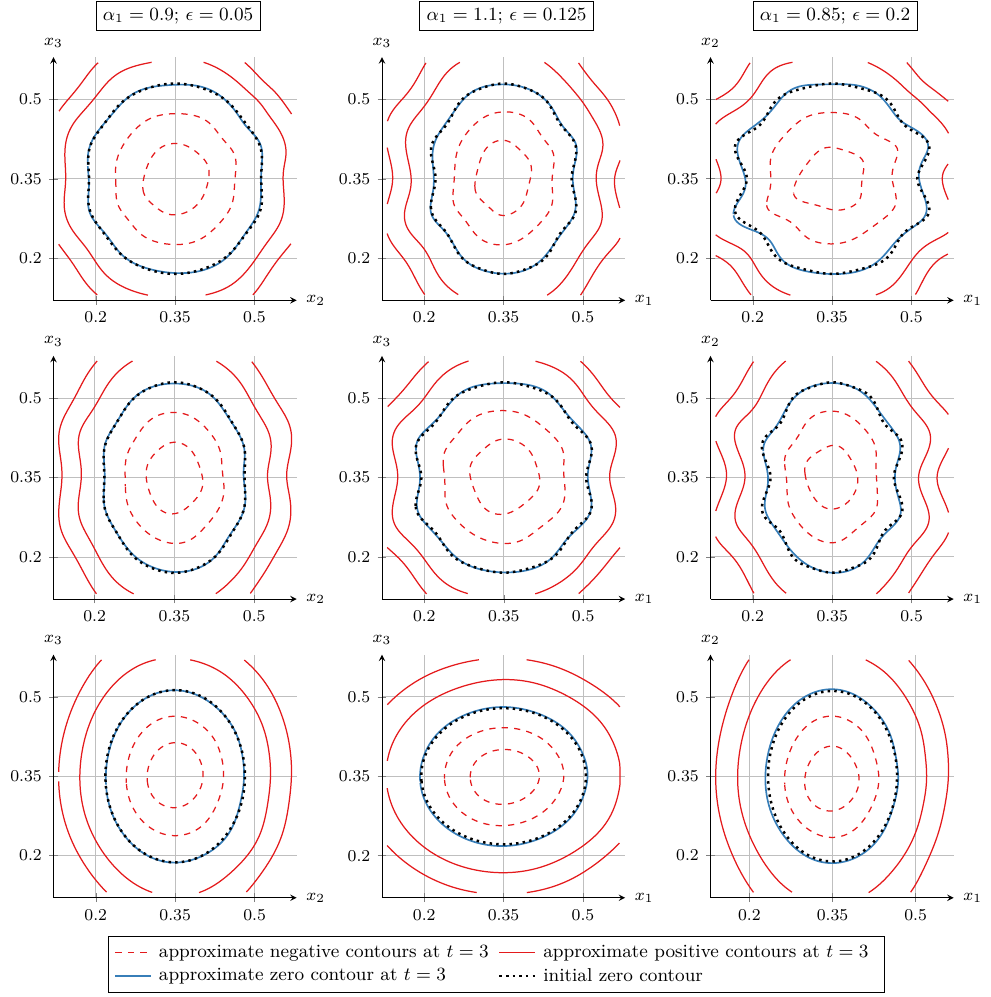}
    \caption{\rallparagraph Deformation of a 3D level-set function with two parameters from \cref{sec:3D_level_set_param}: comparison of the exact solution (dotted black line) and the NSL solution (solid blue line; solid and dashed red lines) at the final time $t = 3$. The solid blue line corresponds to the zero contour of the final solution; the solid and dashed red lines respectively correspond to positive and negative contours of the final solution. From top to bottom, we display the sphere sliced by planes $x_1 = 0.35$, $x_2 = 0.35$ and $x_3=0.35$. From left to right, three values of the parameters $\mu$ are considered.}
    \label{fig:3D_parametric_level_set_comparison}
\end{figure}

The results are depicted in
\cref{fig:3D_parametric_level_set_comparison}.
First, like in the previous case,
we observe a very good agreement
between zero contours of
the approximate and exact solutions,
as evidenced by the closeness
of the solid blue line and the dotted black line.
Of course, this problem is made harder by the
parametrized perturbation of the solution.
Hence, the approximate and exact solutions
are not as nicely superimposed as in
\cref{fig:3D_level_set_comparison},
especially for highly perturbed solutions
(represented in the right panels).
Second, the positive and negative contours
(solid and dashed red lines)
emphasize that the approximate solution
remains a level-set function,
even after having undergone the deformation.
Third, the volume preservation remains good:
across \num{128} instances of parameters in $\mathbb{M}$,
the average relative error on the volume
is $\qty{0.468}{\percent}$,
with a standard deviation of $\num{3.39e-3}$
and a maximal error of $\qty{1.24}{\percent}$.
The method is thus about half as accurate
(in terms of volume preservation)
as for the non-parametric problem,
which remains well within an acceptable range
for such a jump in complexity.
}

\subsection{High-dimensional advection-diffusion equations}
\label{sec:advection_diffusion_experiments}

This last section is dedicated to the approximation
of advection-diffusion equations,
following \cref{sec:advection_diffusion}.
In a $d$-dimensional space domain
$\Omega \subset \mathbb{R}^d$,
we define a constant advection field
$a = (1, 1, \ldots, 1)^\intercal$
and a constant diffusion coefficient $\sigma$.
The advection-diffusion equation
\eqref{eq:advection_diffusion}
then reads
\begin{equation*}
    \label{eq:advection_diffusion_experiments}
    \partial_t u + a \cdot \nabla u - \sigma \Delta u = 0,
    \quad x \in \Omega, \;
    t \geqslant 0.
\end{equation*}
We first check the convergence of the scheme
in \cref{sec:advection_diffusion_convergence},
and then we test it on two high-dimensional problems:
a periodic solution in \cref{sec:advection_diffusion_periodic}
and a Gaussian solution in \cref{sec:advection_diffusion_gaussian}.
\rabis{In particular, we perform comparisons
with a classical semi-Lagrangian scheme to
evaluate the performance of the method,
both in terms of accuracy and computational cost.
Note, however, that the comparison will naturally favor the proposed method
in a high-dimensional setting,
since classical methods such as the semi-Lagrangian scheme
suffer from the curse of dimensionality.}
\rall{A discussion of the results is finally
    provided in \cref{sec:advection_diffusion_discussion}.}

\subsubsection{Convergence study}
\label{sec:advection_diffusion_convergence}

For this study,
the space domain $\Omega = (-1, 1)^d$ is
supplemented with periodic boundary conditions.
\rb{We define the exact solution, for all $x \in \Omega$
    and $t \geqslant 0$, by
    \begin{equation*}
        u_\text{ex}(t, x) =
        2 + \sin\big( \pi \langle x - a t - s, a \rangle\big) \,
        e^{-\sigma \pi^2 \langle a, a \rangle t},
    \end{equation*}
    with \smash{$s = (s_i)_{i=1}^d \in \mathbb{R}^d$},
    where for all $i = 1, \ldots, d$,
    $s_i = (i - 1)/d$,
    is used to make the problem harder by
    introducing a different phase shift in each dimension.
    The rank of the solution, given by trigonometric formulas, is~$2^{d-1}$.}
In practice, the diffusion coefficient is set to $\sigma = 0.1$.
For this experiment,
natural gradient preconditioning is used,
but not adaptive sampling,
since there are no obvious areas in which
additional points should be sampled.
Moreover, \rall{$\sin$} activation functions are employed.

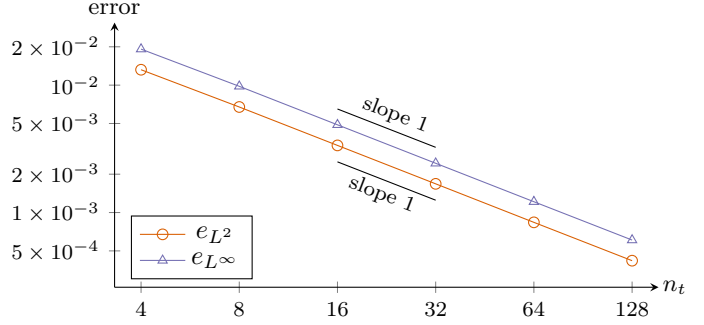
\begin{figure}[!ht]
    \centering
    \begin{subfigure}{0.35\textwidth}
        \centering
        \rallparagraph \begin{tabular}{ccc}
            \toprule
            $n_t$ & $e_{L^2}$      & $e_{L^\infty}$ \\
            \cmidrule(lr){1-3}
            4     & \num{1.32e-02} & \num{1.91e-02} \\
            8     & \num{6.73e-03} & \num{9.79e-03} \\
            16    & \num{3.36e-03} & \num{4.88e-03} \\
            32    & \num{1.68e-03} & \num{2.43e-03} \\
            64    & \num{8.37e-04} & \num{1.21e-03} \\
            128   & \num{4.18e-04} & \num{6.08e-04} \\
            \bottomrule
        \end{tabular}
        \caption{Errors for advection-diffusion (in dimension $d=3$) across different numbers of time steps.}
        \label{tab:advection_diffusion_time_steps}
    \end{subfigure}%
    \hfill
    \begin{subfigure}{0.55\textwidth}
        \centering
        \begin{tikzpicture}
            \begin{loglogaxis}[
                    axis lines = left,
                    enlarge x limits={abs=10pt},
                    enlarge y limits={abs=10pt},
                    xlabel style={at={(ticklabel* cs:1.01)},anchor=west},
                    ylabel style={at={(ticklabel* cs:1.01)},anchor=west},
                    label style={font=\small},
                    tick label style={font=\footnotesize},
                    xlabel = {\makebox[1pt][c]{$n_t$}},
                    ylabel = {\makebox[1pt][c]{\rotatebox{270}{error}}},
                    hide obscured x ticks=false,
                    hide obscured y ticks=false,
                    scaled y ticks = false,
                    scaled x ticks = false,
                    xtick={2, 4, 8, 16, 32, 64, 128},
                    xticklabels={2, 4, 8, 16, 32, 64, 128},
                    ytick={5e-4, 1e-3, 2e-3, 5e-3, 1e-2, 2e-2},
                    yticklabels={$\num{5e-4}$, $\num{1e-3}$, $\num{2e-3}$, $\num{5e-3}$, $10^{-2}$, $\num{2e-2}$},
                    grid=minor,
                    legend pos=south west,
                    width=0.95\textwidth,
                    height=0.55\textwidth,
                ]
                \addplot[mark=o, color=graph_2] coordinates {
                        (4,1.321e-02) (8,6.737e-03)
                        (16,3.365e-03) (32,1.680e-03) (64,8.375e-04) (128,4.189e-04)
                    };
                \addlegendentry{$e_{L^2}$}
                \addplot[mark=triangle, color=graph_1] coordinates {
                        (4,1.918e-02) (8,9.790e-03)
                        (16,4.883e-03) (32,2.436e-03) (64,1.216e-03) (128,6.082e-04)
                    };
                \addlegendentry{$e_{L^\infty}$}
                \addplot[color=black] coordinates {
                        (16,2.5e-3) (32, 1.25e-3)
                    };
                \node[anchor=north, rotate=-21, font=\footnotesize] at (axis cs:23,1.9e-3) {slope $1$};
                \addplot[color=black] coordinates {
                        (16, 6.5e-3) (32, 3.25e-3)
                    };
                \node[anchor=south, rotate=-21, font=\footnotesize] at (axis cs:23,4.25e-3) {slope $1$};
            \end{loglogaxis}
        \end{tikzpicture}
        \caption{Relative errors vs. number of time steps (in log-log scale).}
        \label{fig:error_vs_time_steps}
    \end{subfigure}
    \caption{Advection-diffusion equation from \cref{sec:advection_diffusion_experiments}: convergence of the error with respect to the number of time steps $n_t$ in dimension $d=3$.}
    \label{fig:advection_diffusion_convergence}
\end{figure}

To test the convergence of the scheme as $\dt$ decreases,
we select dimension $d=3$ to avoid
unnecessary computational costs but still
have a reasonably high-dimensional problem.
According to \cref{rmk:higher_cpu_time},
we choose to use hyperparameters
corresponding to a higher accuracy,
at the cost of computation time.
\rall{This configuration is used to ensure a
    good spatial accuracy, thus making it possible to
    observe the convergence in time.}
This leads us to taking $\ell = [30, 30]$ layers,
$N_e = \num{400}$ epochs for the initialization,
$N_e = \num{15}$ epochs per time step,
and $N_c = \num{100000}$ collocation points.
We take the final time $T=1$, and we vary $\dt$
by running the scheme
with several numbers $n_t$ of time steps.
\ra{To that end, we compute the following
    $L^2$ and $L^\infty$ relative errors
    \begin{equation*}
        e_{L^2} =
        \sqrt{\frac 1 {N_c'} \sum_{k=1}^{N_c'}
            \frac{\left( u(T, x_i) - u_\text{ex}(T, x_i) \right)^2}{u_\text{ex}(T, x_i)^2}}
        \text{\qquad and \qquad}
        e_{L^\infty} =
        \max_{i\in\{1, \dots, N_c'\}}
        \frac{\left| u(T, x_i) - u_\text{ex}(T, x_i) \right|}{\left| u_\text{ex}(T, x_i) \right|},
    \end{equation*}
    for $(x_i)_{i\in\{1, \dots, N_c'\}}$ a set of $N_c'$ collocation points
    uniformly drawn in \smash{$(-1, 1)^d$}.
    In practice, we use about $10$ times as many collocation points to compute the error than to compute the approximate solution: $N_c' \approx 10 N_c$.}
Note that the true \smash{$L^2$} error would have
required a multiplication by the volume
of the domain, which is \smash{$2^d$}.
The errors are reported in
\cref{fig:advection_diffusion_convergence},
where we observe, as expected, that the error
decreases with order \rall{one}
as the number of time steps increases.
Moreover, the scheme is indeed stable
for all time steps,
even very large ones.

\subsubsection{High-dimensional periodic solution}
\label{sec:advection_diffusion_periodic}

Then, we check the error and the computation time
with respect to the dimension $d$,
making sure that they do not explode with the dimension.
We still consider the periodic solution from
\cref{sec:advection_diffusion_convergence}.
We fix the dimension-dependent final time
\begin{equation*}
    \rall{T = \frac{\log 2}{\sigma \pi^2 \langle a, a \rangle},}
\end{equation*}
which corresponds to the time at which the diffusion
halves the amplitude of the solution.
The interval $[0, T]$ is discretized into $20$ time steps.

Compared to the previous case,
we still use \rall{$\sin$} activation functions,
but the number of layers depends on the dimension $d$:
namely, we take $\ell = [7d, 7d, 7d]$.
\rall{The other hyperparameters
    are the same with respect to the dimension:
    the number of collocation points is $N_c = \num{75000}$,
    the number of epochs is $N_e = \num{250}$
    for the initialization and $N_e = \num{5}$
    for each iteration.}

\rall{This time, we also compare our neural
    semi-Lagrangian scheme to a classical one,
    also implemented on GPU
    using the \texttt{PyTorch} library.
    The classical SL scheme uses a uniform mesh
    with $n_x$ points per spatial direction.
    The characteristic curves are then computed
    to shift the mesh,
    which is then interpolated on the original mesh
    in a directionwise manner,
    using a Lagrange interpolator of degree $3$.
    Contrary to the NSL case, the error is computed
    directly at the mesh points rather than
    at randomly sampled collocation points.}

\begin{table}[!ht]
    \centering
    \rallparagraph
    \begin{tabular}{c
            cccc
            cccc}
        \toprule
        \multirow{3}{*}[-4pt]{\makecell{dimen-                                                      \\ sion $d$}}
          & \multicolumn{4}{c}{classical SL}
          & \multicolumn{4}{c}{neural SL}                                                           \\
        \cmidrule(lr){2-5}
        \cmidrule(lr){6-9}
          & \multirow{2}{*}[-3pt]{$n_x$}
          & \multirow{2}{*}[-3pt]{Time}
          & \multicolumn{2}{c}{relative error}
          & \multicolumn{2}{c}{Time}
          & \multicolumn{2}{c}{relative error}                                                      \\
        \cmidrule(lr){4-5}
        \cmidrule(lr){6-7}
        \cmidrule(lr){8-9}
          &                                    &
          & $e_{L^2}$
          & $e_{L^\infty}$
          & init.
          & iter.
          & $e_{L^2}$
          & $e_{L^\infty}$                                                                          \\
        \cmidrule(lr){1-9}
        1 & 36                                 & \qty{0.96}{\s}   & \num{1.57e-03} & \num{2.71e-03}
          & \qty{20.88}{\s}                    & \qty{8.96}{\s}   & \num{1.56e-03} & \num{2.73e-03} \\
        2 & 38                                 & \qty{2.19}{\s}   & \num{1.63e-03} & \num{2.70e-03}
          & \qty{22.02}{\s}                    & \qty{11.05}{\s}  & \num{1.54e-03} & \num{2.79e-03} \\
        3 & 40                                 & \qty{3.31}{\s}   & \num{1.56e-03} & \num{2.56e-03}
          & \qty{25.49}{\s}                    & \qty{13.96}{\s}  & \num{1.54e-03} & \num{3.40e-03} \\
        4 & 40                                 & \qty{17.09}{\s}  & \num{1.49e-03} & \num{2.47e-03}
          & \qty{31.98}{\s}                    & \qty{18.31}{\s}  & \num{1.55e-03} & \num{3.06e-03} \\
        5 & 30                                 & \qty{50.64}{\s}  & \num{9.11e-03} & \num{1.33e-02}
          & \qty{49.33}{\s}                    & \qty{27.43}{\s}  & \num{1.52e-03} & \num{6.17e-03} \\
        6 & 19                                 & \qty{110.62}{\s} & \num{3.91e-03} & \num{8.72e-03}
          & \qty{72.50}{\s}                    & \qty{42.72}{\s}  & \num{1.53e-03} & \num{1.11e-02} \\
        7 & 13                                 & \qty{158.15}{\s} & \num{1.00e-02} & \num{2.99e-02}
          & \qty{113.48}{\s}                   & \qty{59.05}{\s}  & \num{1.54e-03} & \num{6.60e-03} \\
        8 & 10                                 & \qty{271.33}{\s} & \num{1.66e-02} & \num{8.85e-02}
          & \qty{170.27}{\s}                   & \qty{84.11}{\s}  & \num{1.57e-03} & \num{1.94e-02} \\
        \bottomrule
    \end{tabular}
    \caption{\rallparagraph Advection-diffusion equation \eqref{eq:advection_diffusion_experiments} in a periodic domain from \cref{sec:advection_diffusion_periodic}: comparison between the classical and neural semi-Lagrangian schemes. Some performance metrics are reported in different dimensions. For the classical SL scheme, we report the number $n_x$ of points per spatial direction, the computation time, and the $L^2$ and $L^\infty$ errors at the final time $T$. For the neural SL scheme, we report the computation time for all the time iterations (iter.) and initialization (init.), with associated $L^2$ and $L^\infty$ errors.}
    \label{tab:advection_diffusion}
\end{table}

\begin{figure}[!ht]
    \centering
    \begin{tikzpicture}
        \begin{axis}[
                compare SL and NSL,
                name=left_plot,
                ylabel = {\makebox[1pt][c]{\rotatebox{270}{Time (\unit{\s})}}},
                ymode=log,
                legend columns=2,
                legend to name=leg:legend_SL_NSL_sine,
            ]

            \addplot[
                mark=square,
                thick,
                graph_3,
            ] coordinates {
                    (1,0.96) (2,2.19) (3,3.31) (4,17.09) (5,50.64) (6,110.62) (7,158.15) (8,271.33)
                };

            \addlegendentry{\,Classical SL scheme \;}

            \addplot[
                mark=o,
                thick,
                graph_2,
            ] coordinates {
                    (1,20.88 + 8.96) (2,22.02 + 11.05) (3,25.49 + 13.96) (4,31.98 + 18.31) (5,49.33 + 27.43) (6,72.50 + 42.72) (7,113.48 + 59.05) (8,170.27 + 84.11)
                };

            \addlegendentry{\,Neural SL scheme \;}

            \addplot [name path=top, draw=none, forget plot] coordinates {(1,63.525) (4.25,63.525)};
            \addplot [name path=bottom, draw=none, forget plot] coordinates {(1,0.96) (4.25,0.96)};
            \addplot [red, opacity=0.2] fill between[of=top and bottom];
            \addlegendentry{\,Similar $L^2$ error for SL and NSL ($d \leqslant 4$)\;}

            \addplot [name path=top, draw=none, forget plot] coordinates {(4.75,271.33) (8,271.33)};
            \addplot [name path=bottom, draw=none, forget plot] coordinates {(4.75,45) (8,45)};
            \addplot [green, opacity=0.2] fill between[of=top and bottom];
            \addlegendentry{\,Similar computation time for SL and NSL ($d \geqslant 5$)\;}

        \end{axis}

        \begin{axis}[
                compare SL and NSL,
                name=middle_plot,
                at=(left_plot.outer east),
                anchor=outer west,
                ylabel = {\makebox[1pt][c]{\rotatebox{270}{$e_{L^2}$}}},
                ymode=log,
                ymin = 1.5e-3,
                ymax = 1.3e-2,
                ytick={1e-3, 2e-3, 5e-3, 1e-2},
                yticklabels={\num{1e-3}, \num{2e-3}, \num{5e-3}, \num{1e-2}},
            ]

            \addplot[
                mark=square,
                thick,
                graph_3,
            ] coordinates {
                    (1,1.57e-3) (2,1.63e-3) (3,1.56e-3) (4,1.49e-3) (5,9.11e-3) (6,3.91e-3) (7,1.00e-2) (8,1.66e-2)
                };

            \addplot[
                mark=o,
                thick,
                graph_2,
            ] coordinates {
                    (1,1.56e-3) (2,1.54e-3) (3,1.54e-3) (4,1.55e-3) (5,1.52e-3) (6,1.53e-3) (7,1.54e-3) (8,1.57e-3)
                };

            \addplot [name path=top, draw=none] coordinates {(1,1.85e-3) (4.25,1.85e-3)};
            \addplot [name path=bottom, draw=none] coordinates {(1,1.4e-3) (4.25,1.4e-3)};
            \addplot [red, opacity=0.2] fill between[of=top and bottom];

            \addplot [name path=top, draw=none] coordinates {(4.75,1.66e-2) (8,1.66e-2)};
            \addplot [name path=bottom, draw=none] coordinates {(4.75,1.4e-3) (8,1.4e-3)};
            \addplot [green, opacity=0.2] fill between[of=top and bottom];

        \end{axis}

        \begin{axis}[
                compare SL and NSL,
                name=right_plot,
                at=(middle_plot.outer east),
                anchor=outer west,
                ylabel = {\makebox[1pt][c]{\rotatebox{270}{memory (\unit{\giga\byte})}}},
                ymin=1,
            ]

            \addplot[
                mark=square,
                thick,
                graph_3,
            ] coordinates {
                    (1,0.19) (2,0.22) (3,0.35) (4,0.57) (5,2.58) (6,4.81) (7,7.09) (8,11.48)
                };

            \addplot[
                mark=o,
                thick,
                graph_2,
            ] coordinates {
                    (1,0.79) (2,1.64) (3,2.77) (4,4.37) (5,6.37) (6,8.82) (7,11.56) (8,14.88)
                };

            \addplot [name path=top, draw=none] coordinates {(1,5.37) (4.25,5.37)};
            \addplot [name path=bottom, draw=none] coordinates {(1,0.19) (4.25,0.19)};
            \addplot [red, opacity=0.2] fill between[of=top and bottom];

            \addplot [name path=top, draw=none] coordinates {(4.75,14.88) (8,14.88)};
            \addplot [name path=bottom, draw=none] coordinates {(4.75,2.25) (8,2.25)};
            \addplot [green, opacity=0.2] fill between[of=top and bottom];

        \end{axis}

        \node[yshift=-20pt] at (middle_plot.outer south) {\pgfplotslegendfromname{leg:legend_SL_NSL_sine}};

    \end{tikzpicture}
    \caption{\rallparagraph Advection-diffusion equation \eqref{eq:advection_diffusion_experiments} in a periodic domain from \cref{sec:advection_diffusion_periodic}: comparison between the classical (green line) and neural (orange line) semi-Lagrangian schemes. Some performance metrics (from left to right: computation time, $L^2$ error, and memory consumption) are displayed with respect to the dimension $d$ for both schemes. The faded red and green areas represent the regions where the classical SL scheme has roughly the same $L^2$ error and computation time as the NSL scheme, respectively.}
    \label{fig:advection_diffusion_periodic_sine_comparison}
\end{figure}
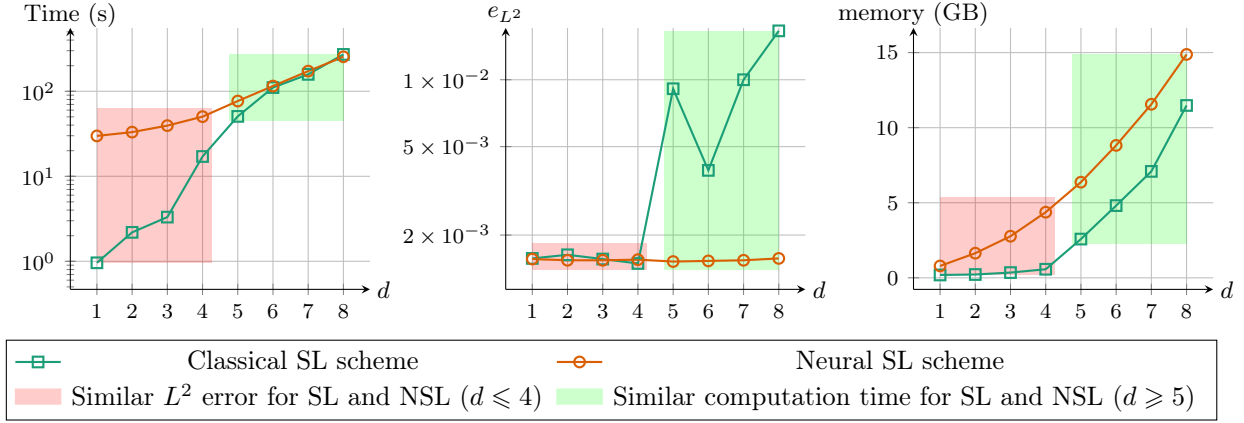

{\rallparagraph The results are reported in \cref{tab:advection_diffusion,fig:advection_diffusion_periodic_sine_comparison}.
To determine the number $n_x$ of space points per direction
in the classical SL scheme,
we have used the following rule:
\begin{itemize}[nosep]
    \item for lower dimensions ($d \leqslant 4$), choose $n_x$ to have roughly the same $L^2$ error as the NSL scheme (represented by a red area on \cref{fig:advection_diffusion_periodic_sine_comparison});
    \item for higher dimensions ($d \geqslant 5$), choose $n_x$ to have roughly the same computation time as the NSL scheme (represented by a green area on \cref{fig:advection_diffusion_periodic_sine_comparison}).
\end{itemize}
This leads to the memory consumption and computation time
being comparable for high dimensions,
where the classical scheme struggles the most.
Otherwise, reaching similar errors in high dimensions
would have been prohibitively expensive in terms of computation time,
or even unreachable in terms of
the \qty{64}{\giga\byte} memory of our~GPU.

With this configuration,
we observe that the $L^2$ error remains
roughly constant with the dimension for the~NSL scheme,
while it sharply increases for the classical SL scheme
when higher dimensions are considered.
For both~$L^2$ and $L^\infty$ errors
with roughly the same resource usage
(computation time and memory),
we observe that the NSL scheme
outperforms the classical SL scheme
by a factor of around $10$ in dimensions $d\geqslant5$.
Moreover, for dimensions $d \leq 4$,
we note that the NSL scheme
has an error similar to that of a classical SL scheme with
around~$40$ points per spatial direction.
Reaching this level of refinement
becomes untenable in higher dimensions.
These findings are in line with
e.g.~\cite{GroKomLatSch2024},
where PINNs were proven to be inferior
to standard finite element methods
in low dimensions ($d \leqslant 3$)}.

\subsubsection{High-dimensional Gaussian solution}
\label{sec:advection_diffusion_gaussian}

We now consider a more localized solution.
Namely, in a non-periodic domain,
the exact solution is given by
\begin{equation*}
    u_\text{ex}(t, x) =
    1 + \sqrt{\frac {\det \Sigma(0)} {\det \Sigma(t)}}
    \exp \left(
    - \frac 1 2 \,
    \bigl(x - M(t)\bigr)^\intercal \,
    \Sigma(t)^{-1} \,
    \bigl(x - M(t)\bigr)
    \right),
\end{equation*}
where the mean $M$ is defined by
$\rallparagraph M(t) = a t$
and the covariance $\Sigma$ satisfies
$\Sigma(t) = \Sigma(0) + 2 \sigma t I_d$,
with $I_d$ the identity matrix in dimension $d$.
The symmetric positive-definite initial covariance
$\Sigma(0)$ satisfies
\begin{equation*}
    \rallparagraph
    \Sigma(0) = 2 d \, \sigma_0^2 \,
    \bigl( 2 d \, I_d + \tilde {\Sigma} \bigr)
    \text{, \qquad with \quad} \sigma_0 = 0.05
    \text{\quad and \quad} \tilde {\Sigma}_{ij} = d - |i - j|,
\end{equation*}
and we take the diffusion coefficient $\rallparagraph\sigma = \num{5e-2}$.
This solution corresponds to a nonsymmetric Gaussian pulse
being advected and diffused
along the diagonal of the domain.
To have a good representation of this solution,
we take the space domain $\rallparagraph\Omega = (-3, 3)^d$,
and the final time is $\rallparagraph T = 1$.
\rall{Note that, in this section,
    we only consider dimensions $d > 1$,
    to have a true non-separable and anisotropic solution.}

\rb{Contrary to an isotropic Gaussian bump, the tensor rank of this exact solution is not one. For example, in dimension~2, its rank is around 15 to 20 (depending on time). In this case, an approximation with precision $10^{-3}$ requires a rank of $5$ or $6$. So, like in the previous case, a low rank method is not well-suited to this test case. Indeed, the approximate rank value combined with the dimension, which can go up to 8, means that the number of 1D interpolations that must be performed at each step of the time-stepping algorithm becomes significant, making the whole computation much more expensive.}

This time, we use both the natural gradient preconditioning
and the adaptive sampling.
For the latter, we sample using the gradient of the solution,
like in \eqref{eq:3D_cylinder_advection_vector_field},
and we set $\sigma_1 = \num{20}$, $\sigma_2 = \num{100}$ and $\sigma_3 = \num{5000}$.
Regarding the other hyperparameters,
we take $\rallparagraph 5$ time steps,
regularized hat activation functions
as defined in~\eqref{eq:regularized_hat_function},
and the number of layers is $\ell = [7d, 7d, 7d]$.
The number of epochs is $\rallparagraph N_e = \num{150}$ for the initialization
and $\rallparagraph N_e = \num{15}$ for each iteration,
and we take $\rallparagraph N_c = \num{75000}$ collocation points.

\begin{table}[!ht]
    \centering
    \rallparagraph
    \begin{tabular}{c
            cccc
            cccc}
        \toprule
        \multirow{3}{*}[-4pt]{\makecell{dimen-                                                      \\ sion $d$}}
          & \multicolumn{4}{c}{classical SL}
          & \multicolumn{4}{c}{neural SL}                                                           \\
        \cmidrule(lr){2-5}
        \cmidrule(lr){6-9}
          & \multirow{2}{*}[-3pt]{$n_x$}
          & \multirow{2}{*}[-3pt]{Time}
          & \multicolumn{2}{c}{relative error}
          & \multicolumn{2}{c}{Time}
          & \multicolumn{2}{c}{relative error}                                                      \\
        \cmidrule(lr){4-5}
        \cmidrule(lr){6-7}
        \cmidrule(lr){8-9}
          &                                    &
          & $e_{L^2}$
          & $e_{L^\infty}$
          & init.
          & iter.
          & $e_{L^2}$
          & $e_{L^\infty}$                                                                          \\
        \cmidrule(lr){1-9}
        2 & 64                                 & \qty{0.86}{\s}   & \num{9.39e-04} & \num{1.15e-02}
          & \qty{21.27}{\s}                    & \qty{10.84}{\s}  & \num{8.75e-04} & \num{1.14e-02} \\
        3 & 64                                 & \qty{1.27}{\s}   & \num{3.08e-04} & \num{1.09e-02}
          & \qty{23.97}{\s}                    & \qty{14.23}{\s}  & \num{3.44e-04} & \num{1.10e-02} \\
        4 & 64                                 & \qty{27.16}{\s}  & \num{1.15e-04} & \num{9.42e-03}
          & \qty{29.76}{\s}                    & \qty{17.86}{\s}  & \num{1.47e-04} & \num{9.26e-03} \\
        5 & 32                                 & \qty{62.91}{\s}  & \num{1.03e-04} & \num{1.58e-02}
          & \qty{41.45}{\s}                    & \qty{25.12}{\s}  & \num{7.35e-05} & \num{8.02e-03} \\
        6 & 23                                 & \qty{85.59}{\s}  & \num{1.12e-04} & \num{1.75e-02}
          & \qty{55.86}{\s}                    & \qty{33.86}{\s}  & \num{3.93e-05} & \num{4.09e-03} \\
        7 & 15                                 & \qty{106.87}{\s} & \num{1.19e-04} & \num{3.82e-02}
          & \qty{79.70}{\s}                    & \qty{47.69}{\s}  & \num{3.62e-05} & \num{5.34e-03} \\
        8 & 11                                 & \qty{139.95}{\s} & \num{2.47e-04} & \num{1.08e-01}
          & \qty{109.83}{\s}                   & \qty{64.34}{\s}  & \num{3.69e-05} & \num{3.70e-03} \\
        \bottomrule
    \end{tabular}

    \caption{\rallparagraph Advection-diffusion equation \eqref{eq:advection_diffusion_experiments} of a Gaussian pulse \cref{sec:advection_diffusion_gaussian}: comparison between the classical and neural semi-Lagrangian schemes. Some performance metrics are reported in different dimensions. For the classical SL scheme, we report the number $n_x$ of points per spatial direction, the computation time, and the $L^2$ and $L^\infty$ errors at the final time $T$. For the neural SL scheme, we report the computation time for all the time iterations (iter.) and initialization (init.), with associated $L^2$ and $L^\infty$ errors.}
    \label{tab:advection_diffusion_Gaussian}
\end{table}

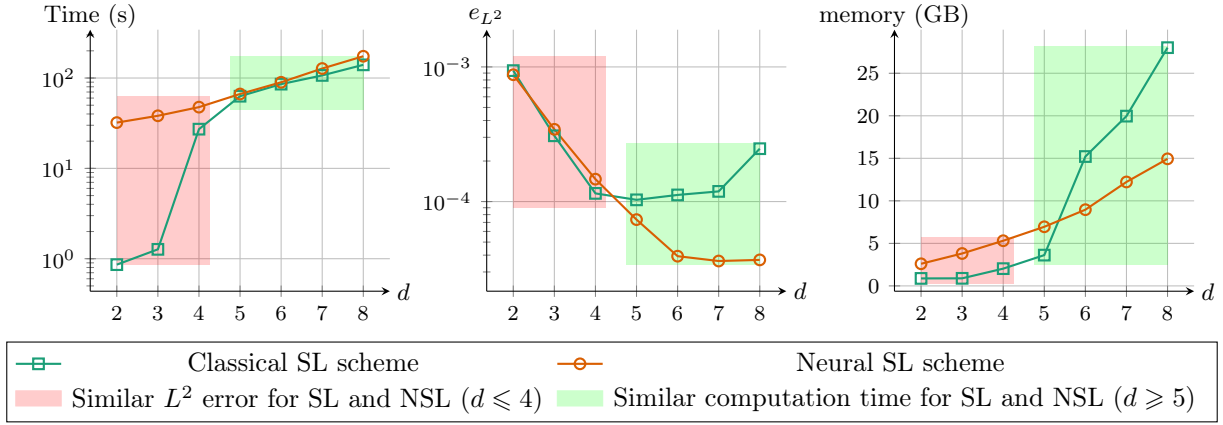
\begin{figure}[!ht]
    \centering
    \begin{tikzpicture}
        \begin{axis}[
                compare SL and NSL,
                name=left_plot,
                ylabel = {\makebox[1pt][c]{\rotatebox{270}{Time (\unit{\s})}}},
                ymode=log,
                xtick={2, 3, 4, 5, 6, 7, 8},
                legend columns=2,
                legend to name=leg:legend_SL_NSL_gaussian,
            ]

            \addplot[
                mark=square,
                thick,
                graph_3,
            ] coordinates {
                    (2,0.86) (3,1.27) (4,27.16) (5,62.91) (6,85.59) (7,106.87) (8,139.95)
                };

            \addlegendentry{\,Classical SL scheme \;}

            \addplot[
                mark=o,
                thick,
                graph_2,
            ] coordinates {
                    (2,21.27 + 10.84) (3,23.97 + 14.23) (4,29.76 + 17.86) (5,41.45 + 25.12) (6,55.86 + 33.86) (7,79.70 + 47.69) (8,109.83 + 64.34)
                };

            \addlegendentry{\,Neural SL scheme \;}

            \addplot [name path=top, draw=none, forget plot] coordinates {(2,63.525) (4.25,63.525)};
            \addplot [name path=bottom, draw=none, forget plot] coordinates {(2,0.86) (4.25,0.86)};
            \addplot [red, opacity=0.2] fill between[of=top and bottom];
            \addlegendentry{\,Similar $L^2$ error for SL and NSL ($d \leqslant 4$)\;}

            \addplot [name path=top, draw=none, forget plot] coordinates {(4.75,109.83 + 64.34) (8,109.83 + 64.34)};
            \addplot [name path=bottom, draw=none, forget plot] coordinates {(4.75,45) (8,45)};
            \addplot [green, opacity=0.2] fill between[of=top and bottom];
            \addlegendentry{\,Similar computation time for SL and NSL ($d \geqslant 5$)\;}

        \end{axis}

        \begin{axis}[
                compare SL and NSL,
                name=middle_plot,
                at=(left_plot.outer east),
                anchor=outer west,
                ylabel = {\makebox[1pt][c]{\rotatebox{270}{$e_{L^2}$}}},
                ymode=log,
                xtick={2, 3, 4, 5, 6, 7, 8},
            ]

            \addplot[
                mark=square,
                thick,
                graph_3,
            ] coordinates {
                    (2,9.39e-4) (3,3.08e-4) (4,1.15e-4) (5,1.03e-4) (6,1.12e-4) (7,1.19e-4) (8,2.47e-4)
                };

            \addplot[
                mark=o,
                thick,
                graph_2,
            ] coordinates {
                    (2,8.749e-4) (3,3.440e-4) (4,1.467e-4) (5,7.345e-5) (6,3.930e-5) (7,3.618e-5) (8,3.688e-5)
                };

            \addplot [name path=top, draw=none] coordinates {(2,1.2e-3) (4.25,1.2e-3)};
            \addplot [name path=bottom, draw=none] coordinates {(2,9e-5) (4.25,9e-5)};
            \addplot [red, opacity=0.2] fill between[of=top and bottom];

            \addplot [name path=top, draw=none] coordinates {(4.75,2.7e-4) (8,2.7e-4)};
            \addplot [name path=bottom, draw=none] coordinates {(4.75,3.4e-5) (8,3.4e-5)};
            \addplot [green, opacity=0.2] fill between[of=top and bottom];

        \end{axis}

        \begin{axis}[
                compare SL and NSL,
                name=right_plot,
                at=(middle_plot.outer east),
                anchor=outer west,
                ylabel = {\makebox[1pt][c]{\rotatebox{270}{memory (\unit{\giga\byte})}}},
                ymin=2.5,
                ymax=27,
                ytick={0,5,10,15,20,25},
                xtick={2, 3, 4, 5, 6, 7, 8},
            ]

            \addplot[
                mark=square,
                thick,
                graph_3,
            ] coordinates {
                    (2,0.88) (3,0.89) (4,2.03) (5,3.62) (6,15.21) (7,19.96) (8,28.01)
                };

            \addplot[
                mark=o,
                thick,
                graph_2,
            ] coordinates {
                    (2,2.60) (3,3.81) (4,5.31) (5,6.95) (6,8.96) (7,12.22) (8,14.94)
                };

            \addplot [name path=top, draw=none] coordinates {(2,5.7) (4.25,5.7)};
            \addplot [name path=bottom, draw=none] coordinates {(2,0.19) (4.25,0.19)};
            \addplot [red, opacity=0.2] fill between[of=top and bottom];

            \addplot [name path=top, draw=none] coordinates {(4.75,28.2) (8,28.2)};
            \addplot [name path=bottom, draw=none] coordinates {(4.75,2.5) (8,2.5)};
            \addplot [green, opacity=0.2] fill between[of=top and bottom];

        \end{axis}

        \node[yshift=-20pt] at (middle_plot.outer south) {\pgfplotslegendfromname{leg:legend_SL_NSL_gaussian}};

    \end{tikzpicture}
    \caption{\rallparagraph Advection-diffusion equation \eqref{eq:advection_diffusion_experiments} of a Gaussian pulse \cref{sec:advection_diffusion_gaussian}: comparison between the classical (green line) and neural (orange line) semi-Lagrangian schemes. Some performance metrics (from left to right: computation time, $L^2$ error, and memory consumption) are displayed with respect to the dimension $d$ for both schemes. The faded red and green areas represent the regions where the classical SL scheme has roughly the same $L^2$ error and computation time as the NSL scheme, respectively.}
    \label{fig:advection_diffusion_gaussian_comparison}
\end{figure}

{\rallparagraph The results are reported in
\cref{tab:advection_diffusion_Gaussian,fig:advection_diffusion_gaussian_comparison}.
Similar comments as in the previous case apply.
Namely, we observe that the NSL scheme
is much more efficient than the classical SL scheme
in high dimensions.
This effect is further increased by the fact
that the NSL scheme is roughly equivalent
to a classical SL scheme with $n_x = 64$
points per spatial direction,
which becomes even more prohibitive in high dimensions
than the previous case.
Therefore, the NSL scheme reaches
much lower $L^2$ and $L^\infty$ errors
than the classical SL scheme,
up to a factor of $50$ on the $L^\infty$ error
in dimension $d=8$ for a similar computation time.
Moreover, the memory consumption of the classical SL scheme
increases much faster in this case,
mostly when computing the initial condition.
Indeed, taking $n_x=12$ points per spatial direction
in dimension $d=8$ leads to almost
\qty{64}{\giga\byte} of memory,
almost filling the GPU and taking
over five minutes to compute the approximate solution.

\subsubsection{Discussion}
\label{sec:advection_diffusion_discussion}

The results of the two previous sections allow us to draw some conclusions
regarding the comparison between the classical and neural semi-Lagrangian schemes.

\paragraph*{Memory consumption}

The analysis of memory complexity is particularly challenging in our setting. In the classical SL method, memory usage scales predictably with the number of dofs, which grows exponentially with the dimension. For example, using $N$ dofs (i.e., grid points) per dimension results in $N^d$ dofs in $d$ dimensions, quickly making storage intractable.
In contrast, our approach relies on small neural networks (here, with between $161$ dofs for $d=1$ and $6776$ dofs for $d=8$), leading to negligible storage cost even if each time step is stored. However, the runtime memory footprint and number of operations per dof can still be significant, primarily due to the training process. This makes training significantly more memory-intensive than inference. Estimating the theoretical memory usage for neural networks is very challenging, and thus we only provided an empirical comparison.

\paragraph*{Computation time}

For both considered test cases,
we observed that the neural semi-Lagrangian scheme
is more efficient than the classical one
for dimensions $d \geqslant 5$,
leading to lower errors for the same computation time.
Moreover, even if one was willing to pay an extremely high computational cost,
the classical scheme would probably not be able to reach
the same accuracy as the NSL one
due to the aforementioned memory constraints.
Indeed, to decrease the error by a factor of $8$,
one needs to increase the number of points
per spatial direction by a factor of $2$
(since the classical scheme is third-order accurate in space,
and provided the time-stepping is accurate enough).
In dimension $d=8$, this means increasing the total number of dofs
from~$11^8$ to $22^8 = 10^8 N^8 \approx \num{5e10}$,
which is clearly unreachable on a single GPU.
Even if it were reachable,
the computation time would be multiplied by a factor of roughly $2^8$,
leading to a total computation time of around~\qty{10}{\hour}.

\paragraph*{Applicability of the neural semi-Lagrangian method}

As shown in our comparisons, we are able to solve problems in high dimensions, where, \rabis{as expected}, the standard semi-Lagrangian method become intractable due to the explosion in dofs, and associated memory and computational costs. This confirms that, while our method involves a higher per-iteration memory and computational cost, its global complexity in high dimensions is often far more manageable thanks to its mesh-free nature.}

{\rallparagraph

\subsection{Extension to the Vlasov-Poisson equations}
\label{sec:vlasov_poisson}

As a last experiment, we propose to extend the NSL method
to a coupled nonlinear system: the Vlasov-Poisson equations.
In this experiment, the goal is not to solve these equations efficiently.
Rather, it is to show that the proposed method
is able to handle such a system,
despite it being naturally well-suited to linear problems.
We first describe the system and how to solve it
in \cref{sec:vlasov_poisson_equations}.
Then, we present an approximation of the well-known
bump-on-tail instability in \cref{sec:vlasov_poisson_validation}.

\subsubsection{Governing equations and solving strategy}
\label{sec:vlasov_poisson_equations}

The Vlasov-Poisson equations (see e.g.~\cite{KorReuRam2019})
are a system of kinetic equations,
governed in one space and one velocity dimension by
\begin{equation*}
    \begin{dcases}
        \partial_t u(t,x,v) +
        v \, \partial_x u(t,x,v) +
        E(t,x) \, \partial_v u(t,x,v) = 0,                  &
        \forall \, t \in [0, T],
        \forall \, x \in \Omega_x,
        \forall \, v \in \Omega_v, \vphantom{\int{\Omega_v}}  \\
        \Delta \Psi(t,x) =
        \int_{\Omega_v} u(t,x,v) \, dv - \langle u \rangle, &
        \forall \, t \in [0, T],
        \forall \, x \in \Omega_x,                            \\
        E(t,x) = \partial_x \Psi(t,x),                      &
        \forall \, t \in [0, T],
        \forall \, x \in \Omega_x, \vphantom{\int{\Omega_v}}
    \end{dcases}
\end{equation*}
where $t \in [0, T]$ is the time variable (with $T$ a final time),
$x \in \Omega_x$ is the space variable,
and $v \in \Omega_v$ is the velocity variable.
Periodic boundary conditions are prescribed everywhere.
The function~$u$ represents the distribution function,
\smash{$\langle u \rangle = \frac {1}{|\Omega_x|} \int_{\Omega_x} \int_{\Omega_v} u(t,x,v) \, dv \, dx$} its average,
and~$E$ is an electric field, with~$\Psi$ its potential.
This equation comprises the Vlasov equation from \cref{sec:2D_fake_Vlasov}
and a Poisson equation giving $E$, the $v$-component of the advection field.

\paragraph{Classical semi-Lagrangian method}
Because of this coupling, classical semi-Lagrangian strategies are not directly applicable.
Usual, mesh-based numerical methods rely on an operator splitting, alternating
between advection in the $x$-direction and advection in the~$v$-direction,
effectively dividing the dimension by two.
Strang operator splitting makes it possible to achieve second-order accuracy in time.
To describe the algorithm, assume that values $u^n$ at time $t^n$ are known on the grid.
We seek updated values $u^{n+1}$ at time $t^{n+1}$, still on the grid.
The algorithm relies on the following operations, presented here in semi-discrete form:
\begin{enumerate}[label=\arabic*.]
    \item solve $\partial_t u + v \, \partial_x u = 0$
          with initial data $u^n$, from $t^n$ to $t^n + 0.5 \, \dt$,
          to get intermediate values $u^*$;
    \item solve \smash{$\Delta \Psi^* = \displaystyle \int_{\Omega_v} u^* \, dv - \langle u^* \rangle$}
          to get an intermediate electric field $E^* = \partial_x \Psi^*$;
    \item solve $\partial_t u + E^* \, \partial_v u = 0$
          with initial data $u^*$, from $t^n$ to $t^{n+1}$,
          to get new intermediate values $u^{**}$;
    \item solve $\partial_t u + v \, \partial_x u = 0$
          with initial data $u^{**}$, from $t^n + 0.5 \, \dt$ to $t^{n+1}$,
          to get the updated values $u^{n+1}$.
\end{enumerate}

\paragraph{Neural semi-Lagrangian method}
Similarly, our neural semi-Lagrangian method has to be adapted to this coupling.
In our case, there is no need to split the two advection operators to lower the dimension,
since our meshless algorithm can handle higher dimensions than classical, mesh-based ones.
In our case, the algorithm uses a predictor-corrector method to reach second-order accuracy in time.
It reads as follows:
\begin{enumerate}[label=\arabic*.]
    \item solve \smash{$\Delta \Psi^n = \displaystyle \int_{\Omega_v} u^n \, dv - \langle u^n \rangle$}
          to get the electric field $E^n = \partial_x \Psi^n$ at time $t^n$;
    \item solve $\partial_t u + v \, \partial_x u + E^n \, \partial_v u = 0$
          with initial data $u^n$, from $t^n$ to $t^{n+1}$,
          to get a prediction $u^*$;
    \item solve \smash{$\Delta \Psi^* = \displaystyle \int_{\Omega_v} u^* \, dv - \langle u^* \rangle$}
          to get a predicted electric field $E^* = \partial_x \Psi^*$,
          and set \smash{$\bar E = \dfrac{E^n + E^*}{2}$};
    \item solve $\partial_t u + v \, \partial_x u + \bar E \, \partial_v u = 0$
          with initial data $u^n$, from $t^n$ to $t^{n+1}$,
          to get the updated values $u^{n+1}$.
\end{enumerate}
The second and fourth steps are
the same as in the previous sections,
but the first and third ones
involve solving a Poisson equation.
Naturally, a PINN is well-suited for this task,
since it gives the potential $\Psi^n$ as a function that can be
automatically differentiated to get the electric field $E^n$.
Concretely, to solve the Poisson equation,
the PINN uses hard-constrained homogeneous
Dirichlet boundary conditions on $\Psi$,
to fix the free constant in the Poisson equation
with periodic boundary conditions.
Classical schemes (based on e.g.\ a discrete Fourier transform)
can also be employed to solve the Poisson equation
in dimension $d=1$.
We still decide to use a PINN for this task,
as a preparation for higher dimensions,
and since both approaches gave the same results.
All in all, our algorithm consists in alternating a PINN
with the NSL method at each time step.
More complex algorithms, potentially higher-order accurate in time,
could also be employed.
We do not investigate this direction further here,
since our goal is to provide a proof of concept
for such coupled nonlinear systems.

\subsubsection{Validation: the bump-on-tail instability}
\label{sec:vlasov_poisson_validation}

The bump-on-tail instability consists in perturbing
an initial equilibrium distribution function,
thus creating vortices that grow over time.
In lower dimensions, it is well-solved by classical schemes.
The space domain is~$\Omega_x = [0, 10 \pi]$
and the velocity domain is $\Omega_v = [-9, 9]$,
while the final time is $T = 32$.
The initial condition is given by
\begin{equation*}
    u(0,x,v) = \left(
    \frac{1 - \varepsilon_\text{bot}}{\sqrt{2 \pi}} e^{-\frac{v^2}{2}}
    + \frac{2 \, \varepsilon_\text{bot}}{\sqrt{2 \pi T_\text{bot}}} e^{-\frac{(v - v_\text{bot})^2}{2 T_\text{bot}}}
    \right)
    \left(1 + \varepsilon_\text{pert} \cos(k_x x)\right),
\end{equation*}
where $\varepsilon_\text{bot} = 0.1$,
$v_\text{bot} = 3.8$,
$T_\text{bot} = 0.2$,
$\varepsilon_\text{pert} = 0.03$,
and $k_x = 0.4$.

Since $\Psi$ is a PINN approximating the periodic solution to a one-dimensional equation,
we use the rectangle method for integrating the loss function, for added efficiency.
Moreover, the integral over $\Omega_v$ in the Poisson equation
is also approximated using the rectangle method (using 196 points)
since $u$ is a periodic function.
In addition, to keep the distribution function $u$ nonnegative,
we add a nonnegative activation function
(ReLU\textsuperscript{4}) on the last layer.
We take $\dt = 0.2$, leading to $160$ time steps.
With this choice of hyperparameters and time discretization,
the computation time is around four hours.
This is prohibitively long for such a two-dimensional test case
(the reference solution runs in two hours on a single CPU core),
but we recall that this simulation is merely a proof of concept,
intended to show that the proposed NSL method
is able to handle coupled and nonlinear problems.

\begin{figure}[!ht]
    \centering
    \includegraphics[width=\textwidth]{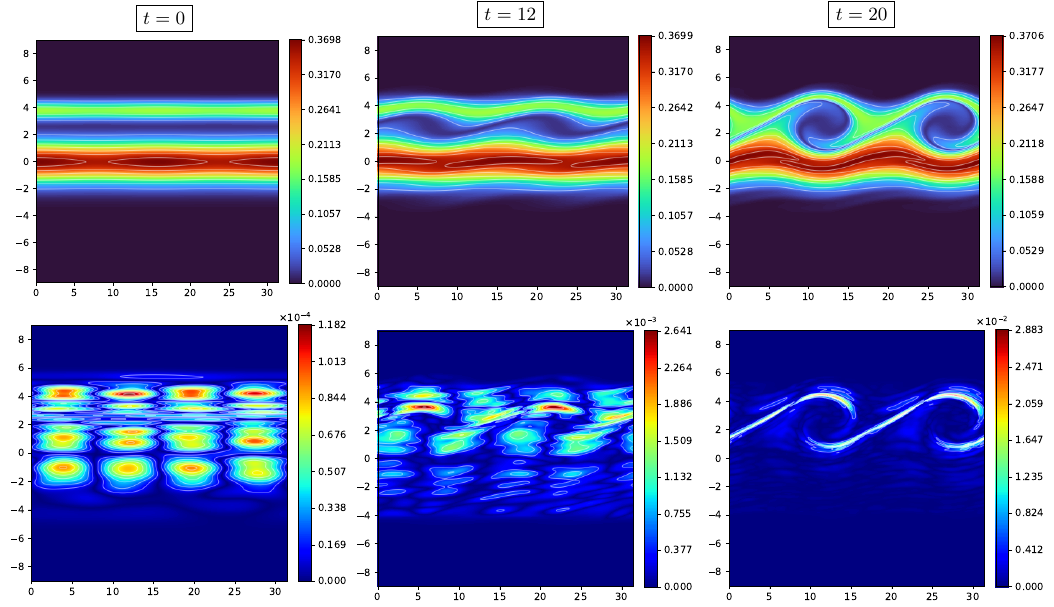}
    \caption{\rallparagraph Bump-on-tail instability from \cref{sec:vlasov_poisson_validation}: evolution of the distribution function $u$ at different times (from left to right: $t = 0$, $t = 12$ and $t = 20$) The top panels display $u$ (in color) and its contours (in white). The bottom panels display the pointwise error between $u$ and the reference solution.}
    \label{fig:vlasov_poisson_bump_on_tail_1}
\end{figure}

\begin{figure}[!ht]
    \centering
    \includegraphics[width=\textwidth]{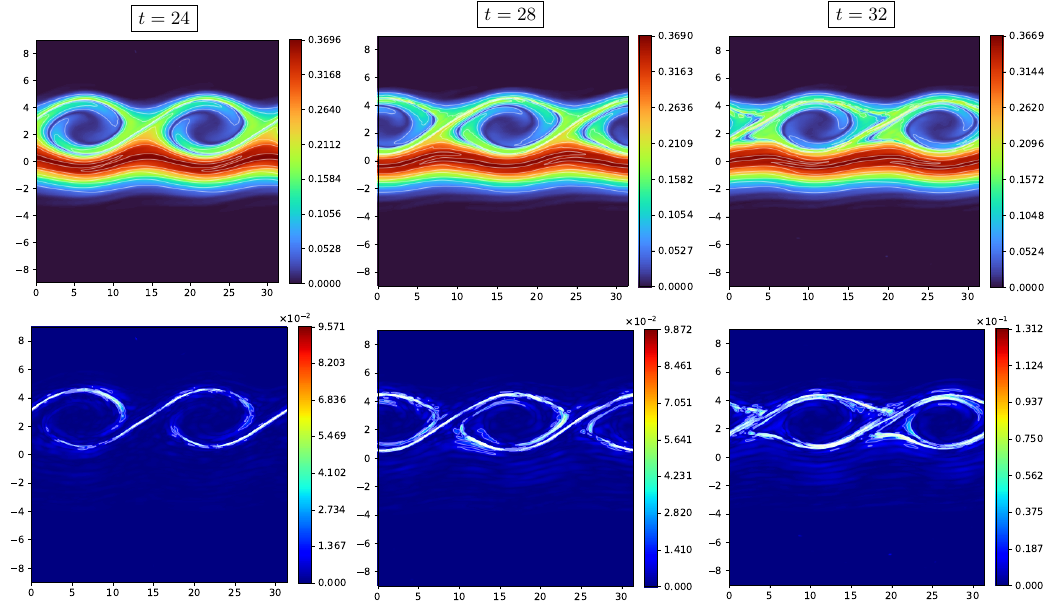}
    \caption{\rallparagraph Bump-on-tail instability from \cref{sec:vlasov_poisson_validation}: evolution of the distribution function $u$ at different times (from left to right: $t = 24$, $t = 28$ and $t = 32$) The top panels display $u$ (in color) and its contours (in white). The bottom panels display the pointwise error between $u$ and the reference solution.}
    \label{fig:vlasov_poisson_bump_on_tail_2}
\end{figure}

\begin{table}[!ht]
    \centering
    \rallparagraph
    \begin{tabular}{ccccccc}
        \toprule
         & $t=0$            & $t=12$           & $t=20$
         & $t=24$           & $t=28$           & $t=32$           \\
        \cmidrule(lr){1-7}
        relative $L^2$ error
         & $\num{2.25e-04}$ & $\num{3.41e-03}$ & $\num{1.49e-02}$
         & $\num{2.56e-02}$ & $\num{3.97e-02}$ & $\num{6.23e-02}$ \\
        relative $L^\infty$ error
         & $\num{3.20e-04}$ & $\num{7.14e-03}$ & $\num{7.79e-02}$
         & $\num{2.59e-01}$ & $\num{2.67e-01}$ & $\num{3.55e-01}$ \\
        \bottomrule
    \end{tabular}
    \caption{\rallparagraph Bump-on-tail instability from \cref{sec:vlasov_poisson_validation}: relative $L^2$ and $L^\infty$ errors at different times.}
    \label{tab:errors_wrt_time_BOT}
\end{table}

The results are displayed in
\cref{fig:vlasov_poisson_bump_on_tail_1} for $t \in \{0, 12, 20\}$ and in
\cref{fig:vlasov_poisson_bump_on_tail_2} for $t \in \{24, 28, 32\}$.
The reference solution, used to compute the error,
is obtained using the algorithm described in \cref{sec:vlasov_poisson_equations}
on a very fine space-time grid of $2048 \times 1024$ points and $\num{8000}$ time steps
(leading to $\dt = \num{5e-3}$).
Relative $L^2$ and $L^\infty$ errors are reported in \cref{tab:errors_wrt_time_BOT}.
In all cases, we observe a good agreement with the reference solution.
On the one hand, for times $t \leqslant 20$,
the relative $L^2$ error remains below roughly $1\%$,
and increases slowly over time.
Of course, as time advances, the instability grows and becomes harder to approximate,
and thus the relative~$L^2$ error also increases over time,
up to around $6\%$ at $t = 32$.
On the other hand, the $L^\infty$ error is naturally larger than usual,
and may not be a very relevant measure of error.
This is due to two main reasons:
first, the time grid of the NSL method is quite coarse;
second, even with a finer grid, there would be an $\mathcal{O}(\dt^2)$ dispersion error.
Indeed, the classical and neural methods do not use the same time splitting,
even though both are second-order accurate in time.
This leads to a slight shift between the two solutions,
exacerbated by the long-time simulation.
We nevertheless included the pointwise error in the figures,
as it provides a more detailed view of the discrepancies.
We also remark that the optimization problems become harder and harder to solve with time.
Indeed, at the beginning of the simulation,
the optimization problems are solved up to a mean squared error of $10^{-10}$
for the problem on $u$, and $10^{-12}$ for the problem on $\Psi$.
This increases to $10^{-7}$ and $\num{5e-9}$, respectively, at the final time.}

\section{Conclusion}
\label{sec:conclusion}

In this work, we have introduced a neural version
of the well-known semi-Lagrangian method.
This method has the advantage of being able to solve
advection-diffusion in high dimensions,
without any time step requirements for stability.
It is based on modelling the approximate solution
using a neural network,
and then using nonlinear optimization
to solve the resulting backwards transport equation.
It is an iterative process creating a sequence of neural networks
$(u_{\theta^n})_n$ approximating the solution at time $t^n = n \dt$.
This sequence is constructed using a simple algorithm,
whose steps each consist in solving a nonlinear optimization problem.
The rough algorithm is as follows:
\begin{itemize}[nosep]
    \item \emph{step $0$}: fit $u_{\theta^0}$ to the initial condition $u_0$;
    \item \emph{step $n+1$}: given the network $u_{\theta^n}$ at time $t^n$,
          fit $u_{\theta^{n+1}}$ to
          $u_{\theta^n}(\widetilde{\mathcal{X}}(t^n; t^{n+1}, \cdot))$,
          where $\widetilde{\mathcal{X}}(t^n; t^{n+1}, \cdot)$
          is an approximation of the backwards characteristic curve.
\end{itemize}
Numerical results show the performance of the method,
especially compared to (albeit more general)
methods from the literature.
Namely, we show that a good accuracy is retained in
high dimensions,
for a comparatively low computational cost.

To improve Vlasov simulations,
the next step consists in exploring
improved network architectures,
e.g.\ using Fourier features~\cite{TanSriMilFriRagSinRamBarNg2020},
PirateNets~\cite{WanLiChePer2024},
or Kolmogorov-Arnold Networks (KAN, see~\cite{LiuWan2025}).
Domain decomposition and parallelism will also be employed
to optimize computational efficiency,
in the spirit of e.g.~\cite{DolHeiMisMos2023}.
We will also propose a better adaptive sampling,
to further reduce the computational cost,
especially at high dimensions.
Work like~\cite{MaoMen2023} could be adapted.
\rall{We already provided a proof of concept on the Vlasov-Poisson equations,
    but our neural strategy will also be extended to more complex systems.
    For instance, Vlasov-Poisson simulations will be improved,
    and Vlasov-Maxwell~\cite{LiuCaiLapCao2021} or gyrokinetic~\cite{Gra2016} models will be tackled,
    going towards high-dimensional
    applications in plasma physics using advanced splitting methods like in~\cite{CroEinFao2015}.}
\ra{Our method could also be adapted to further
    improve the results obtained on the
    level-set deformation,
    by introducing strong or weak constraints
    on the volume preservation.}

\section*{Acknowledgments}

As part of the ``France 2030'' initiative, this work has benefited from a national grant managed by the French National Research Agency (Agence Nationale de la Recherche) attributed to the Exa-MA project of the NumPEx PEPR program, under the reference ANR-22-EXNU-0002.
This work was also supported by the French National Research Agency through projects
ANR-23-PEIA-0004 (PDE-AI, all four authors),
ANR-21-CE46-0014 (Milk, E. F. and L. N.),
and ANR-22-CE25-0017 (OptiTrust, V. M.-D.).

\bibliographystyle{plain}
\bibliography{references}

\begin{thebibliography}{10}

\bibitem{AmaDou1998}
S.~Amari and S.~C. Douglas.
\newblock Why natural gradient?
\newblock In {\em Proceedings of the 1998 IEEE International Conference on Acoustics, Speech and Signal Processing, ICASSP '98 (Cat. No.98CH36181)}, volume~2 of {\em ICASSP-98}, pages 1213--1216. IEEE, 1998.

\bibitem{PyTorch2}
J.~Ansel and E.~Yang et~al.
\newblock {PyTorch 2: Faster Machine Learning Through Dynamic Python Bytecode Transformation and Graph Compilation}.
\newblock In {\em Proceedings of the 29th ASPLOS, Volume 2}, volume~5 of {\em ASPLOS '24}, pages 929--947. ACM, 2024.

\bibitem{9788008}
A.~Beltran-Pulido, I.~Bilionis, and D.~Aliprantis.
\newblock {Physics-Informed Neural Networks for Solving Parametric Magnetostatic Problems}.
\newblock {\em IEEE Trans. Energy Convers.}, 37(4):2678--2689, 2022.

\bibitem{BerPeh2023}
J.~Berman and B.~Peherstorfer.
\newblock {Randomized Sparse Neural Galerkin Schemes for Solving Evolution Equations with Deep Networks}.
\newblock In {\em Thirty-seventh Conference on Neural Information Processing Systems}, 2023.

\bibitem{besse2008convergence}
N.~Besse and M.~Mehrenberger.
\newblock {Convergence of classes of high-order semi-Lagrangian schemes for the Vlasov--Poisson system}.
\newblock {\em Math. Comput.}, 77(261):93--123, 2008.

\bibitem{besse2003semi}
N.~Besse and E.~Sonnendr{\"u}cker.
\newblock {Semi-Lagrangian schemes for the Vlasov equation on an unstructured mesh of phase space}.
\newblock {\em J. Comput. Phys.}, 191(2):341--376, 2003.

\bibitem{biesek2023burgers}
V.~Biesek and P.~H. de~Almeida~Konzen.
\newblock {Burgers' PINNs with implicit Euler Transfer Learning}.
\newblock {\em Rev. Mundi Eng., Tecnol. Gest.}, 9(4), 2024.

\bibitem{BokSim2016}
O.~Bokanowski and G.~Simarmata.
\newblock {Semi-Lagrangian discontinuous Galerkin schemes for some first- and second-order partial differential equations}.
\newblock {\em ESAIM: M2AN}, 50(6):1699--1730, 2016.

\bibitem{BonCalCarFer2021}
L.~Bonaventura, E.~Calzola, E.~Carlini, and R.~Ferretti.
\newblock {Second Order Fully Semi-Lagrangian Discretizations of Advection-Diffusion-Reaction Systems}.
\newblock {\em J. Sci. Comput.}, 88(1), 2021.

\bibitem{BonFerRoc2018}
L.~Bonaventura, R.~Ferretti, and L.~Rocchi.
\newblock A fully semi-{L}agrangian discretization for the {2D} incompressible {N}avier-{S}tokes equations in the vorticity-streamfunction formulation.
\newblock {\em Appl. Math. Comput.}, 323:132--144, 2018.

\bibitem{BRUNA2024112588}
J.~Bruna, B.~Peherstorfer, and E.~Vanden-Eijnden.
\newblock Neural {G}alerkin schemes with active learning for high-dimensional evolution equations.
\newblock {\em J. Comput. Phys.}, 496:112588, 2024.

\bibitem{BuiDapFre2011}
C.~Bui, C.~Dapogny, and P.~Frey.
\newblock An accurate anisotropic adaptation method for solving the level set advection equation.
\newblock {\em Int. J. Numer. Meth. Fl.}, 70(7):899--922, 2011.

\bibitem{ChaDesMeh2013}
F.~Charles, B.~Després, and M.~Mehrenberger.
\newblock {Enhanced Convergence Estimates for Semi-Lagrangian Schemes Application to the Vlasov--Poisson Equation}.
\newblock {\em SIAM J. Numer. Anal.}, 51(2):840--863, 2013.

\bibitem{Chen2023ALC}
Y.~Chen, W.~Guo, and X.~Zhong.
\newblock A learned conservative semi-{L}agrangian finite volume scheme for transport simulations.
\newblock {\em J. Comput. Phys.}, 490:112329, 2023.

\bibitem{chen2024conservative}
Y.~Chen, W.~Guo, and X.~Zhong.
\newblock Conservative semi-{L}agrangian finite difference scheme for transport simulations using graph neural networks.
\newblock {\em J. Comput. Phys.}, 526:113768, 2025.

\bibitem{chen2024teng}
Z.~Chen, J.~McCarran, E.~Vizcaino, L.~Soljacic, and D.~Luo.
\newblock {TENG: Time-Evolving Natural Gradient for Solving PDEs With Deep Neural Nets Toward Machine Precision}.
\newblock In R.~Salakhutdinov, Z.~Kolter, K.~Heller, A.~Weller, N.~Oliver, J.~Scarlett, and F.~Berkenkamp, editors, {\em Proceedings of the 41st International Conference on Machine Learning}, volume 235 of {\em Proceedings of Machine Learning Research}, pages 7143--7162. PMLR, 21--27 Jul 2024.

\bibitem{CohMig2017}
A.~Cohen and G.~Migliorati.
\newblock Optimal weighted least-squares methods.
\newblock {\em SMAI J. Comput. Math.}, 3:181--203, 2017.

\bibitem{CroEinFao2015}
N.~Crouseilles, L.~Einkemmer, and E.~Faou.
\newblock Hamiltonian splitting for the {V}lasov-{M}axwell equations.
\newblock {\em J. Comput. Phys.}, 283:224--240, 2015.

\bibitem{CROUSEILLES20101927}
N.~Crouseilles, M.~Mehrenberger, and \'E. Sonnendrücker.
\newblock Conservative semi-{L}agrangian schemes for {V}lasov equations.
\newblock {\em J. Comput. Phys.}, 229(6):1927--1953, 2010.

\bibitem{CroMehVec2011}
N.~Crouseilles, M.~Mehrenberger, and F.~Vecil.
\newblock {Discontinuous Galerkin semi-Lagrangian method for Vlasov-Poisson}.
\newblock {\em ESAIM: Proceedings}, 32:211--230, 2011.

\bibitem{CroResSon2009}
N.~Crouseilles, T.~Respaud, and \'E. Sonnendrücker.
\newblock A forward semi-{L}agrangian method for the numerical solution of the {V}lasov equation.
\newblock {\em Comput. Phys. Commun.}, 180(10):1730--1745, 2009.

\bibitem{De-Ryck:2022aa}
T.~De~Ryck and S.~Mishra.
\newblock {Error analysis for physics-informed neural networks (PINNs) approximating Kolmogorov PDEs}.
\newblock {\em Adv. Comput. Math.}, 48(6), 2022.

\bibitem{DolHeiMisMos2023}
V.~Dolean, A.~Heinlein, S.~Mishra, and B.~Moseley.
\newblock Multilevel domain decomposition-based architectures for physics-informed neural networks.
\newblock {\em Comput. Method. Appl. M.}, 429:117116, 2024.

\bibitem{douglas1982numerical}
J.~Douglas, Jr and T.~F. Russell.
\newblock Numerical methods for convection-dominated diffusion problems based on combining the method of characteristics with finite element or finite difference procedures.
\newblock {\em SIAM J. Numer. Anal.}, 19(5):871--885, 1982.

\bibitem{EYu2018}
W.~E and B.~Yu.
\newblock {The Deep Ritz Method: A Deep Learning-Based Numerical Algorithm for Solving Variational Problems}.
\newblock {\em Commun. Math. Stat.}, 6(1):1--12, 2018.

\bibitem{EinJos2021}
L.~Einkemmer and I.~Joseph.
\newblock A mass, momentum, and energy conservative dynamical low-rank scheme for the vlasov equation.
\newblock {\em J. Comput. Phys.}, 443:110495, 2021.

\bibitem{EinKorKusMcCQiu2024}
L.~Einkemmer, K.~Kormann, J.~Kusch, R.~G. McClarren, and J.-M. Qiu.
\newblock A review of low-rank methods for time-dependent kinetic simulations.
\newblock {\em J. Comput. Phys.}, 538:114191, 2025.

\bibitem{Gra2016}
V.~Grandgirard et~al.
\newblock A 5{D} gyrokinetic full-$f$ global semi-{L}agrangian code for flux-driven ion turbulence simulations.
\newblock {\em Comput. Phys. Commun.}, 207:35--68, 2016.

\bibitem{falcone1998convergence}
M.~Falcone and R.~Ferretti.
\newblock {Convergence analysis for a class of high-order semi-Lagrangian advection schemes}.
\newblock {\em SIAM J. Numer. Anal.}, 35(3):909--940, 1998.

\bibitem{Fer2010}
R.~Ferretti.
\newblock {A Technique for High-Order Treatment of Diffusion Terms in Semi-Lagrangian Schemes }.
\newblock {\em Commun. Comput. Phys.}, 8(2):445--470, 2010.

\bibitem{ferretti2020stability}
R.~Ferretti and M.~Mehrenberger.
\newblock {Stability of semi-Lagrangian schemes of arbitrary odd degree under constant and variable advection speed}.
\newblock {\em Math. Comput.}, 89(324):1783--1805, 2020.

\bibitem{FILBET2001166}
F.~Filbet, \'E. Sonnendrücker, and P.~Bertrand.
\newblock {Conservative Numerical Schemes for the Vlasov Equation}.
\newblock {\em J. Comput. Phys.}, 172(1):166--187, 2001.

\bibitem{finzi2023stable}
M.~A. Finzi, A.~Potapczynski, M.~Choptuik, and A.~G. Wilson.
\newblock {A Stable and Scalable Method for Solving Initial Value PDEs with Neural Networks}.
\newblock In {\em The Eleventh International Conference on Learning Representations}, 2023.

\bibitem{GroKomLatSch2024}
T.~G. Grossmann, U.~J. Komorowska, J.~Latz, and C.-B. Schönlieb.
\newblock Can physics-informed neural networks beat the finite element method?
\newblock {\em IMA J. Appl. Math.}, 89(1):143--174, 2024.

\bibitem{GUO2025118101}
J.~Guo, G.~Domel, C.~Park, H.~Zhang, O.~Can~Gumus, Y.~Lu, G.~J. Wagner, D.~Qian, J.~Cao, T.~J.~R. Hughes, and W.~K. Liu.
\newblock {Tensor-decomposition-based A Priori Surrogate (TAPS) modeling for ultra large-scale simulations}.
\newblock {\em Comput. Methods Appl. Mech. Eng.}, 444, 2025.

\bibitem{guo2025interpolating}
J.~Guo, X.~Xie, C.~Park, H.~Zhang, M.~J. Politis, G.~Domel, and W.~K. Liu.
\newblock Interpolating neural network-tensor decomposition ({INN}-{TD}): a scalable and interpretable approach for large-scale physics-based problems.
\newblock In {\em Forty-second International Conference on Machine Learning}, 2025.

\bibitem{HamMehSelSon2016}
A.~Hamiaz, M.~Mehrenberger, H.~Sellama, and \'E. Sonnendrücker.
\newblock The semi-{L}agrangian method on curvilinear grids.
\newblock {\em Commun. Appl. Ind. Math.}, 7(3):99--137, 2016.

\bibitem{HeZhaRenSun2016}
K.~He, X.~Zhang, S.~Ren, and J.~Sun.
\newblock {Deep Residual Learning for Image Recognition}.
\newblock In {\em 2016 IEEE Conference on Computer Vision and Pattern Recognition (CVPR)}, pages 770--778. IEEE, 2016.

\bibitem{HirSmaDev2013}
M.~W. Hirsch, S.~Smale, and R.~L. Devaney.
\newblock {\em {Differential Equations, Dynamical Systems, and an Introduction to Chaos}}.
\newblock Elsevier, 2013.

\bibitem{HonSieTanXu2022}
Q.~Hong, J.~W. Siegel, Q.~Tan, and J.~Xu.
\newblock {On the Activation Function Dependence of the Spectral Bias of Neural Networks}.
\newblock {\em preprint}, 2022.

\bibitem{hu2024physics}
S.~Hu, M.~Liu, S.~Zhang, S.~Dong, and R.~Zheng.
\newblock Physics-informed neural network combined with characteristic-based split for solving {N}avier-{S}tokes equations.
\newblock {\em Eng. Appl. Artif. Intel.}, 128:107453, 2024.

\bibitem{Hu:2024aa}
Z.~Hu, K.~Shukla, G.~E. Karniadakis, and K.~Kawaguchi.
\newblock Tackling the curse of dimensionality with physics-informed neural networks.
\newblock {\em Neural Netw.}, 176:106369, 2024.

\bibitem{JinAPNN}
S.~Jin, Z.~Ma, and K.~Wu.
\newblock {Asymptotic-Preserving Neural Networks for Multiscale Time-Dependent Linear Transport Equations}.
\newblock {\em J. Sci. Comput.}, 94(3), 2023.

\bibitem{KAST2024112986}
M.~Kast and J.~S. Hesthaven.
\newblock Positional embeddings for solving {PDE}s with evolutional deep neural networks.
\newblock {\em J. Comput. Phys.}, 508:112986, 2024.

\bibitem{kolda2009tensor}
T.~G. Kolda and B.~W. Bader.
\newblock {Tensor Decompositions and Applications}.
\newblock {\em SIAM Rev.}, 51(3):455--500, 2009.

\bibitem{Kor2015}
K.~Kormann.
\newblock {A Semi-Lagrangian Vlasov Solver in Tensor Train Format}.
\newblock {\em SIAM J. Sci. Comput.}, 37(4):B613--B632, 2015.

\bibitem{KorReuRam2019}
K.~Kormann, K.~Reuter, and M.~Rampp.
\newblock A massively parallel semi-{L}agrangian solver for the six-dimensional {V}lasov-{P}oisson equation.
\newblock {\em Int. J. High Perform. C.}, 33(5):924--947, 2019.

\bibitem{NEURIPS2021_df438e52}
A.~Krishnapriyan, A.~Gholami, S.~Zhe, R.~Kirby, and M.~W. Mahoney.
\newblock Characterizing possible failure modes in physics-informed neural networks.
\newblock In {\em Thirty-fourth Conference on Neural Information Processing Systems}, 2021.

\bibitem{LagLikFot1998}
I.~E. Lagaris, A.~Likas, and D.~I. Fotiadis.
\newblock Artificial neural networks for solving ordinary and partial differential equations.
\newblock {\em IEEE Trans. Neural Netw.}, 9(5):987--1000, 1998.

\bibitem{LARIOSCARDENAS2022111623}
L.~Á. Larios-Cárdenas and F.~Gibou.
\newblock Error-correcting neural networks for semi-{L}agrangian advection in the level-set method.
\newblock {\em J. Comput. Phys.}, 471:111623, 2022.

\bibitem{LeV1996}
R.~J. LeVeque.
\newblock {High-Resolution Conservative Algorithms for Advection in Incompressible Flow}.
\newblock {\em SIAM J. Numer. Anal.}, 33(2):627--665, 1996.

\bibitem{LI2022110958}
L.~Li and C.~Yang.
\newblock {APFOS-Net: A}symptotic preserving scheme for anisotropic elliptic equations with deep neural network.
\newblock {\em J. Comput. Phys.}, 453:110958, 2022.

\bibitem{LiuCaiLapCao2021}
H.~Liu, X.~Cai, G.~Lapenta, and Y.~Cao.
\newblock {Conservative semi-Lagrangian kinetic scheme coupled with implicit finite element field solver for multidimensional Vlasov Maxwell system}.
\newblock {\em Commun. Nonlinear Sci.}, 102:105941, 2021.

\bibitem{LiuWan2025}
Z.~Liu and Y.~Wang et~al.
\newblock {KAN}: {K}olmogorov-{A}rnold {N}etworks.
\newblock In {\em The Thirteenth International Conference on Learning Representations}, 2025.

\bibitem{MaoMen2023}
Z.~Mao and X.~Meng.
\newblock Physics-informed neural networks with residual/gradient-based adaptive sampling methods for solving partial differential equations with sharp solutions.
\newblock {\em Appl. Math. Mech.}, 44(7):1069--1084, 2023.

\bibitem{MenLiZhaKar2020}
X.~Meng, Z.~Li, D.~Zhang, and G.~E. Karniadakis.
\newblock {PPINN: Parareal physics-informed neural network for time-dependent PDEs}.
\newblock {\em Comput. Method. Appl. M.}, 370:113250, 2020.

\bibitem{MISHRA2021107705}
S.~Mishra and R.~Molinaro.
\newblock Physics informed neural networks for simulating radiative transfer.
\newblock {\em J. Quant. Spectrosc. Radiat. Transfer}, 270:107705, 2021.

\bibitem{mojgani2022lagrangian}
R.~Mojgani, M.~Balajewicz, and P.~Hassanzadeh.
\newblock Kolmogorov $n$-width and {L}agrangian physics-informed neural networks: {A} causality-conforming manifold for convection-dominated {PDE}s.
\newblock {\em Comput. Method. Appl. M.}, 404:115810, 2023.

\bibitem{MueZei2023}
J.~M\"{u}ller and M.~Zeinhofer.
\newblock Achieving high accuracy with {PINN}s via energy natural gradient descent.
\newblock In A.~Krause, E.~Brunskill, K.~Cho, B.~Engelhardt, S.~Sabato, and J.~Scarlett, editors, {\em Proceedings of the 40th International Conference on Machine Learning}, volume 202 of {\em Proceedings of Machine Learning Research}, pages 25471--25485. PMLR, 23--29 Jul 2023.

\bibitem{NurLeiYan2023}
L.~Nurbekyan, W.~Lei, and Y.~Yang.
\newblock {Efficient Natural Gradient Descent Methods for Large-Scale PDE-Based Optimization Problems}.
\newblock {\em SIAM J. Sci. Comput.}, 45(4):A1621--A1655, 2023.

\bibitem{park2024interpolating}
C.~Park, S.~Saha, J.~Guo, H.~Zhang, X.~Xie, M.~A. Bessa, D.~Qian, W.~Chen, G.~J. Wagner, J.~Cao, and W.~K. Liu.
\newblock Interpolating neural network: {A} novel unification of machine learning and interpolation theory.
\newblock {\em arXiv preprint arXiv:2404.10296}, 2024.

\bibitem{pironneau1982transport}
O.~Pironneau.
\newblock {On the transport-diffusion algorithm and its applications to the Navier-Stokes equations}.
\newblock {\em Numer. Math.}, 38:309--332, 1982.

\bibitem{qiu2011positivity}
J.-M. Qiu and C.-W. Shu.
\newblock {Positivity preserving semi-Lagrangian discontinuous Galerkin formulation: theoretical analysis and application to the Vlasov--Poisson system}.
\newblock {\em J. Comput. Phys.}, 230(23):8386--8409, 2011.

\bibitem{QuaManNeg2016}
A.~Quarteroni, A.~Manzoni, and F.~Negri.
\newblock {\em {Reduced Basis Methods for Partial Differential Equations}}.
\newblock Springer International Publishing, 2016.

\bibitem{quarteroni2008numerical}
A.~Quarteroni and A.~Valli.
\newblock {\em Numerical approximation of partial differential equations}, volume~23.
\newblock Springer Science \& Business Media, 2008.

\bibitem{RAISSI2019686}
M.~Raissi, P.~Perdikaris, and G.~E. Karniadakis.
\newblock Physics-informed neural networks: {A} deep learning framework for solving forward and inverse problems involving nonlinear partial differential equations.
\newblock {\em J. Comput. Phys.}, 378:686--707, 2019.

\bibitem{restelli2006semi}
M.~Restelli, L.~Bonaventura, and R.~Sacco.
\newblock {A semi-Lagrangian discontinuous Galerkin method for scalar advection by incompressible flows}.
\newblock {\em J. Comput. Phys.}, 216(1):195--215, 2006.

\bibitem{rossmanith2011positivity}
J.~A. Rossmanith and D.~C. Seal.
\newblock {A positivity-preserving high-order semi-Lagrangian discontinuous Galerkin scheme for the Vlasov--Poisson equations}.
\newblock {\em J. Comput. Phys.}, 230(16):6203--6232, 2011.

\bibitem{SchFur2025}
N.~Schwencke and C.~Furtlehner.
\newblock {AN}a{GRAM}: A natural gradient relative to adapted model for efficient {PINN}s learning.
\newblock In {\em The Thirteenth International Conference on Learning Representations}, 2025.

\bibitem{SONNENDRUCKER1999201}
\'E. Sonnendrücker, J.~Roche, P.~Bertrand, and A.~Ghizzo.
\newblock {The Semi-Lagrangian Method for the Numerical Resolution of the Vlasov Equation}.
\newblock {\em J. Comput. Phys.}, 149(2):201--220, 1999.

\bibitem{SemiLagrangianIntegrationSchemesforAtmosphericModelsAReview}
A.~Staniforth and J.~Côté.
\newblock {Semi-Lagrangian Integration Schemes for Atmospheric Models---A Review}.
\newblock {\em Mon. Weather Rev.}, 119(9):2206--2223, 1991.

\bibitem{StiCha2023}
J.~Stiasny and S.~Chatzivasileiadis.
\newblock Physics-informed neural networks for time-domain simulations: {A}ccuracy, computational cost, and flexibility.
\newblock {\em Electr. Pow. Syst. Res.}, 224:109748, 2023.

\bibitem{SUKUMAR2022114333}
N.~Sukumar and A.~Srivastava.
\newblock Exact imposition of boundary conditions with distance functions in physics-informed deep neural networks.
\newblock {\em Comput. Method. Appl. M.}, 389:114333, 2022.

\bibitem{SarSol2019}
S.~Särkkä and A.~Solin.
\newblock {\em {Applied Stochastic Differential Equations}}.
\newblock Cambridge University Press, 2019.

\bibitem{TanSriMilFriRagSinRamBarNg2020}
M.~Tancik, P.~Srinivasan, B.~Mildenhall, S.~Fridovich-Keil, N.~Raghavan, U.~Singhal, R.~Ramamoorthi, J.~Barron, and R.~Ng.
\newblock {Fourier Features Let Networks Learn High Frequency Functions in Low Dimensional Domains}.
\newblock In H.~Larochelle, M.~Ranzato, R.~Hadsell, M.~F. Balcan, and H.~Lin, editors, {\em Advances in Neural Information Processing Systems}, volume~33, pages 7537--7547. Curran Associates, Inc., 2020.

\bibitem{WanLiChePer2024}
S.~Wang, B.~Li, Y.~Chen, and P.~Perdikaris.
\newblock {PirateNets: Physics-informed Deep Learning with Residual Adaptive Networks}.
\newblock {\em J. Mach. Learn. Res.}, 25:1--51, 2024.

\bibitem{WANG2024116813}
S.~Wang, S.~Sankaran, and P.~Perdikaris.
\newblock Respecting causality for training physics-informed neural networks.
\newblock {\em Comput. Method. Appl. M.}, 421:116813, 2024.

\bibitem{WuZhuTanKarLu2023}
C.~Wu, M.~Zhu, Q.~Tan, Y.~Kartha, and L.~Lu.
\newblock A comprehensive study of non-adaptive and residual-based adaptive sampling for physics-informed neural networks.
\newblock {\em Comput. Methods Appl. Mech. Engrg.}, 403:115671, 2023.

\bibitem{XIU2001658}
D.~Xiu and G.~E. Karniadakis.
\newblock {A Semi-Lagrangian High-Order Method for Navier-Stokes Equations}.
\newblock {\em J. Comput. Phys.}, 172(2):658--684, 2001.

\bibitem{YANG2021110632}
C.~Yang and M.~Mehrenberger.
\newblock {Highly accurate monotonicity-preserving Semi-Lagrangian scheme for Vlasov-Poisson simulations}.
\newblock {\em J. Comput. Phys.}, 446:110632, 2021.

\bibitem{Zhang_2023}
B.~Zhang, G.~Cai, H.~Weng, W.~Wang, L.~Liu, and B.~He.
\newblock Physics-informed neural networks for solving forward and inverse {V}lasov-{P}oisson equation via fully kinetic simulation.
\newblock {\em Mach. Learn.: Sci. Technol.}, 4(4):045015, 2023.

\bibitem{ZheHayChrQiu2025}
N.~Zheng, D.~Hayes, A.~Christlieb, and J.-M. Qiu.
\newblock {A Semi-Lagrangian Adaptive-Rank (SLAR) Method for Linear Advection and Nonlinear Vlasov-Poisson System}.
\newblock {\em J. Comput. Phys.}, page 113970, 2025.

\end{thebibliography}
\addcontentsline{toc}{section}{References}

\appendix

\section{Algorithm to compute the approximate characteristic curve}
\label{sec:alg_characteristic_curve}

In this section, we explain how to actually compute the
approximate characteristic curve
\smash{$\widetilde{\mathcal{X}}$},
in \cref{alg:characteristic_curve}.
In this algorithm, a time-stepping algorithm must be chosen
to approximate the characteristic curve.
We elect to use the RK4 (fourth-order Runge-Kutta) method,
a well-known and widely used ODE solver.
For some time step $\delta \tau$,
some initial condition $x$ at time $\tau$,
and some parameters $\mu$,
we denote by
$\text{RK4}(\tau, \delta \tau, x, \mu)$
the RK4 approximation of the solution
to the ODE~\eqref{eq:characteristic_curve}
at time $\tau - \delta \tau$.

\begin{algorithm}[!ht]
    \caption{Computation of the characteristic curve for advection equations}
    \label{alg:characteristic_curve}
    \begin{algorithmic}[1]
        \State \textbf{Input:} Final position $x$, parameters $\mu$, final time $t^{n+1}$, initial time $t^n = t^{n+1} - \dt$, advection field $a$
        \State \textbf{Output:} Approximate foot \smash{$\widetilde{\mathcal{X}}(t^n; t^{n+1}, x, \mu)$} of the characteristic curve

        \medskip

        \If{$a(\mu)$ is constant in space}

        \State Set $\widetilde{\mathcal{X}}(t^n; t^{n+1}, x, \mu) = x - a(\mu) \dt$

        \ElsIf{$a(x, \mu)$ is such that the ODE \eqref{eq:characteristic_curve} has a closed-form solution $\mathcal{X}$}

        \State Set $\widetilde{\mathcal{X}}(t^n; t^{n+1}, x, \mu) = \mathcal{X}(t^n; t^{n+1}, x, \mu)$

        \Else

        \State Numerically solve the ODE \eqref{eq:characteristic_curve} with initial condition $x$ at time $t^{n+1}$

        \State \textbf{Initialization:} set $\widetilde{\mathcal{X}} = x$, set $\delta \tau$ and $n_\tau$ such that $n_\tau \delta \tau = \dt$, set $\tau = t^{n+1}$

        \While{$\tau > t^n$}

        \State Solve \eqref{eq:characteristic_curve}
        using the RK4 method on the
        sub-time step $\delta \tau$:
        $\widetilde{\mathcal{X}} \leftarrow \text{RK4}(\tau, \delta \tau, \widetilde{\mathcal{X}}, \mu)$

        \State Update the sub-time step: $\tau \leftarrow \tau - \delta \tau$

        \EndWhile

        \EndIf

        \medskip

        \State \textbf{Return:} $\widetilde{\mathcal{X}}(t^n; t^{n+1}, x, \mu)$
    \end{algorithmic}
\end{algorithm}

\section{Hyperparameters for the numerical experiments}
\label{sec:hyperparameters}

This section regroups all the hyperparameters used in the numerical experiments.
We denote by:
\begin{itemize}[nosep]
    \item $N_c$ the number of collocation points, i.e.,
          the number of points uniformly sampled to discretize the space domain $\Omega$
          and approximate the integrals in the optimization problems,
          see \cref{alg:NSL_known_characteristic};
    \item $N_e$ the number of epochs;
    \item $\Xi$ the activation function;
    \item $\omega^\text{BC}$ (resp. $\omega^\text{IC}$) the
          boundary (resp. initial) loss function coefficient,
          i.e., the coefficient by which the boundary (resp. initial) loss function is multiplied
          when weakly imposing the boundary (resp. initial) conditions in the PINN;
    \item $N_c^\text{BC}$ (resp. $N_c^\text{IC}$) the number of
          collocation points in the boundary (resp. initial) integrals in the PINN;
    \item $\ell$ the list of layer sizes,
    \item $\lambda$ the learning rate (note that no learning rate is reported when natural gradient preconditioning is deployed\rc{, since the default one of $1 \times 10^{-3}$ is used}).
\end{itemize}

Note that $\omega^\text{BC}$ or $\omega^\text{IC}$
being set to $0$ means that the corresponding loss function is not used,
and that the boundary or initial conditions are strongly imposed.
Moreover, the activation function $\hat h$ corresponds
to the following regularized hat function,
see~\cite{HonSieTanXu2022}:
\begin{equation}
    \label{eq:regularized_hat_function}
    \hat{h}: x \mapsto \exp\left(-12 \tanh\left(\frac{x^2}{2}\right)\right).
\end{equation}

For each experiment, we report the hyperparameters used:
\begin{itemize}[nosep]
    \item to train the PINN;
    \item to compute the initial dofs $\theta^0$ for the dPINN, NG and NSL methods (in the column called ``initialization'');
    \item to iterate the dPINN method, i.e.,
          to solve the nonlinear optimization
          problem~\eqref{eq:dPINN_optimization_problem} at each time step;
    \item to iterate the NSL method, i.e.,
          to solve the nonlinear optimization
          problem~\eqref{eq:NGSL_nonlinear_optimization} at each time step.
\end{itemize}

\rc{While we did not run a
    full ablation or sensitivity analysis,
    we nevertheless provide some insights into
    the choice of hyperparameters
    for the NSL method.
    For both the initialization and the iterations,
    standard neural network training strategies are used
    to affect the different error terms
    $\varepsilon_\text{int}$,
    $\varepsilon_\text{opt}$ and
    $\varepsilon_\text{approx}$.
    More collocation points will lower $\varepsilon_\text{int}$,
    while more epochs will lower $\varepsilon_\text{opt}$
    and additional layers will lower $\varepsilon_\text{approx}$.
    \begin{itemize}[nosep]
        \item Correctly training the network to approximate the initial condition
              is the most important step for the NSL method.
              If the initial condition is not well-approximated,
              then the time iterations will also not provide a good approximation.
              This is why a larger number $N_e$ of epochs is often used
              for the initialization than for the iterations.
        \item The number of collocation points $N_c$ is usually similar
              for the initialization and the iterations,
              except when the solution becomes more complex as time progresses.
              In this case, $N_c$ may be increased for the iterations.
        \item Especially when natural gradient preconditioning is used,
              the number of layers can be kept small,
              as even small networks are very expressive if they are well-trained.
        \item Finally, the activation function $\Xi$ is usually chosen
              to be a $\tanh$ function,
              except for oscillatory solutions where a $\sin$ function is used,
              or for localized solutions benefiting from
              a regularized hat function $\hat{h}$.
    \end{itemize}
}

\begin{table}[!ht]
    \centering
    \begin{tabular}{ccccc}
        \toprule
        Hyperparameter     & PINN           & initialization                                & dPINN iterations & NSL iterations \\
        \cmidrule(lr){1-5}
        $\lambda$          & $\num{9e-3}$   & $\num{1.5e-2}$                                & $\num{1.5e-2}$   & $\num{1.5e-2}$ \\
        $N_e$              & $\num{1000}$   & $\num{1500}$                                  & $200$            & $150$          \\
        $N_c$              & $\num{2000}$   & $\num{3000}$ $(\num{5000} \text{ for dPINN})$ & $\num{5000}$     & $\num{1000}$   \\
        $\ell$             & $[40, 40, 40]$ & $[40, 40, 40]$                                & N/A              & N/A            \\
        $\Xi$              & $\tanh$        & $\tanh$                                       & N/A              & N/A            \\
        $\omega^\text{BC}$ & $50$           & N/A                                           & N/A              & N/A            \\
        $N_c^\text{BC}$    & $\num{1000}$   & N/A                                           & N/A              & N/A            \\
        $\omega^\text{IC}$ & $0$            & N/A                                           & N/A              & N/A            \\
        $N_c^\text{IC}$    & $\num{0}$      & N/A                                           & N/A              & N/A            \\
        \bottomrule
    \end{tabular}
    \caption{Hyperparameters used in \cref{sec:1D_transport_constant_non_parametric}.}
    \label{tab:hyperparameters_1D_transport_constant_non_parametric}
\end{table}

\begin{table}[!ht]
    \centering
    \begin{tabular}{ccccc}
        \toprule
        Hyperparameter     & PINN                   & initialization                                 & dPINN iterations & NSL iterations \\
        \cmidrule(lr){1-5}
        $\lambda$          & $\num{9e-3}$           & $\num{1.5e-2}$                                 & $\num{5e-3}$     & $\num{5e-3}$   \\
        $N_e$              & $\num{3000}$           & $\num{2500}$                                   & $\num{2500}$     & $\num{2500}$   \\
        $N_c$              & $\num{10000}$          & $\num{3000}$ $(\num{10000} \text{ for dPINN})$ & $\num{1000}$     & $\num{10000}$  \\
        $\ell$             & $[40, 60, 80, 60, 40]$ & $[40, 60, 80, 60, 40]$                         & N/A              & N/A            \\
        $\Xi$              & $\tanh$                & $\tanh$                                        & N/A              & N/A            \\
        $\omega^\text{BC}$ & $50$                   & N/A                                            & N/A              & N/A            \\
        $N_c^\text{BC}$    & $\num{1000}$           & N/A                                            & N/A              & N/A            \\
        $\omega^\text{IC}$ & $0$                    & N/A                                            & N/A              & N/A            \\
        $N_c^\text{IC}$    & $\num{0}$              & N/A                                            & N/A              & N/A            \\
        \bottomrule
    \end{tabular}
    \caption{Hyperparameters used in \cref{sec:1D_transport_constant_parametric}.}
    \label{tab:hyperparameters_1D_transport_constant_parametric}
\end{table}

\begin{table}[!ht]
    \centering
    \begin{tabular}{cccc}
        \toprule
        Hyperparameter     & PINN               & NG and NSL initialization & NSL iterations \\
        \cmidrule(lr){1-4}
        $\lambda$          & $\num{1e-2}$       & $\num{1e-3}$              & $\num{5e-3}$   \\
        $N_e$              & $\num{5000}$       & $\num{700}$               & $\num{500}$    \\
        $N_c$              & $\num{25000}$      & $\num{15000}$             & $\num{15000}$  \\
        $\ell$             & $[80, 80, 80, 80]$ & $[40, 40, 40, 40, 40]$    & N/A            \\
        $\Xi$              & $\hat h$           & $\hat h$                  & N/A            \\
        $\omega^\text{BC}$ & $0$                & N/A                       & N/A            \\
        $N_c^\text{BC}$    & $\num{0}$          & N/A                       & N/A            \\
        $\omega^\text{IC}$ & $50$               & N/A                       & N/A            \\
        $N_c^\text{IC}$    & $\num{15000}$      & N/A                       & N/A            \\
        \bottomrule
    \end{tabular}
    \caption{Hyperparameters used in \cref{sec:2D_transport_rotating_parametric}.}
    \label{tab:hyperparameters_2D_transport_rotating_parametric}
\end{table}

\begin{table}[!ht]
    \centering
    \begin{tabular}{cccc}
        \toprule
        Hyperparameter     & PINN               & NG and NSL initialization & NSL iterations \\
        \cmidrule(lr){1-4}
        $\lambda$          & $\num{1e-2}$       & $\num{1e-3}$              & $\num{5e-3}$   \\
        $N_e$              & $\num{1000}$       & $\num{250}$               & $\num{500}$    \\
        $N_c$              & $\num{15000}$      & $\num{5000}$              & $\num{15000}$  \\
        $\ell$             & $[80, 80, 80, 80]$ & $[40, 40, 40, 40, 40]$    & N/A            \\
        $\Xi$              & $\hat h$           & $\hat h$                  & N/A            \\
        $\omega^\text{BC}$ & $50$               & N/A                       & N/A            \\
        $N_c^\text{BC}$    & $\num{5000}$       & N/A                       & N/A            \\
        $\omega^\text{IC}$ & $50$               & N/A                       & N/A            \\
        $N_c^\text{IC}$    & $\num{1000}$       & N/A                       & N/A            \\
        \bottomrule
    \end{tabular}
    \caption{Hyperparameters used in \cref{sec:2D_fake_Vlasov}.}
    \label{tab:hyperparameters_2D_fake_Vlasov}
\end{table}

\begin{table}[!ht]
    \centering
    \rallparagraph
    \begin{tabular}{cccc}
        \toprule
        Hyperparameter & \cref{sec:2D_level_set} & \cref{sec:3D_level_set} & \cref{sec:3D_level_set_param} \\
        \cmidrule(lr){1-4}
        $N_e$          & $\num{250}$             & $\num{250}$             & $\num{500}$                   \\
        $N_c$          & $\num{60000}$           & $\num{64}^3$            & $\num{48}^3$                  \\
        $\ell$         & $[35, 50, 35]$          & $[35, 50, 50,35]$       & see legend                    \\
        $\Xi$          & $\tanh$                 & $\tanh$                 & $\tanh$                       \\
        \bottomrule
    \end{tabular}
    \caption{\rallparagraph Hyperparameters used in \cref{sec:2D_level_set,sec:3D_level_set,sec:3D_level_set_param}. For the 3D level-set deformation without and with parameters (\cref{sec:3D_level_set,sec:3D_level_set_param}), the number of epochs in the first time step is increased to $500$ and $1000$ respectively. For the 3D level-set deformation with parameters, a residual network (see~\cite{HeZhaRenSun2016}) is employed, with $9$ layers of $26$ neurons each, with skip connections between layers $1$ and $3$, $4$ and $6$, and $7$ and $9$.}
    \label{tab:hyperparameters_level_set}
\end{table}

\begin{table}[!ht]
    \centering
    \rallparagraph
    \begin{tabular}{ccc}
        \toprule
        Hyperparameter     & Network for $u$ & Network for $\Psi$ \\
        \cmidrule(lr){1-3}
        $N_e^\text{init.}$ & $\num{500}$     & $\num{350}$        \\
        $N_e^\text{iter.}$ & $\num{100}$     & $\num{50}$         \\
        $N_c$              & $\num{112}^2$   & $\num{768}$        \\
        $\ell$             & $[30] \times 9$ & $[20] \times 6$    \\
        $\Xi$              & $\sin$          & $\sin$             \\
        \bottomrule
    \end{tabular}
    \caption{\rallparagraph Hyperparameters used in \cref{sec:vlasov_poisson_validation}. For both residual networks, skip connections skip one layer (they run from layers $1$ to $3$ and $4$ to $6$ for both $u$ and $\Psi$, and $7$ to $9$ for $u$).}
    \label{tab:hyperparameters_bump_on_tail}
\end{table}

\end{document}